\documentclass[11pt,reqno,a4paper]{amsart}
\usepackage{amssymb}
\usepackage[cp1251]{inputenc}
\usepackage[russian]{babel}
\usepackage{graphicx}
\newcommand{\reff}[1]{\mbox{\rm (\ref{#1})}}
\newcommand{\om}{\omega}
\newcommand{\omp}{{\omega{\scriptstyle'}}}
\newcommand{\etap}{{\eta{\scriptstyle'}}}
\newcommand{\wpp}{{\wp{\scriptstyle'}}}

\newcommand{\ri}{\mathrm{i}}
\newcommand{\re}{\mathrm{e}}
\newcommand{\ds}{\displaystyle}
\newcommand{\sss}{\scriptscriptstyle}
\newcommand{\ded}{{\boldsymbol\eta}}

\newcommand{\elambda}{{e_{\sss\!\lambda}^{}}}
\newcommand{\emu}{{e_{\sss\!\mu}}}
\newcommand{\sigmalambda}{{\sigma_{\sss\!\lambda}^{}}}
\newcommand{\Kappa}{\boldsymbol{\symbol{26}}}

\newcommand{\mfrac}[2] 
{\raisebox{0.045em}{\mbox{\footnotesize$\displaystyle
\frac{#1}{#2}$}}}

\newcommand{\spin}[1]{\mbox{\footnotesize$\langle#1\rangle$}}

\newcommand{\thetaAB}[2]
{\theta \raisebox{0.04em}{\mbox{$ \scalebox{0.5}[0.9]{\big[}
\!\!\raisebox{0.1em}{\scalebox{0.8}{\mbox{\tiny$
\begin{array}{c}#1\\#2\end{array}$}}}
\!\! \scalebox{0.5}[0.9]{\big]}$}}}

\newcommand{\varthetaAB}[2]
{\vartheta \raisebox{0.04em}{\mbox{$ \scalebox{0.5}[0.9]{\big[}
\!\!\raisebox{0.1em}{\scalebox{0.8}{\mbox{\tiny$
\begin{array}{c}#1\\#2\end{array}$}}}
\!\! \scalebox{0.5}[0.9]{\big]}$}}}

\newcommand{\bbig}[1]
{\raisebox{-0.1em}{\scalebox{1.1}[1.5]{\mbox{$#1$}}}}

\begin{document}
\hfill {\tt http://arXiv.org/math.CA/0808.3486}

\vskip 1cm \centerline{\Large\textbf{\textsc{О функциях Якоби и
Вейерштрасса (I)}}} \centerline{\Large\textbf{\textsc{и уравнение
Пенлеве}}}

\vskip 0.5cm
\centerline{\large Ю.\,В. Брежнев}
\vskip 0.5cm

\centerline{{\tt   brezhnev@mail.ru}} \thispagestyle{empty}

\vskip 1.5cm \centerline{ {\hsize =11 true cm \vbox{ \small
\noindent Статья является переработанной и существенно расширенной
версией работы
{\tt math.CA/0601371}, дополненной приложением.\\
\phantom{aaa}Предлагаются новые результаты в теории классических
$\theta$-функций Якоби и $\sigma$-функций Вейерштрасса: обыкновенные
дифференциальные уравнения  и разложения в ряды. Мы даем также
расширение канонических $\theta$-функций и рассматриваем приложение
к шестому уравнению Пенлеве (P6). Общее решение Пикара--Хитчина
уравнения P6 явно представлено в виде логарифмической производной от
соответствующей $\tau$-функции (форма Пенлеве).}}}

\vskip 1.5cm \centerline{ {\hsize =10 true cm \vbox{ \small
\noindent The paper is an essentially extended version of the work
{\tt
math.CA/0601371}, supplemented with an application.\\
\phantom{aaa} We present new results in the theory of classical
$\theta$-functions of Jacobi and $\sigma$-functions of Weierstrass:
ordinary differential equations and series expansions. We also give
the extension of canonical $\theta$-functions and consider an
application to the sixth Painlev\'e equation (P6). Picard--Hitchin's
general solution of P6 is represented explicitly in a form of
logarithmic derivative of a corresponding $\tau$-function
(Painlev\'e's form).}}}


\pagebreak

\hfill \parbox{85mm}{\em \footnotesize \ldots\ мне как читателю
оказывается сложным разобраться в том, что же является новым, а что
заимствовано из классических источников. $(\boldsymbol{!})$

\ldots\ в этой работе все ``новые'' утверждения являются более или
менее тривиальными следствиями известных фактов.

\ldots\ многие формулы из \S 3,\;5--8 рецензенту хорошо известны.
$(\boldsymbol{!!})$

\ldots\ на поверку оказывается тривиальным переписыванием хорошо
известного параметрического решения уравнения Пенлеве.

\ldots\  работа произвела впечатление некоего справочного пособия,
адресованного математикам не знакомым с классической теорией.
$(\boldsymbol{!})$

\ldots\ referee did not see these formulae before but believe that
such formulae can be derived. $(\boldsymbol{!!!})$

\ldots\ the paper can't be considered as a research paper.

\ldots\ it seems that he has not made a really systematic attempt to
check all the most likely sources. If he had done so, the number of
``useful formula'' that could be gleaned from them could probably
easily be multiplied by a significant factor. $(\boldsymbol{!!})$\\ \\

\em (Выдержки из отзывов на работу. Знаком $(\boldsymbol{!})$
отмечены наиболее <<замечательные образцы>>. В качестве
дополнительного введения, после статьи приведено письмо,
характеризующее  уровень компетентности  рецензий.)}

\section{Введение}

\noindent Тэта-функции Якоби   и вейерштрассовский базис функций
$(\sigma,\,\zeta,\,\wp,\,\wpp)$ возникают в  многочисленных теориях
и приложениях. Они были предметом интенсивных исследований начиная с
момента их возникновения и  уже к концу XIX века
теория приняла свой нынешний  вид. Подавляющее число результатов
было получено в работах самого Якоби и Вейерштрасса, а также их
современников: Эрмит, Эннепер, Киперт, Нейман, Альфен, Гурвиц,
Фробениус, Фрике и др. Существует обширная литература, включая
справочную \cite{we2,bateman}, по теории эллиптических, модулярных и
родственных им функций. Большое количество показательных примеров
приведено в классическом {\em ``Курсе современного анализа''\/}
\cite{WW},  а наиболее обстоятельными, как по изложению теории, так
и в формульно-справочном отношении, являются не  часто упоминаемые
четыре тома Таннери--Молька \cite{tannery}, а также учебник по
алгебре Вебера \cite{weber}, 2-х томный труд Альфена с посмертным
изданием 3-го тома \cite{halphen} и др.
\cite{koenig,krause,bateman}. Сюда  же следует отнести и собственно
труды Вейерштрасса \cite{we} и Якоби \cite{jacobi}, которые, являясь
весьма подробными в изложениях, остаются источником важных
наблюдений.

В настоящей работе ({\bf I}) мы приводим  новые результаты
применительно к $\theta$-функ\-ци\-ям Якоби и $\sigma$-функциям
Вейерштрасса, а именно, по разложению их в ряды и дифференциальным
уравнениям, которым они удовлетворяют. В отдельном сообщении ({\bf
II}), как продолжение, мы изложим способ построения <<эллиптической
части>> теории, не связанный  с обращением голоморфного  интеграла
(Якоби--Вейерштрасс) или с подходом, основанным на построении
эллиптических функций как мероморфных и двояко-периодических
(Лиувилль--Эйзенштейн--Вейерштрасс).

Разложения в ряды эллиптических/модулярных и $\theta$-функций до сих
пор являются постоянно эксплуатируемым инструментом во многих
вопросах. В силу известных аналитических свойств, это позволяет
получать точные результаты. Достаточно упомянуть функциональные ряды
для всевозможных $\theta$-отношений \cite{weber,tannery},
тео\-ре\-ти\-ко-числовые $q$-ряды и ряды Ламберта \cite{apostol},
знаменитые ряды МакКая--Томпсона, их следствия типа <<Moonshine
Conjecture>> и ее современные расширения. Многочисленные формы рядов
для вейерштрассовских и якобиевских функций фигурируют в литературе
повсюду \cite{halphen,krause,tannery,bateman,WW}. Несколько
неожиданным выглядит то, что до сих пор отсутствуют формулы
степенных $\theta$-рядов   как аналогов вейерштрассовских
$\sigma$-рядов, хотя Якоби пытался получать их еще до появления его
{\em Fun\-da\-men\-ta Nova\/} \cite[{\bf I}:\,стр.\,259--60]{jacobi}
(см. далее \S\,4).

Предлагаемые ниже динамические системы имеют самостоятельную
ценность, поскольку, как мы увидим, в качестве базовых
дифференциальных уравнений логичнее рассматривать не столько
уравнения на эллиптические функции ($\wp$-уравнение Вейерштрасса или
уравнения Эйлера на якобиевские sn, cn, dn), сколько {\em
обыкновенные дифференциальные уравнения на собственно
$\theta$-функции\/} (\S\,6--7). При этом выстраивается  естественная
картина: 1) эллиптические функции являются подклассом абелевых
эллиптических интегралов всех трех родов, они
--- подклассом канонических $\theta$-функций, а последние ---
подклассом описываемых  (\S\,9) неканонических расширений\footnote{К
началу этой цепочки можно было бы добавить элементарные функции. Под
классом здесь условно понимаются рациональные функции и  логарифмы
от них.}; 2) дифференциально-модулярные свойства всех классов
возникают автоматически (\S\,8--9); 3) алгебраические свойства
($\theta$-тождества) появляются как алгебраические интегралы
дифференциалов (\S\,9.1). Обширные приложения эта тематика имеет,
как известно, в теории интегрируемых нелинейных уравнений.

Дифференциальные свойства <<модулярной части>> функций Якоби более
трансцендентны  и неразрывно связаны с теорией линейных
дифференциальных уравнений с бесконечными  группами фуксовых
монодромий.  За последнее время, в работах Абловица
\cite[стр.\,573--89]{conte}, Тахтаджана, Харнада, МакКая, Хитчина
\cite{hitchin}, Охиямы и других, эта область нашла развитие и
красивые приложения, известные как монополи \cite{hitchin2} и
уравнения типа Альфена \cite{tod,babich},  Шази--Пикара--Фукса
\cite{tah}, космологические метрики Тода \cite{tod}, Хитчина
\cite{hitchin} и др.

Связь между интегрируемыми нелинейными уравнениями и уравнениями
интегрируемыми в модулярных функциях  замечена давно и не раз
отмечалась в литературе. См. например работы \cite{tah},
\cite[стр.\,573--89]{conte} или обсуждение соотношения с уравнениями
типа КдВ на стр.\,126--8 в книге \cite{hitchin2}. Мы дадим
объяснение этому факту и увидим, что  <<модулярная интегрируемость>>
имеет следующий характер. Она является условием совместности
линейного <<$\tau$-уравнения>> теплопроводности и
квадратурно-алгебраически интегрируемых нелинейных <<$x$-уравнений>>
на $\theta$-функции.

Другим и не менее важным объектом приложений излагаемого ниже
является проблема  координатно-аналитического описания пространств
модулей алгебраических кривых и
диф\-фе\-рен\-ци\-аль\-но-гео\-мет\-ри\-чес\-ких структур на них.
Этот вопрос в литературе  не поднимался, в то время как ясно, что
предлагаемая <<дифференциальная техника>> будет необходимой в тех
ситуациях, когда классы алгебраических кривых накрывают
эллиптические. Этим был вызван наш основной интерес в этой области,
но мы не  затрагиваем эту тему, ограничившись только информативно
справочной частью аппарата и указанием работы \cite{br}, где
предлагаемые результаты   используются в нескольких приложениях.
Другие приложения  следуют из самих формул.

Начиная с \S\,3 и до конца работы мы будем придерживаться следующего
правила. Если отдельно не указывается или про формулу ничего не
сказано, то результат (утверждение в тексте) является новым. В
некоторых случаях, несмотря на ожидаемость ответа, следствия
оказываются не очевидными. В частности, это относится к
дифференциальным уравнениям на якобиевские функции по первому
аргументу. Дифференциальные тождества между $\theta$-функциями
иногда встречаются в старой литературе \cite{weber,tannery}, но
важна дифференциально-полиномиальная {\em замкнутость\/}, а наряду с
функциями $\theta$ в теории должна фигурировать функция $\theta_1'$.
Это сводит б\'ольшую часть вычислений к упражнениям в <<differential
calculus>> и позволяет избежать манипуляций с неявными
алгебраическими соотношениями. Способы вывода формул обычно следуют
из контекста и поэтому, за недостатком места, мы опускаем
доказательства, ограничиваясь комментариями к существенным деталям.
Содержание \S\,2 является общеизвестным  и приводится здесь для
фиксации обозначений. Список литературы, чтобы избежать его
чрезмерного увеличения,  во многом сокращен.

\section{Определения и обозначения}

\subsection{Функции Якоби}
Четыре функции $\theta_{1,2,3,4}$  и их эквиваленты при
использовании записи с характеристиками  определяются следующими
каноническими рядами:
$$
-\theta_{11}:\qquad \theta_1(x|\tau)= -\ri\; \sum\limits_{\sss
k=-\infty}^{\sss \infty}\! (-1)^k\,
\re^{(k+\frac12)^2\pi\ri\,\tau}\, \re^{(2k+1)\pi \ri\, x}\,,\!
$$
$$
\theta_{10}:\qquad \theta_2(x|\tau)= \phantom{-\ri\;}
\sum\limits_{\sss k=-\infty}^{\sss \infty}
\re^{(k+\frac12)^2\pi\ri\,\tau}\,\re^{(2k+1)\pi \ri\, x}\,, \quad\;
$$
$$
\theta_{00}:\qquad \theta_3(x|\tau)= \phantom{-\ri\;}
\sum\limits_{\sss k=-\infty}^{\sss \infty}
\re^{k^2\pi\ri\,\tau}\,\re^{2k\pi \ri\, x}\,, \qquad\qquad\quad\!\!
$$
$$
\theta_{01}:\qquad \theta_4(x|\tau)=
\phantom{-\ri\;}\sum\limits_{\sss k=-\infty}^{\sss \infty}\!
(-1)^k\, \re^{k^2\pi\ri\,\tau}\, \re^{2k\pi \ri\, x}\,. \qquad\;\,\!
$$
Будем использовать  обозначения $\theta_k\equiv \theta_k(x|\tau)$ и
иногда $\theta_k(x)$. Значения $\theta$-функций при $x=0$ называются
$\vartheta$-константами:
$\vartheta_k\equiv\vartheta_k(\tau)=\theta_k(0|\tau)$. По
соображениям удобства форматирования формул,  будем употреблять
двоякую нотацию для $\theta$-функций с характеристиками (Эрмит):
$$
\thetaAB{\alpha}{\beta}(x|\tau)= \theta_{\alpha\beta}(x|\tau)=
{\sum\limits_{{}^{\sss k=-\infty}}^{\sss \infty}}^ {\ds\mathstrut}
\re^{\pi\ri \left(\!k+\frac\alpha2\!\right)^{\!2}\tau+
2\pi\ri\left(\!k+\frac\alpha2\!\right)\!
\left(\!x+\frac\beta2\!\right)\! }_{\mathstrut}\,.
$$
Пусть $(n,m)=0,\pm1,\pm2,\ldots\,$. Мы рассматриваем только
целочисленные характеристики, поэтому, в силу формулы
$\thetaAB{\alpha+2m}{\beta+2n}=(-1)^{\alpha
n}\!\cdot\!\thetaAB{\alpha}{\beta}$, функции $\theta_{\alpha\beta}$
всегда сводятся к $\pm\theta_{1,2,3,4}$. При сдвиге на
$\frac12$-периоды у $\theta$-функции сдвигаются характеристики:
\begin{equation}\label{shift}
\thetaAB{\alpha}{\beta} \!\left({\textstyle x+\frac n2+\frac
m2\,\tau}\big|\tau\right)=(-\ri)^{(\beta+n)m}\:
\thetaAB{\alpha\mbox{\tiny+}m}{\beta\mbox{\tiny+}n}(x|\tau)
\!\cdot\!\re^{\mbox{\tiny--}\frac{\pi\ri}{4}
\,m(4x+m\tau)}_{\mathstrut}\,.
\end{equation}
Двукратные сдвиги на $\frac12$-периоды дают закон преобразования
$\theta$-функции  в себя:
$$
\theta_{\alpha\beta} (x+n+m\,\tau|\tau)=
(-1)^{n\alpha-m\beta}\,\theta_{\alpha\beta} (x|\tau)\!\cdot\!
\re^{\mbox{\tiny--}\pi\ri\,m(2x+m\tau)}\,.
$$
Значение любой $\theta$-функции в любом $\frac12$-периоде есть
некоторая $\vartheta$-константа с экспоненциальным множителем:
$$
\thetaAB{\alpha}{\beta} \!\left({\textstyle\frac n2+\frac
m2\,\tau}\big|\tau\right)=(-\ri)^{(\beta+n)m}\:
\varthetaAB{\alpha\mbox{\tiny+}m}{\beta\mbox{\tiny+}n}(\tau)
\!\cdot\!\re^{\mbox{\tiny--}\frac{\pi\ri}{4}\,
m^2\tau}_{\mathstrut}\,.
$$
В настоящей работе мы используем <<$\tau$-представление>>
$\vartheta,\theta$-функций. Переход к часто используемому
<<$q$-представлению>> $\big(q=\re^{\pi\ri\tau}_{}\big)$
осуществляется по формуле $\partial_\tau=\frac{-\ri}{\pi\,
q}\partial_q$.

\subsection{Функции Вейерштрасса}
Будем  использовать общепринятые обозначения Вейерштрасса:
$\sigma(x|\om,\omp)=\sigma(x;g_2^{},g_3^{})$,
$\wp(x|\om,\omp)=\wp(x;g_2^{},g_3^{})$ и т.\,д. Инварианты
$(g_2^{},\,g_3^{})$  являются функциями от периодов $(2\om,2\omp)$
(и наоборот) или модуля $\tau=\frac{\omp}{\om}$ и определяются по
известным рядам Вейерштрасса--Эйзенштейна. Эти ряды совершенно
непригодны для вычислений. Гурвиц, в своей диссертации (1881), нашел
красивый переход к функциональным рядам Ламберта, которые
используются в теориях, имеют приложения и являются наиболее
эффективными в числовых расчетах:
$$
\begin{array}{l}
\ds g_2^{}(\tau) =20\,\pi^4\,\bbig{\Big\{}\frac{1}{240} + \sum_{\sss
k=1}^{\sss \infty}
\frac{k^3\,\re^{2k\pi\ri\,\tau}}{1-\re^{2k\pi\ri\,\tau}}
\bbig{\Big\}},\qquad g_3^{}(\tau) = \mfrac73\,\pi^6\,
\bbig{\Big\{}\frac{1}{504} - \sum_{\sss k=1}^{\sss \infty}
\frac{k^5\,\re^{2k\pi\ri\,\tau}}{1-\re^{2k\pi\ri\,\tau}}\bbig{\Big\}}\,.
\end{array}
$$

Отыскание периодов $(2\om,2\omp)$ по произвольно заданным параметрам
$(a,b)$ эллиптической кривой $w^2=4\,z^3-a\,z-b$ известно как
эллиптическая модулярная задача обращения. Она
решается\footnote{Заметим, что решения этой задачи в 3-м томе {\em
Матем.\;Энцикл.\/} на стр.\,789 и в книге \cite{a} на стр.\,51
изложены не верно, что не связано с опечатками.} через решение
трансцендентного уравнения $J(\tau)=A$ с привлечением $J$-инварианта
Клейна:
$$
J(\tau)=\frac{a^3}{a^3-27\,b^2}\quad \to\quad\om=\sqrt{\frac
ab\frac{g_3^{}(\tau)}{g_2^{}(\tau)}}\quad \to\quad\omp=\tau\,\om\,.
$$
В вырожденных,  лемнискатическом  $(b=0)$ и эквиангармоническом
$(a=0)$ случаях, необходимы отдельные формулы. Они могут быть явно
получены  и мы приведем выражения (новые) для гауссовской
лемнискатической константы $\om_{\mbox{\tiny L}}^{\mathstrut}$ и
эквиангармонической константы $\om_{\mbox{\tiny E}}^{\mathstrut}$:
$$
\left\{
\begin{array}{ll} \ds
\om_{\mbox{\tiny L}}^{\mathstrut}&\!\!\!\!=
\!\!\sqrt[-4]{8\,a\,}\,\pi\cdot\vartheta_{4_{\ds\mathstrut}}^2(2\ri)\\
\ds\omp&\ds\!\!\!\!=\ri\,\om_{\mbox{\tiny L}}^{\mathstrut}
\end{array}\right.,
\qquad \left\{
\begin{array}{ll} \ds
\om_{\mbox{\tiny E}}^{\mathstrut}&
\ds\!\!\!\!=\!\!\sqrt[-12]{\!-27\,b^2\,}\,\pi\cdot\ded^2(\epsilon)\\
\ds\omp&\ds\!\!\!\!=\epsilon\,\om_{\mbox{\tiny E}}^{\mathstrut}\,,
\quad\epsilon=-\mfrac{1^{\ds\mathstrut}}{2}+\mfrac{\sqrt{3}}{2}\,\ri
\end{array}\right.,
$$
где $\ded$ --- функция Дедекинда (\S\,2.3). В силу соотношений
однородности для функций $\sigma,\,\zeta,\,\wp,\,\wpp$, полупериоды
$(\om,\omp)$ или инварианты $(g_2^{},g_3^{})$ могут быть заменены
одной величиной ---  модулем  $\tau$. Такие функции будем обозначать
следующим образом:
$$
\sigma(x|\tau)\equiv\sigma(x|1,\tau),\quad
\zeta(x|\tau)\equiv\zeta(x|1,\tau)\quad
\wp(x|\tau)\equiv\wp(x|1,\tau)\quad
\wpp(x|\tau)\equiv\wpp(x|1,\tau)\,.
$$
$\eta$-функция Вейерштрасса определяется формулой
$\eta(\tau)=\zeta(1|1,\tau)$, а для ее вычисления используется
следующий ряд и его модулярное свойство:
$$
\eta(\tau)=2\pi^2\,\bbig{\Big\{} \frac{1}{24}- \sum_{\sss k=1}^{\sss
\infty} \frac{\re^{2k\pi\ri\,\tau}}
{(1-\re^{2k\pi\ri\,\tau})^2}\bbig{\Big\}}, \qquad
\ds\eta\!\left(\mfrac{a\,\tau+b}{c\,\tau+d}\right) =\ds
(c\tau+d)^2\,\eta(\tau)-\frac{\pi\ri}{2}\,c\,(c\,\tau+d)\,,
$$
где числа $(a,b,c,d)$ являются целыми и $a\,d-b\,c=1$. Модулярные
преобразования представляют не только теоретический интерес. От
значения модуля сильно зависит сходимость функциональных рядов.
Перенося $\tau$  в фундаментальную область модулярной группы
$\boldsymbol{\Gamma}(1)=\mathrm{PSL}_2(\mathbb{Z})$ (процесс легко
автоматизируется), получаем значения $\tau=\epsilon$, у которых
мнимая часть имеет минимально допустимое значение:
$\boldsymbol\Im(\tau)=\frac{\sqrt{3}}{2}$. В такой <<самой худшей>>
точке ряды сходятся чрезвычайно быстро.

Три функции Вейерштрасса $\sigma_1^{},\,\sigma_2^{},\,\sigma_3^{}$
определяются следующими выражениями:
$$
\begin{array}{l}
\ds \sigmalambda\!(x|\om,\omp)= \frac{\theta_{\sss\lambda+1}\!
\big(\frac{x}{2\om}\big|\frac{\om{\sss'}}{\om}\big)}
{\vartheta_{\sss\lambda+1}\! \big(\frac{\om{\sss'}}{\om}\big)}\,
\re^{\eta(\om,\om{\sss'})\,\frac{x^2}{2\om}}_{\mathstrut}\,.
\end{array}
$$
$\sigma$-функция Вейерштрасса, как функция от
$(x,\,g_2^{},\,g_3^{})$, удовлетворяет  линейным дифференциальным
уравнениям, полученных Вейерштрассом:
$$
\left\{
\begin{array}{r}
\ds
x\,\frac{\partial\sigma}{\partial x}-
4\,g_2^{}\,\frac{\partial\sigma}{\partial g_2^{}}-
6\,g_3^{}\,\frac{\partial\sigma}{\partial g_3^{}}-\sigma=0 \\\\
\ds \frac{\partial^2\sigma}{\partial x^2}-
12\,g_3^{}\,\frac{\partial\sigma}{\partial g_2^{}}-
\frac23\,g_2^2\,\frac{\partial\sigma}{\partial g_3^{}}
+\frac{1}{12}\,g_2^{}\,x^2\,\sigma=0
\end{array}\right..
$$
Эти уравнения позволяют получить рекуррентные
соотношения на коэффициенты степенных разложений функции $\sigma$:
\begin{equation}\label{sigma}
\sigma(x;g_2^{},g_3^{})= \displaystyle
C_0\,x+C_1\,\frac{x^3}{3!}+\cdots=
x-\frac{g_2^{}}{240}\,x^5-\frac{g_3^{}}{840}\,x^7 +\cdots
\end{equation}
(нормировка $\sigma(0)=0$, $\sigma'(0)=1$, $\sigma''(0)=0$ является
стандартной). Известны две такие классические рекурренции. Одна
принадлежит Альфену \cite[{\bf I}:\,стр.\,300]{halphen}:
\begin{equation}\label{D}
C_k=\mbox{\large$\widehat{\boldsymbol{\mathfrak D}}$}\,C_{k-1}
-\frac16\,(k-1)(2\,k-1)\,g_2^{}\,C_{k-2}\,, \quad
\mbox{где\;\;\large$\widehat{{\mathfrak D}}$}=
12\,g_3^{}\,\frac{\partial}{\partial g_2^{}}+
\frac23\,g_2^2\,\frac{\partial}{\partial g_3^{}}
\end{equation}
(в других обозначениях ее выписывал также Вейерштрасс \cite[{\bf
V}:\,стр.\,49]{we} и даже Якоби (см.\,\S\,7)). Вторая получена
Вейерштрассом:
$$
\sigma(x;g_2^{},g_3^{})=\mbox{\large$\ds\sum_{\sss
m,n=0}^{\sss\infty}$} \,A_{m,\,n}
\left(\!{\mfrac{g_2^{}}{2}}\!\right)^{\!m} \big( 2g_3^{}\big)^n
{\mfrac{x^{4m+6n+1}} {(4m+6n+1)!}}\,,
$$
$$
\begin{array}{l}
\ds A_{m,\,n}= \mfrac{16}{3}\,(n+1)\,A_{m-2,\,n+1}
+3\,(m+1)\,A_{m+1,\,n-1} -{}\\
\ds\phantom{A_{m,\,n}=}-\mfrac13\,(2m+3n-1)(4m+6\,n-1)\,
A_{m-1,\,n}{}^{\ds{}^{\ds\mathstrut}}\,,
\end{array}
$$
где $A_{\sss 0,0}=1$ и $A_{m,n}=0$ при $(n,m)<0$. Недавно, в связи с
обобщениями на $\boldsymbol{\sigma}$-функции Клейна, была получена
еще одна рекурренция \cite{eilbeck} (существуют и другие
рекурренции). Из всех рекурренций, рекурренция Вейерштрасса наименее
затратная, причем сравнительная степень ее эффективности очень
быстро растет с ростом порядка разложений. Это объясняется тем, что
все коэффициенты уже сгруппированы по параметрам, а сама рекурренция
содержит только умножение  чисел. Вейерштрасс отдельно доказывает,
что числа $A_{m,n}$ являются целыми \cite[{\bf V}:\,стр.\,50]{we}.
Этими же соображениями мы будем руководствоваться при построении
$\sigma$- и $\theta$-рядов в \S\S\,3--4.

\subsection{Функция Дедекинда}
По причине совпадения стандартных обозначений для функции
Вейерштрасса $\eta(\tau)$ и функции Дедекинда, будем использовать
для последней знак $\ded(\tau)$:
$$
\ded(\tau)= \re^{\frac{\pi\ri}{12}\,\tau}_{\mathstrut}\,
{\prod\limits_{\sss k=1}^{\sss \infty}}^
{\ds\mathstrut}_{\ds\mathstrut} \big(1-\re^{2k\pi\ri\,\tau}\big)=
\re^{\frac{\pi\ri}{12}\,\tau}_{\mathstrut} {\sum\limits_{\sss
k=-\infty}^{\sss \infty}}
\!(-1)^k\,\re^{(3k^2+k)\pi\ri\,\tau}\,.\qquad \big(\mbox{Эйлер
(1748)}\big)
$$
Функция Дедекинда связана с функциями Якоби--Вейерштрасса
дифференциально и алгебраически посредством уравнений
$$
\frac{1}{\ded}\,\frac{d\ded}{d\tau}=\frac{\ri}{\pi}\,\eta\,,\qquad
2\,\ded^3=\vartheta_2\,\vartheta_3\,\vartheta_4\,.
$$

\section{Разложения в ряды $\sigma$-функций Вейерштрасса}

\noindent Рекурренция Вейерштрасса $A_{m,n}$ имеет следующую
графическую иллюстрацию: \medskip

\centerline{\unitlength=0.8mm
\begin{picture}(60,60)
\put(0,15){\line(1,0){57}} \put(0,15){\line(0,1){43}}
\multiput(29.7,15)(0,1){24}{\scriptsize .}
\multiput(0,39.75)(1,0){21}{\scriptsize .}
\multiput(23.3,39.75)(1,0){6}{\scriptsize
.}\put(22,40){\circle{0.1}} \put(30,40){\circle{2}}
\put(30,40){\circle*{0.9}} \put(29,42.5){{\scriptsize$(\!m,n\!)$}}
\put(22,40){\circle*{2}}
\put(14.3,36){{\scriptsize$(\!m\!-\!1,n\!)$}}
\multiput(21.5,15)(0,1){20}{\scriptsize .} \put(14,48){\circle*{2}}
\put(10.5,51){{\scriptsize$(\!m\!-\!2,n\!+\!1\!)$}}
\multiput(0,48)(1,0){13}{\scriptsize .}
\multiput(13.5,15)(0,1){33}{\scriptsize .} \put(14,48){\circle{0.1}}
\put(38,32){\circle{0.1}}\put(38,32){\circle*{2}}
\put(36,35){{\scriptsize$(\!m\!+\!1,n\!-\!1\!)$}}
\multiput(37.7,15)(0,1){18}{\scriptsize .}
\multiput(0,31.7)(1,0){37}{\scriptsize .} \put(50,12){$m$}
\put(-4,52){$n$} \put(-2,12){\footnotesize $0$}
\put(28,42){\line(1,-1){28}} \put(-1,55.5){\line(2,-1){34}}
\end{picture}
}\vspace{-0.5cm}

\noindent
Она означает, что для вычисления точки $(m,n)$ необходимы
все точки лежащие внутри указанного четырехугольника.

Обозначив $[n]$ как целую часть числа $n$, не трудно показать, что
формула  классического $\sigma$-ряда \reff{sigma}, явным образом
сгруппированная по $x$, имеет следующий вид:
\begin{equation}\label{wei}
\sigma(x;g_2^{},g_3^{})= \mbox{\large$\ds\sum_{\scriptstyle
k=0}^\infty$}\, \bbig{\Big\{} \sum_{\sss\nu=[k/3]}^{\sss k/2}
\!2^{2k-5\nu} A_{3\nu-k,\,k-2\nu}\,\cdot g_2^{3\nu-k}\, g_3^{k-2\nu}
\bbig{\Big\}} \frac{x^{2k+1}}{(2k+1)!}\ .
\end{equation}

Обозначим $\elambda\equiv \wp(\om_{\!\sss\lambda}^{}|\om,\omp)$.
Тогда функции $\sigmalambda$ удовлетворяют  уравнениям Альфена:
$$
\left\{
\begin{array}{rcl}
\ds
x\,\frac{\partial\sigmalambda}{\partial x}-
2\,\elambda\,\frac{\partial\sigmalambda}{\partial\elambda}-
4\,g_2^{}\,\frac{\partial\sigmalambda}{\partial g_2^{}}\!\!\!&=&\!\!\!0\\ \\
\ds \frac{\partial^2\sigmalambda}{\partial x^2}-
\left(4\,\elambda^{\!\!2}-\mfrac23\,g_2^{}\right)\!
\frac{\partial\sigmalambda}{\partial \elambda}-
12\,\big(4\,\elambda^{\!\!3}-g_2^{}\,\elambda\big)
\frac{\partial\sigmalambda}{\partial g_2^{}}+
\left(\elambda+\mfrac{1}{12}\,g_2^{}\,x^2\right)\!
\sigmalambda\!\!\!&=&\!\!\!0
\end{array}
\right..
$$
Отсюда выводим аналоги  рекурренции \reff{wei} для других степенных
$\sigma$-рядов:
$$
\sigmalambda\!(x;\elambda,g_2^{})=
\mbox{\large$\ds\sum_{\scriptstyle k=0}^\infty$}\, \bbig{\Big\{}
\sum_{\sss \nu=0}^{\sss
k/2}\,2^{-\nu}\,\mathfrak{B}_{k-2\nu,\,\nu}\cdot
\elambda^{\!\!k-2\nu}g_2^\nu \bbig{\Big\}}\,\frac{x^{2k}}{(2k)!}\, ,
$$
\smallskip
$$
\begin{array}{l}
\ds \mathfrak{B}_{m,n}=24\,(n+1)\,\mathfrak{B}_{m-3,\,n+1}+
(4\,m-12\,n-5)\,\mathfrak{B}_{m-1,\,n}-{}\\ \\\ds
\phantom{\mathfrak{B}_{m,n}=}
{}-\mfrac43\,(m+1)\,\mathfrak{B}_{m+1,\,n-1}-
\mfrac13(m+2\,n-1)(2\,m+4\,n-3)\,\mathfrak{B}_{m,\,n-1}\, .
\end{array}
$$
Здесь $\mathfrak{B}_{\sss 0,0}=1$ и $\mathfrak{B}_{m,\,n}=0$ при
$(m,n)<0$, а все коэффициенты $\mathfrak{B}_{m,\,n}$ ---
целочисленны. Вейерштрасс выписывал рекурренции для функций
$S_\lambda=\re^{\frac12\elambda x^2}_{\mathstrut}\sigmalambda$ в
представлении $\big(\elambda,\,\varepsilon_{\sss\lambda}^{}=
3\elambda^{\!\!2}-\frac14g_2^{}\big)$ \cite[{\bf
II}:\,стр.\,253--4]{we}. Выбор таких функций и такого представления
по всей видимости диктовался соображением  получить рекурренцию
простейшего вида. Как и в случае $\sigma$-функции, она является
4-членной. Можно построить универсальный ряд для всех функций
$\sigma$. Они удовлетворяют дифференциальным уравнениям
\begin{equation}\label{epsilon}
\left\{
\begin{array}{rcl}
\ds x\,\frac{\partial\Xi}{\partial x}-
2\,\elambda\,\frac{\partial\Xi}{\partial\elambda}-
4\,g_2^{}\,\frac{\partial\Xi}{\partial g_2^{}}
-(1-\varepsilon)\,\Xi\!\!\!&=&\!\!\!0\\ \\
\ds \frac{\partial^2\Xi}{\partial x^2}-
\left(4\,\elambda^{\!\!2}-\mfrac23\,g_2^{}\right)\!
\frac{\partial\Xi}{\partial \elambda}-
12\,\big(4\,\elambda^{\!\!3}-g_2^{}\,\elambda\big)
\frac{\partial\Xi}{\partial g_2^{}}+
\left(\varepsilon\,\elambda+\mfrac{1}{12}\,g_2^{}\,x^2\right)
\Xi\!\!\!&=&\!\!\!0
\end{array}
\right.,
\end{equation}
где случай $\Xi=\sigmalambda$ соответствует $\varepsilon=1$, а
$\Xi=\sigma$ соответствует $\varepsilon=0$  и произвольному
$\elambda$. Величина $\varepsilon$ совпадает с четностью
$\theta$-характеристики. Универсальный ряд и целочисленная
рекурренция тогда примут следующий вид:
\begin{equation}\label{universal}
\Xi(x;\elambda,g_2^{})=\mbox{\large$\ds\sum_{k=0}^{\infty}$}\,\bbig{\Big\{}
\sum_{\sss\nu=0}^{\sss k/2}\,2^{-\nu}\,
\mathfrak{B}_{k-2\nu,\,\nu}^{{\sss(}\varepsilon{\sss)}}\cdot
\elambda^{\!\!k-2\nu}g_2^\nu\bbig{\Big\}}\,\frac{x^{2k+1-\varepsilon}}
{(2k+1-\varepsilon)!}\,,
\end{equation}
\smallskip
$$
\begin{array}{l}
\ds
\mathfrak{B}_{m,n}^{{\sss(}\varepsilon{\sss)}}=
24\,(n+1)\,\mathfrak{B}_{m-3,\,n+1}^{{\sss(}\varepsilon{\sss)}}+
(4\,m-12\,n-4-\varepsilon)\,\mathfrak{B}_{m-1,\,n}^
{{\sss(}\varepsilon{\sss)}}-{}\\ \\
\ds\phantom{\mathfrak{B}_{m,n}^{{\sss(}\varepsilon{\sss)}}=}
-\mfrac43\,(m+1)\,\mathfrak{B}_{m+1,\,n-1}^{{\sss(}\varepsilon{\sss)}}-
\mfrac13
(m+2\,n-1)(2\,m+4\,n-1-2\,\varepsilon)\,\mathfrak{B}_{m,\,n-1}
^{{\sss(}\varepsilon{\sss)}}\, .
\end{array}
$$
Представление (\ref{universal}) является 5-членным, в отличие от
4-членной рекурренции Вейерштрасса $A_{m,n}$, но представление
$(g_2^{},g_3^{})$  для функций $\sigmalambda$ не существует. Переход
между парами  $(\elambda,g_2^{})$ и $(\elambda,\emu)$ взаимно
однозначен, поэтому универсальная рекурренция может быть записана в
любом из этих представлений. В последнем случае рекурренция
симметрична.

\section{Степенные ряды для $\theta$-функций Якоби}

\noindent Поскольку $\theta$-функции  являются фундаментальными
объектами во многих вопросах,  мы приведем для их степенных
разложений максимально упрощенные (канонические)
формулы\footnote{Исторический комментарий. Красивые рекурренции на
плоскости (аналоги $A_{m,n}$) Вейерштрасс получал еще в 1840-х
годах, когда пользовался своими функциями {\em Al\/}, но не
якобиевскими $\Theta$ и $H$ \cite[{\bf I}]{we}. Тогда же появился
важный множитель $\re^{Ax^2}$, а потом и функция $\sigma$. Примерно
в это же время Якоби рассматривал степенные разложения для своих
функций $\theta$ (см. далее \S\,7).}. Не трудно видеть, что этими
рядами будут ряды с коэффициентами, полиномиальными по переменным
$\eta(\tau)$, $\vartheta(\tau)$. Это следует из очевидных формул
$$
\theta_1(x|\tau)=\pi\,\ded^3(\tau)
\!\cdot\!\re^{-2\eta(\tau)x^2}_{\mathstrut}\sigma(2x|\tau)\;,
\qquad\quad \theta_{\sss
\lambda}\!(x|\tau)=\vartheta_{\sss\!\lambda}^{}\!(\tau)
\!\cdot\!\re^{-2\eta(\tau)x^2}_{\mathstrut} \sigma_{\sss
\!\lambda-1}^{}\!(2x|\tau)\,.
$$
Здесь, вместо вейерштрассовских представлений в переменных
$g_2^{},g_3^{},\elambda$, логичнее использовать их эквиваленты в
$\vartheta$-константах.  Выражения для точек ветвления $\elambda$
через $\vartheta$-константы хорошо известны. В свою очередь,
$\vartheta$-константы связаны тождеством Якоби
$\vartheta_3^4=\vartheta_2^4+\vartheta_4^4$. Это позволяет делать
переходы между представлениями, выбирая произвольную пару. Будем
говорить {\em $(\alpha,\beta)$-представление\/}, если формулы
записаны через константы
$(\vartheta_{\sss\alpha0},\vartheta_{\sss0\beta})$ при
$(\alpha,\beta)\ne(0,0)$.

\subsection{Оператор Альфена}
Оператор \reff{D} имеет как вейерштрассовские представления
$$
\begin{array}{l}
\ds {\mbox{\large$\widehat{\boldsymbol{\mathfrak D}}$}}=
\left(4\,\elambda^{\!\!2}-\mfrac23 g_2^{}\right)\,
\frac{\partial}{\partial\elambda}+
12\,\big(4\,\elambda^{\!\!3}-g_2^{}\,\elambda\big)\frac{\partial}{\partial
g_2^{}}=\\ \\
\ds \phantom{\mbox{\large$\widehat{\boldsymbol{\mathfrak D}}$}}=
\mfrac43\big(\elambda^{\!\!2}-2\,\emu\elambda-2\,\emu^{\!\!2}\big)
\frac{\partial}{\partial \elambda}+
\mfrac43\big(\emu^{\!\!2}-2\,\elambda\emu-2\,\elambda^{\!\!2}\big)
\frac{\partial}{\partial\emu}\,,
\end{array}
$$
так и $\vartheta$-константное. Например в представлении
$(\vartheta_2,\vartheta_4)$ он выглядит так:
$$
\mbox{\large$\widehat{\boldsymbol{\mathfrak D}}$}= \frac{\pi^2}{3}
\bbig(\vartheta_4^8+2\,\vartheta_2^4\,\vartheta_4^4\bbig)
\frac{\partial}{\partial(\vartheta_4^4)} -\frac{\pi^2}{3}
\bbig(\vartheta_2^8+2\,\vartheta_2^4\,\vartheta_4^4\bbig)
\frac{\partial}{\partial(\vartheta_2^4)}\;.
$$
Обозначим
$\spin{\alpha}=(-1)^\alpha$. Общее $(\alpha,\beta)$-представление
оператора Альфена имеет вид
$$
\mbox{\large$\widehat{\boldsymbol{\mathfrak D}}$}=
\frac{\pi^2}{3}\bbig(\spin{\beta}\,
\vartheta_{\sss\alpha0}{}^{\!\!\!\!\!\!8}\;+ 2\,\spin{\alpha}\,
\vartheta_{\sss\alpha0}{}^{\!\!\!\!\!\!4}\;\;
\vartheta_{\sss0\beta}{}^{\!\!\!\!\!\!4}\; \bbig)
\frac{\partial}{\partial(\vartheta_{\sss\alpha0}{}^{\!\!\!\!\!\!4}\;)}
-\frac{\pi^2}{3}\bbig(\spin{\alpha}\,
\vartheta_{\sss0\beta}{}^{\!\!\!\!\!\!8}\;+ 2\,\spin{\beta}\,
\vartheta_{\sss\alpha0}{}^{\!\!\!\!\!\!4}\;\;
\vartheta_{\sss0\beta}{}^{\!\!\!\!\!\!4}\;\bbig)
\frac{\partial}{\partial(\vartheta_{\sss0\beta}{}^{\!\!\!\!\!\!4}\;)}
\;.
$$
Пусть $\varepsilon$ определяется по четности характеристики функции
$\theta_{{\sss\!}\alpha\beta}(x|\tau)$:
$$
\varepsilon=\frac{\spin{\alpha\beta}+1}{2}\quad
\Rightarrow\quad\left\{\!\!\!
\begin{array}{l}
\varepsilon=0\quad \mbox{при \ }
\theta_{{\sss\!}\alpha\beta}=\pm\,\theta_1\\
\varepsilon=1\quad \mbox{при \ }
\theta_{{\sss\!}\alpha\beta}=\pm\,
\theta_{2,3,4}{}^{\ds\mathstrut}
\end{array}\!\!\!
\right.\;.
$$
Тогда уравнение \reff{epsilon} на функции
$\Xi=(\sigma,\sigmalambda)$ принимает вид
$$
\frac{\partial^2\,\Xi}{\partial x^2}-
\mbox{\large$\widehat{\boldsymbol{\mathfrak D}}$}\, \Xi
+\Big\{\varepsilon\,\elambda(\vartheta)+\mfrac{\pi^4}{12^2}
\big(\vartheta_2^8+\vartheta_2^4\,\vartheta_4^4+ \vartheta_4^8 \big)
x^2\Big\}\,\Xi=0\,.
$$
Это же уравнение  в $(\alpha,\beta)$-представлении имеет вид
\begin{equation}\label{eq}
\frac{\partial^2\,\Xi}{\partial
x^2}-\mbox{\large$\widehat{\boldsymbol{\mathfrak D}}$}\, \Xi
+\Big\{e_{\sss\gamma\delta}(\vartheta) + \mfrac{\pi^4}{12^2}\,
\bbig[ \vartheta_{\sss\alpha0}^8 +\spin{\alpha+\beta}\,
\vartheta_{\sss\alpha0}^4\, \vartheta_{\sss0\beta}^4+
\vartheta_{\sss0\beta}^8\, \bbig]\, x^2 \Big\}\,\Xi=0\,,
\end{equation}
где, не зависимо от $(\alpha,\beta)$-представления,
$e_{\sss\gamma\delta}(\vartheta)$ соответствуют   функции $\sigma$
или $\sigmalambda$:
\begin{equation}\label{ee}
e_{\sss\gamma\delta}(\vartheta)= \mfrac{\pi^2}{12} \big(
\spin{\gamma}\,\vartheta_{\sss0\delta}^4-
\spin{\delta}\,\vartheta_{\sss\gamma0}^4\big) \quad\Rightarrow \quad
\left\{
\begin{array}{ll}
\ds e_{\sss 00}\equiv 0,&
e_{\sss 01}\equiv e_{\sss 1}\\
e_{\sss 11}\equiv e_{\sss 2},& e_{\sss 10}\equiv e_{\sss 3}
\end{array}
\right\}\,.
\end{equation}
Представление называется собственным (симметричным) если
$(\gamma,\delta)=(\alpha,\beta)$.

Степенные разложения для $\theta$-функций, которые последуют далее,
получаются, используя разные представления оператора
{\large$\widehat{\boldsymbol{\mathfrak D}}$} в  уравнении \reff{eq}
на функцию $\Xi$. Мы опускаем детали выводов и доказательств.

\subsection{Функция $\theta_1$} Разложение функции
\begin{equation}\label{theta1}
\begin{array}{l}
\ds\theta_1(x|\tau)=\sum_{k=0}^{\infty}\,C_k(\tau)\!\cdot\!
x^{2k+1}=\\
\ds \phantom{\theta_1(x|\tau)}=2\pi\,\ded^3 \,\Big\{
x-2\,\eta\!\cdot\!x^3+\bbig(2\,\eta^2- \mfrac{\pi^4}{180}
\big(\vartheta_2^8+\vartheta_2^4\,\vartheta_4^4+\vartheta_4^8\big)
\bbig) \!\cdot\!x^5+\cdots\Big\}^{{}^{\ds\mathstrut}}
\end{array}
\end{equation}
определяется  аналитическими выражениями
\begin{equation}\label{t1}
\begin{array}{l}
\ds \theta_1(x|\tau)=\phantom{\ded^3}2\pi\,\,
\mbox{\large$\ds\sum_{\scriptstyle k=0}^\infty$}
\,\frac{(4\pi\ri)^k}{(2k+1)!}\,\frac{d^k\ded^3}{d\tau^k}
\!\cdot\!x^{2k+1}=\\\\  \ds \phantom{\theta_1(x|\tau)}=
2\pi\,\ded^3\,\, \mbox{\large$\ds\sum_{\scriptstyle k=0}^\infty$}
\,(-2)^k\,\bbig{\Big\{} \sum_{\sss\nu=0}^{\sss k}
\Big(\!\!-\mfrac{\pi^2}{6}\Big)^\nu
\mfrac{\eta^{k-\nu}\,\boldsymbol{\mathcal{N}}_{\!\nu}(\vartheta)}{(k-\nu)!\,(2\nu+1)!}
\bbig{\Big\}} \,x^{2k+1}\;,
\end{array}
\end{equation}
где полиномы $\boldsymbol{\mathcal{N}}_{\!\nu}(\vartheta)$, в
зависимости от комбинаций $\vartheta$-констант, имеют вид
$$
\boldsymbol{\mathcal{N}}_{\!\nu}(\vartheta)=
\mbox{\large$\ds\sum_{s=0}^{\nu}$}\;
\mbox{\footnotesize$\textstyle\left\{\!\!\!
\begin{array}{r}
\mbox{\small$\mathfrak{G}_{\nu\mbox{\tiny--}s,\,s}$}\cdot\vartheta_4^{4s\mathstrut}\,
\vartheta_2^{4{\sss(}\nu-s{\sss)}}\\
(-1)^s\,\mbox{\small$\mathfrak{G}_{\nu\mbox{\tiny--}s,\,s}$}\cdot
\vartheta_3^{4s\ds\mathstrut}\,
\vartheta_4^{4{\sss(}\nu-s{\sss)}\ds\mathstrut}\\
(-1)^s\,\mbox{\small$\mathfrak{G}_{s,\,\nu\mbox{\tiny--}s}$}\cdot\vartheta_3^{4s\ds\mathstrut}\,
\vartheta_2^{4{\sss(}\nu-s{\sss)}\ds\mathstrut}
\end{array}\!\!\!
\right\}$}\,.
$$

\noindent
Целочисленная рекурренция $\mathfrak{G}_{m,n}$
$(\mathfrak{G}_{\sss 0,\,0}=1)$ здесь выглядит следующим образом:
$$
\begin{array}{l}
\ds \mathfrak{G}_{m,\,n}=4\,(n-2\,m-1)\,\mathfrak{G}_{m,\,n-1}-
4\,(m-2\,n-1)\,\mathfrak{G}_{m-1,\,n}-{}
\\ \\
\ds\phantom{\mathfrak{G}_{m,\,n}=}
-2\,(m+n-1)(2\,m+2\,n-1)\big(\mathfrak{G}_{m-2,\,n}+
\mathfrak{G}_{m-1,\,n-1}+\mathfrak{G}_{m,\,n-2}\big)
\end{array}.
$$

\noindent
Рекурренция $\mathfrak{G}_{m,n}$ антисимметрична:
$\mathfrak{G}_{m,n}=(-1)^{m+n}\,\mathfrak{G}_{n,m}$. Можно выписать
представление  типа $\theta_1=\sum
C_{mnp}\,g_2^m\,g_3^n\,\eta^p\,x^k$ через рекурренцию Вейерштрасса,
но  $\mathfrak{G}_{m,n}$ эффективнее рекурренции $A_{m,n}$, так как
полиномы уже сгруппированы по $\vartheta$-конс\-тан\-там. Отметим
тот факт, что  производные $\theta_1^{\sss(2k+1)}(0|\tau)$, т.\,е.
выражения перед $x^{2k+1}$ в (\ref{theta1}--\ref{t1}), генерируют
полиномы по $(\eta,\vartheta)$, которые точно интегрируются $k$ раз
по $\tau$.

\subsection{Функции $\theta_{2,3,4}$}
Разложение функций $\thetaAB{\alpha}{\beta}=\pm\theta_{2,3,4}$ вида
\begin{equation}\label{theta234}
\thetaAB{\alpha}{\beta}(x|\tau) =
\sum_{k=0}^{\infty}\,C_k^{\sss(\alpha,\beta)}(\tau)\!\cdot\!x^{2k}=
\varthetaAB{\alpha}{\beta}- \varthetaAB{\alpha}{\beta}\, \Big\{
2\,\eta+\mfrac{\pi^2}{6}
\bbig(\spin{\beta}\,\varthetaAB{\alpha\mbox{\tiny--}1}{0}^4-
\spin{\alpha}\,\varthetaAB{0}{\beta\mbox{\tiny--}1}^4\bbig)\Big\}\,x^2+\cdots
\end{equation}
определяется  аналитическим выражением
\begin{equation}\label{series}
\thetaAB{\alpha}{\beta}(x|\tau) = \mbox{\large$\ds\sum_{\scriptstyle
k=0}^\infty$} \,\frac{(4\pi\ri)^k}{(2k)!}\,
\frac{d^k\varthetaAB{\alpha}{\beta}}{d\tau^k} \!\cdot\!x^{2k}\;.
\end{equation}
В собственном представлении ряд \reff{series} имеет сгруппированный
вид
\begin{equation}\label{t234}
\thetaAB{\alpha}{\beta}(x|\tau)= \varthetaAB{\alpha}{\beta}\,\,
\mbox{\large$\ds\sum_{\scriptstyle k=0}^\infty$}
\,(-2)^k\,\bbig{\Big\{} \sum_{\sss\nu=0}^{\sss k}
\Big(\!\!-\mfrac{\pi^2}{6}\Big)^\nu
\mfrac{\eta^{k-\nu}\,\boldsymbol{\mathcal{N}}_\nu^{\sss(\alpha,\beta)}
(\vartheta)}{(k-\nu)!\,(2\nu)!} \bbig{\Big\}}\, x^{2k}
\end{equation}
с универсальной целочисленной рекурренцией $(\mathfrak{G}_{\sss
0,\,0}^{\sss(\alpha,\beta)}=1)$ следующего вида:
$$
\boldsymbol{\mathcal{N}}_\nu^{\sss(\alpha,\beta)}(\vartheta)=
\mbox{\large$\ds\sum_{s=0}^{\nu}$}\;
\mathfrak{G}_{s,\,\nu\mbox{\tiny--}s}^{\sss(\alpha,\beta)} \cdot
\varthetaAB{\alpha\mbox{\tiny--}1}{0}^{4s}\,
\varthetaAB{0}{\beta\mbox{\tiny--}1}^{4{\sss(}\nu-s{\sss)}\mathstrut}
\;,
$$
$$
\begin{array}{rl}
\ds \mathfrak{G}_{m,\,n}^{\sss(\alpha,\beta)}\!\!&=\,
\spin{\alpha}\,(4\,n-8\,m-3)\,
\mathfrak{G}_{m,\,n-1}^{\sss(\alpha,\beta)}-
\spin{\beta}\,(4\,m-8\,n-3)\,
\mathfrak{G}_{m-1,\,n}^{\sss(\alpha,\beta)}-
\\ \\&\ds\phantom{=\,}
-2\,(m+n-1)(2\,m+2\,n-3)
\big(\mathfrak{G}_{m-2,\,n}^{\sss(\alpha,\beta)}+
\spin{\alpha+\beta}\, \mathfrak{G}_{m-1,\,n-1}^{\sss(\alpha,\beta)}+
\mathfrak{G}_{m,\,n-2}^{\sss(\alpha,\beta)}\big)\;.
\end{array}
$$

Хотя эти целочисленные рекурренции являются вполне эффективными,
имеются  дополнительные свойства, сокращающие вычисления наполовину.
Симметрия рекурренций $\mathfrak{G}_{m,\,n}^{\sss(\alpha,\beta)}$ по
перестановке индексов определяется формулами
$$
\mathfrak{G}_{n,\,m}^{\sss(\alpha,\beta)}=
(-1)^{\sss(m+n)(\alpha+\beta+1)}\,
\mathfrak{G}_{m,\,n}^{\sss(\alpha,\beta)}\;,\qquad
\mathfrak{G}_{m,\,n}^{\sss(\beta,\alpha)}=
(-1)^{\sss(m+n)(\alpha+\beta)\phantom{+1}}\,
\mathfrak{G}_{m,\,n}^{\sss(\alpha,\beta)}\;.
$$
Это означает, что фактически имеются только две рекурренции: при
$\alpha=\{1,0\}$ и $\beta=0$. Обозначая
$\mathfrak{G}_{m,\,n}^{\sss(\alpha,0)}=
\mathfrak{G}_{m,\,n}^{\sss(\alpha)}$,
мы будем иметь:
$$
\begin{array}{rl}
\ds \mathfrak{G}_{m,\,n}^{\sss(\alpha)}\!\!&=\,
\spin{\alpha}\,(4\,n-8\,m-3)\,
\mathfrak{G}_{m,\,n-1}^{\sss(\alpha)}- (4\,m-8\,n-3)\,
\mathfrak{G}_{m-1,\,n}^{\sss(\alpha)}-{}
\\ \\&\phantom{=\,}
-2\,(m+n-1)(2\,m+2\,n-3) \big(\mathfrak{G}_{m-2,\,n}^{\sss(\alpha)}+
\spin{\alpha}\, \mathfrak{G}_{m-1,\,n-1}^{\sss(\alpha)}+
\mathfrak{G}_{m,\,n-2}^{\sss(\alpha)}\big)
\end{array}.
$$
Вышеуказанная перестановка индексов сводится тогда к формулам
$$
\mathfrak{G}_{n,\,m}^{\sss(0)}= (-1)^{\sss(m+n)}\,
\mathfrak{G}_{m,\,n}^{\sss(0)}\;,\qquad
\mathfrak{G}_{n,\,m}^{\sss(1)}= \mathfrak{G}_{m,\,n}^{\sss(1)}\;.
$$
Разложения \reff{t1} и \reff{t234} отличаются лишь множителем и
видом рекурренций. Их можно объединить в одну (вводя четность
$\varepsilon$), но  величина $\spin{\alpha}$ остается, хотя
$(m,n)$-элементы матриц $\mathfrak{G}^{\sss(\beta,\alpha)}$
отличаются лишь знаком.

Четные производные $\theta_{\sss\alpha\beta}^{\sss(2k)}(0|\tau)$,
т.\,е. выражения перед $x^{2k}$ в \reff{theta234}, \reff{t234} дают
бесконечное семейство полиномов по $(\eta,\vartheta)$, которые точно
интегрируются $k$ раз по $\tau$. Их интегрируемость есть следствие
динамической системы, рассматриваемой в \S\,7.

Используя эти разложения можно строить  разложения в точках
$x=\big\{\pm\frac12,\,\pm\frac\tau2\big\}$. Ряды, с точностью до
очевидных модификаций, будут переходить друг в друга.

\section{Дифференцирования функций $(\sigma,\zeta,\wp,\wpp)$}

\noindent Дифференцирования функций Вейерштрасса по инвариантам
известны \cite[{\bf I}]{halphen}, \cite[\bf IV]{tannery}. Правила
перехода между дифференцированиями по $(g_2^{},g_3^{})$ и
$(\om,\omp)$ рассматривались Вейерштрассом и одновременно с ним ---
Фробениусом и Шти\-кель\-бер\-гером (1882) \cite{halphen}. Обозначим
знак отношения периодов как
$\boldsymbol{\mathfrak{s}}=\mathfrak{sign}\bbig(\boldsymbol\Im\big(\frac{\omp}{\om}\big)\bbig)$.
Тогда, прилагая эти правила и упрощая, получаем следующие формулы:
$$
\left\{
\begin{array}{rcl}
\ds\boldsymbol{\mathfrak{s}}
\,\frac{\partial\sigma}{\partial\om}\!\!\!&=&\!\!\! \ds
-\frac{\ri}{\pi}\, \Big\{
\omp\bbig(\wp-\zeta^2-\mfrac{1}{12}\,g_2^{}\,x^2 \bbig)
+2\,\etap\,\big(x\,\zeta-1\big)\Big\}\,\sigma\\ \\
\ds\boldsymbol{\mathfrak{s}}\,
\frac{\partial\sigma}{\partial\omp}\!\!\!&=&\!\!\! \ds
\phantom{-}\frac{\ri}{\pi}\, \Big\{
\om\bbig(\wp-\zeta^2-\mfrac{1}{12}\,g_2^{}\,x^2 \bbig)
\;+2\,\eta\,\big(x\,\zeta-1\big)\Big\}\,\sigma
\end{array}
\right.,
$$
$$
\;\;\,\left\{
\begin{array}{rcl}
\ds\boldsymbol{\mathfrak{s}}\,\frac{\partial\zeta}{\partial\om}\!\!\!&=&\!\!\!
\ds -\frac{\ri}{\pi}\, \Big\{
2\,(\omp\zeta-x\,\etap)\,\wp+\omp\,\bbig(\wpp-\mfrac
16\,x\,g_2^{}\bbig) +2\,\etap\zeta\Big\} \\\\
\ds\boldsymbol{\mathfrak{s}}\,
\frac{\partial\zeta}{\partial\omp}\!\!\!&=&\!\!\! \ds\phantom{-}
\frac{\ri}{\pi}\, \Big\{
2\,(\om\,\zeta-x\,\eta)\,\wp\;+\om\,\bbig(\wpp-\mfrac
16\,x\,g_2^{}\bbig) \;+2\,\eta\,\zeta\Big\}
\end{array}
\right.,
$$
$$
\;\;\;\,\left\{
\begin{array}{rcl}
\ds\boldsymbol{\mathfrak{s}}\,
\frac{\partial\wp}{\partial\om}\!\!\!&=&\!\!\! \ds\phantom{-}
\frac{\ri}{\pi}\, \Big\{ 2\,(\omp\zeta-x\,\etap)\,\wpp
+4\,(\omp\wp-\etap)\,\wp-\mfrac23\,\omp g_2^{}\Big\} \\\\
\ds\boldsymbol{\mathfrak{s}}\,
\frac{\partial\wp}{\partial\omp}\!\!\!&=&\!\!\! \ds
-\frac{\ri}{\pi}\, \Big\{ 2\,(\om\,\zeta-x\,\eta)\,\wpp
\;+4\,(\om\,\wp-\eta)\,\wp{\sss\;\:}-\mfrac23\,\om\, g_2^{}\Big\}
\end{array}
\right.,
$$
$$
\;\;\;\left\{
\begin{array}{rcl}
\ds\boldsymbol{\mathfrak{s}}\,
\frac{\partial\wpp}{\partial\om}\!\!\!&=&\!\!\! \ds\phantom{-}
\frac{\ri}{\pi}\, \bbig\{ 6\,(\omp\wp-\etap)\,\wpp
+(\omp\zeta-x\,\etap)(12\,\wp^2-g_2^{}) \bbig\} \\\\
\ds\boldsymbol{\mathfrak{s}}\,
\frac{\partial\wpp}{\partial\omp}\!\!\!&=&\!\!\! \ds
-\frac{\ri}{\pi}\, \bbig\{ 6\,(\om\,\wp-\eta)\,\wpp
\:\,+(\om\,\zeta-x\,\eta)(12\,\wp^2-g_2^{}) \bbig\}
\end{array}
\right.\,.
$$
Полагая в этих уравнениях $(\om=1, \,\omp=\tau)$ и
$\boldsymbol{\mathfrak{s}}=1$, мы приходим к динамической системе
уравнений с параметром $x$:
$$
\left\{
\begin{array}{l}
\ds\;\frac{\partial\sigma}{\partial\tau}= \phantom{-}\frac{\ri}{\pi}
\Big\{ \wp-\zeta^2
+2\,\eta\,(x\,\zeta-1)-\mfrac{1}{12}\,g_2^{}\,x^2\Big\}\,\sigma\\\\
\left. \!\!\!
\begin{array}{l}
\ds\;\frac{\partial\zeta}{\partial\tau}= \phantom{-} \frac{\ri}{\pi}
\Big\{ \wpp+2\,(\zeta-x\,\eta)\,\wp+ 2\,\eta\,\zeta
-\mfrac16\,g_2^{}\,x\Big\}\\\\
\ds\;\frac{\partial\wp}{\partial\tau}=-\frac{\ri}{\pi} \Big\{
2\,(\zeta-x\,\eta)\,\wpp+4\,(\wp-\eta)\,\wp
-\mfrac23\,g_2^{}\Big\}\\\\
\ds\frac{\partial\wpp}{\partial\tau}= -\frac{\ri}{\pi}\, \bbig\{
6\,(\wp-\eta)\,\wpp+(\zeta-x\,\eta)(12\,\wp^2-g_2^{}) \bbig\}
\end{array}\right\}
\end{array}
\right..
$$
Отсюда следует, что дифференциальная замкнутость трех
вейерштрассовских функций $\zeta,\,\wp,\,\wpp$, известная по
переменной $x$, имеется также и по $\tau$, что дополнительно
обозначено правой скобкой. В свою очередь, дифференциальная
замкнутость коэффициентов $g_2^{},g_3^{},\eta$ тоже известна (см.
далее уравнения \reff{g2g3}).

Дифференциальные уравнения по переменным $x$, $\om$, $\omp$, $\tau$
не содержат величину $g_3^{}$. Она, в зависимости от представления
$x$ или $\tau$, является либо алгебраическим интегралом этих
уравнений $g_3^{}(\wp,\wpp)=4\,\wp^3-g_2^{}\,\wp-\wpp^2$, либо, по
этой же формуле, определяет фиксированную алгебраическую связь между
переменными $\wp$ и $\wpp$.

Функция $\sigma(x|\tau)$ удовлетворяет дифференциальному уравнению в
частных производных, которое является аналогом уравнения
теплопроводности для функций $\theta(x|\tau)$:
$$
\pi\,\ri\, \frac{\partial\sigma}{\partial \tau}=
\frac{\partial^2\sigma}{\partial x^2}-2\,x\,\eta\,
\frac{\partial\sigma}{\partial x}+
\left(2\,\eta+\mfrac{1}{12}\,g_2^{}\,x^2 \right)\sigma\,.
$$
Каждая из функций $\zeta,\,\wp,\,\wpp(x|\tau)$ удовлетворяет
обыкновенному дифференциальному уравнению второго порядка по
переменной $\tau$ с переменными коэффициентами $g_{2,3}^{},\, \eta$,
а функция $\sigma$
--- уравнению порядка 3. Например функция $Z=\zeta-x\,\eta$ удовлетворяет уравнению
2-го порядка, получаемого исключением переменной $\wp$  из двух
полиномов
$$
\begin{array}{c}
\ds \big(\pi\ri\, Z_\tau+2\,(\wp+\eta)\,Z
\big)^2=4\,\wp^3-g_2^{}\,\wp-g_3^{}\,,\\\\\ds
-\frac{\pi^2}{8}\,\frac{Z_{\tau\tau}}{Z}=
\frac{\pi\ri}{2}\,\big(Z^2+\wp-2\,\eta\big)\,\frac{Z_\tau}{Z}+
(\wp+\eta)\,Z^2-\wp^2+\eta\,\wp-\eta^2+\frac14\,g_2^{}\,.
\end{array}
$$
Можно показать,  что уравнения на $Z,\wp,\wpp$ не содержат $x$.
Полный список уравнений и их $\theta$-эквивалентов  мы не приводим,
так как они не достаточно компактны (см. \S\,9).

\section{Динамические системы на $\theta$-функции}
\noindent В этом и следующем параграфе мы описываем новое и
фундаментальное свойство $\vartheta$, $\theta$ и $\theta'$-функций
Якоби --- они обладают замкнутым дифференциальным исчислением.
\subsection{Дифференциальные уравнения по $x$}
Пять функций
$$
\theta_1\,,\,\theta_2\,,\,\theta_3\,,\,\theta_4
\quad\mbox{и}\quad\theta_1'\equiv \frac{\partial\theta_1}{\partial
x}
$$
удовлетворяют  замкнутым обыкновенным автономным дифференциальным
уравнениям с коэффициентами $\vartheta_{2,3,4}$ и $\eta$:
\begin{equation}\label{X}
\left\{
\begin{array}{ll}
\ds\frac{\partial\theta_2}{\partial x}=
\frac{\theta_1'}{\theta_1}\,\theta_2- \pi\,\vartheta_2^2\!\cdot\!
\frac{\theta_3\theta_4}{\theta_1},& \ds
\frac{\partial\theta_4}{\partial x}=
\frac{\theta_1'}{\theta_1}\,\theta_4^{}- \pi\,\vartheta_4^2\!\cdot\!
\frac{\theta_2\theta_3}{\theta_1}\,,\qquad
\frac{\partial\theta_1}{\partial x}=\theta_1'
\\\\
\ds\frac{\partial\theta_3}{\partial x}=
\frac{\theta_1'}{\theta_1}\,\theta_3- \pi\,\vartheta_3^2\!\cdot\!
\frac{\theta_2\theta_4}{\theta_1},&\ds
\frac{\partial\theta_1'}{\partial x}=
\frac{\theta_1'^2}{\theta_1}-\pi^2\vartheta_3^2\,\vartheta_4^2
\!\cdot\! \frac{\theta_2^2}{\theta_1}-
4\,\Big\{\eta+\mfrac{\pi^2}{12}\big(\vartheta_3^4+\vartheta_4^4\big)
\Big\} \!\cdot\!\theta_1
\end{array}\right.\!\!\!.
\end{equation}
Эти уравнения следуют из  соотношений между вейерштрассовскими
функциями $(\sigma,\zeta,\wp)(x|\tau)$, взятыми в различных
$\frac12$-периодах. Учитывая то, что $\vartheta_1\equiv 0$, общая
запись $x$-дифференцирования функций $\theta_{1,2,3,4}$ имеет вид:
$$
\frac{\partial\theta_k}{\partial x}= \frac{\theta_1'}{\theta_1}
\,\theta_k- \pi\,\vartheta_k^2\!\cdot\!
\frac{\theta_\nu\,\theta_\mu}{\theta_1}\,,\qquad
\mbox{где}\quad\nu=\frac{8\,k-28}{3\,k-10}\,,\;\;
\mu=\frac{10\,k-28}{3\,k-8}\,.
$$
Редко встречающиеся  дифференциальные соотношения между якобиевскими
рядами\footnote{Эти важные соотношения неявно присутствуют у Якоби
\cite[\bf I]{jacobi}, но не попали даже в подробный справочник по
эллиптическим функциям \cite{we2}, составленный Шварцем на основе
лекций Вейерштрасса, и сами лекции \cite{we}. Дифференциальные
уравнения на отношения $\theta$-функций хорошо известны \cite{we2},
\cite[\S\,21.6]{WW}. Это собственно дифференциальные уравнения на
эллиптические функции.}
\begin{equation}\label{diff}
\frac{\theta_\nu'}{\theta_\nu}-\frac{\theta_\mu'}{\theta_\mu}=
\pi\vartheta_k^2\!\cdot\!\frac{\theta_1\theta_k}{\theta_\nu\theta_\mu}
\,\mathfrak{sign}(\nu-\mu) \qquad (k=2,3,4)
\end{equation}
являются следствиями уравнений \reff{X}, если между их
$\theta$-решениями {\em дополнительно наложить\/} (см. \S\,9.1--2)
известные полиномиальные $\theta$-тождества
$$
\mathfrak{sign}(\nu-\mu)\!\cdot\!\vartheta_k^2\,\theta_1^2=
\vartheta_\mu^2\,\theta_\nu^2-\vartheta_\nu^2\,\theta_\mu^2 \qquad
(k=2,3,4)\,.
$$
Тождество Якоби
$\vartheta_1'=\pi\,\vartheta_2\vartheta_3\vartheta_4$, об известных
доказательствах которого в \cite[стр.\,345]{WW} сказано, что ``ни
одно из них не является простым'', есть автоматическое следствие
уравнений \reff{X}, взятых в точке  $x=0$. В равной степени это
относится и к обобщениям  этого тождества, представленных далее
формулами (\ref{jacobi}--\ref{jac}).

В $\zeta$-форме,  дифференциальные соотношения \reff{diff}
Вейерштрасс получал в обратном порядке \cite[\S\,24--25]{we2}:
дифференцируя $\sigma$-тождества и используя дифференциальные
уравнения на отношения $\sigma$-функций. Отметим, что существенным
является не только порядок рассмотрений, но и появление в формулах,
помимо точек ветвления $\elambda$ (т.\,е. $\vartheta$-констант),
величины $\eta$ и пятого уравнения. Иными словами, в
дифференциальном аппарате с неизбежностью появляется функция
$\theta'$ и период мероморфного эллиптического интеграла  $\eta$.

Как отмечено во введении, предыдущие результаты естественно было
ожидать и с более общей точки зрения. А именно, привлекая абелевы
эллиптические интегралы. В самом деле, дифференциальную замкнутость
образуют мероморфные и логарифмические эллиптические интегралы.
Эллиптические функции являются их частными случаями, когда
эллиптический дифференциал является точным. Логарифмический интеграл
выражается через логарифмы отношений $\theta$-функций, а
канонический мероморфный интеграл (т.\,е. вейерштрассовская функция
$\zeta$) содержит $\theta'$.

\subsection{Дифференциальные уравнения по $\tau$}

Функции $\theta_{1,2,3,4}$ и $\theta_1'$ удовлетворяют замкнутым
обыкновенным неавтономным дифференциальным уравнениям по $\tau$:
$$
\frac{\partial \theta_1}{\partial\tau}=
\mfrac{-\ri}{4\,\pi}\,\frac{\theta_1'^2}{\theta_1}+
\phantom{xxxxxxxxxx\;\;\;\,}+ \mfrac{\pi
\ri}{4}\,\vartheta_3^2\,\vartheta_4^2\!\cdot\!
\frac{\theta_2^2}{\theta_1}\;\;\;\;+
\mfrac{\ri}{\pi}\Big\{\eta+\mfrac{\pi^2}{12}
\big(\vartheta_3^4+\vartheta_4^4\big) \Big\}\!\cdot\!\theta_1\,,
\;\,\,
$$
$$
\frac{\partial \theta_2}{\partial\tau}= \mfrac{-\ri}{4\,\pi}\!
\left\{\frac{\theta_1'}{\theta_1}- \pi\,\vartheta_2^2\!\cdot\!
\frac{\theta_3\theta_4} {\theta_1\theta_2}\right\}^{\!2}\theta_2 +
\mfrac{\pi \ri}{4}\,\vartheta_3^2\,\vartheta_4^2\!\cdot\!
\frac{\theta_1^2}{\theta_2}\;\;\;\;+
\mfrac{\ri}{\pi}\Big\{\eta+\mfrac{\pi^2}{12}
\big(\vartheta_3^4+\vartheta_4^4\big) \Big\}\!\cdot\!\theta_2\,,
\;\,\,
$$\
$$
\frac{\partial \theta_3}{\partial\tau}= \mfrac{-\ri}{4\,\pi}\,
\frac{\theta_1'^2}{\theta_1^2}\,\theta_3
+\mfrac{\ri}{2}\vartheta_3^2\!\cdot\!
\theta_1'\,\frac{\theta_2\theta_4}{\theta_1^2} {\sss \;\;\;}-
\mfrac{\pi\ri}{4}\,\vartheta_2^2\,\vartheta_3^2\!\cdot\!
\frac{\theta_4^2}{\theta_1^2}\,\theta_3+
\mfrac{\ri}{\pi}\Big\{\eta+\mfrac{\pi^2}{12}
\big(\vartheta_3^4+\vartheta_4^4\big) \Big\}\!\cdot\!\theta_3\,,
\;\,\,
$$\
$$
\frac{\partial \theta_4}{\partial\tau}=
\mfrac{-\ri}{4\,\pi}\,\frac{\theta_1'^2}{\theta_1^2}\,\theta_4
+\mfrac{\ri}{2}\vartheta_4^2\!\cdot\!
\theta_1'\,\frac{\theta_2\theta_3}{\theta_1^2} {\sss \;\;\;}-
\mfrac{\pi\ri}{4}\,\vartheta_2^2\,\vartheta_4^2\!\cdot\!
\frac{\theta_3^2}{\theta_1^2}\,\theta_4+
\mfrac{\ri}{\pi}\Big\{\eta+\mfrac{\pi^2}{12}
\big(\vartheta_3^4+\vartheta_4^4\big) \Big\}\!\cdot\!\theta_4\,,
\;\,\,
$$\
$$
\frac{\partial \theta_1'}{\partial\tau}=
\mfrac{-\ri}{4\,\pi}\,\frac{\theta_1'^3}{\theta_1^2}
+\mfrac{3\,\ri}{\pi}\!\left\{
\mfrac{\pi^2}{4}\vartheta_3^2\,\vartheta_4^2\!\cdot\!
\frac{\theta_2^2}{\theta_1^2} +\eta+
\mfrac{\pi^2}{12}\big(\vartheta_3^4+\vartheta_4^4\big)
\right\}\theta_1'
-\mfrac{\pi^2}{2}\ri\,\vartheta_2^2\,\vartheta_3^2\,\vartheta_4^2\!\cdot\!
\frac{\theta_2\theta_3\theta_4}{\theta_1^2}\,.
$$
Общая запись правил $\tau$-дифференцирования $\theta$-функций имеет
вид:
\begin{equation}\label{Dtau}
\left\{
\begin{array}{rl}
\ds \frac{\partial\theta_k}{\partial \tau}\!\!\!&=\ds
\mfrac{-\ri}{4\,\pi}\, \frac{\theta_1'^2}{\theta_1^2} \,\theta_k+
\mfrac{\ri}{2}\, \vartheta_k^2\!\cdot\!\theta_1'\,
\frac{\theta_\nu\,\theta_\mu}{\theta_1^2}+ \mfrac{\pi\ri}{4}\,\Big\{
\vartheta_3^2\,\vartheta_4^2\!\cdot\!\theta_2^2-
\vartheta_k^2\,\vartheta_\mu^2\!\cdot\!\theta_\nu^2-
\vartheta_k^2\,\vartheta_\nu^2\!\cdot\!\theta_\mu^2
\Big\}\,\frac{\theta_k}{\theta_1^2}
+{}\\ \\
&\ds\phantom{=\;} {}+\mfrac{\ri}{\pi}\,\Big\{ \eta+\mfrac{\pi^2}{12}
\big(\vartheta_3^4+ \vartheta_4^4\big)\Big\}\!\cdot\!
\theta_k\,,\qquad\quad
\mbox{где}\quad\nu=\mfrac{8\,k-28}{3\,k-10}\,,\;\;
\mu=\mfrac{10\,k-28}{3\,k-8}\\\\
\ds \frac{\partial \theta_1'}{\partial\tau}\!\!\!&=\ds
\mfrac{-\ri}{4\,\pi}\,\frac{\theta_1'^3}{\theta_1^2}
+\mfrac{3\,\ri}{\pi}\!\left\{
\mfrac{\pi^2}{4}\vartheta_3^2\,\vartheta_4^2\!\cdot\!
\frac{\theta_2^2}{\theta_1^2} +\eta+
\mfrac{\pi^2}{12}\big(\vartheta_3^4+\vartheta_4^4\big)
\right\}\theta_1'
-\mfrac{\pi^2}{2}\ri\,\vartheta_2^2\,\vartheta_3^2\,
\vartheta_4^2\!\cdot\!
\frac{\theta_2\theta_3\theta_4}{\theta_1^2}
\end{array}\right..
\end{equation}
Не лишним будет отметить, что эти дифференцирования не вполне
симметричны.

\section{Дифференциальное исчисление  $\vartheta,\eta$-констант}

\noindent Дифференциальную замкнутость образуют константы
$\vartheta_{2,3,4}$, дополненные $\eta$:
\begin{equation}\label{var}
\begin{array}{l}
\ds \frac{d\vartheta_2}{d\tau} _{\ds\mathstrut}=
\mfrac{\ri}{\pi}\,\Big\{\eta+\mfrac{\pi^2}{12}\,
\big(\vartheta_3^4+\vartheta_4^4 \big)\Big\}\,\vartheta_2\,, \qquad
\ds \frac{d\vartheta_4}{d\tau} _{\ds\mathstrut}=
\mfrac{\ri}{\pi}\,\Big\{\eta-\mfrac{\pi^2}{12}\,
\big(\vartheta_2^4+\vartheta_3^4 \big)\Big\}\,\vartheta_4\,,
\\\\
\ds \frac{d\vartheta_3}{d\tau} _{\ds\mathstrut}=
\mfrac{\ri}{\pi}\,\Big\{\eta+\mfrac{\pi^2}{12}\,
\big(\vartheta_2^4-\vartheta_4^4 \big)\Big\}\,\vartheta_3\,,
\qquad\;\;
\frac{d\eta}{d\tau}=\mfrac{\ri}{\pi}\,\Big\{2\,\eta^2-\mfrac{\pi^4}{12^2}
\big(\vartheta_2^8+\vartheta_3^8+\vartheta_4^8 \big) \Big\}\,.
\end{array}
\end{equation}
Дифференциальные соотношения на логарифмы отношений якобиевских
рядов $\vartheta_2\!:\!\vartheta_3\!:\!\vartheta_4$, а также
известная динамическая система Альфена и ее многочисленные
разновидности часто встречаются в  литературе. Они порождаются
системой уравнений
\begin{equation}\label{g2g3}
\frac{d g_2^{}}{d\tau} = \frac{\ri}{\pi}
\Big(8\,g_2^{}\,\eta-12\,g_3^{}\Big)\,,\quad \frac{d g_3^{}}{d\tau}=
\frac{\ri}{\pi} \Big(12\,g_3^{}\,\eta-\mfrac23\,g_2^2\Big)\,,\quad
\frac{d\eta}{d\tau}=
\frac{\ri}{\pi}\Big(2\,\eta^2-\mfrac16\,g_2^{}\Big)\,,
\end{equation}
которая неявно появлялась у Вейерштрасса \cite[{\bf
II}:\,стр.\,249]{we} и выписывалась Альфеном \cite[{\bf
I}:\,стр.\,331, 449--50]{halphen}. Эквивалент системы уравнений
\reff{g2g3} получал также Рамануджан для своих
тео\-ре\-ти\-ко-чис\-ло\-вых $q$-рядов. См. например работу
\cite[\S\,1]{zudilin} и дополнительную информацию в ней. Иными
словами, для дифференциального замыкания необходимо расширить
систему Альфена \reff{g2g3} до ее четно-мерной версии \reff{var}.

Добавим сюда, что на последнее уравнение системы  \reff{g2g3} можно
смотреть как на уравнение Риккати и, следовательно, получить
нетривиальный пример соответствующего решаемого линейного
дифференциального уравнения 2-го порядка.  Коэффициент этого
уравнения пропорционален $g_2^{}(\tau)$ и является всюду голоморфной
в $\mathbb{H}^+$ автоморфной формой относительно группы
$\boldsymbol{\Gamma}(1)$. Выполняя вычисления, указанному уравнению
и его решению можно придать следующий вид \cite{br}:
$$
\Psi_{\tau\tau}+\frac{g_2^{}(\tau)}{3\pi^2}\,\Psi=0\,,\qquad
\Psi=\frac{1}{\ded^2(\tau)}\Big(A+B
\mbox{\footnotesize$\ds\int\limits^{\,\,\tau}$}
\!\ded^4(\tau)\,d\tau\Big)\,.
$$

Возвращаясь к системе \reff{var} заметим, что вероятно неизвестным
является то, что ее другие версии уже выписывались Якоби в связи с
разложениями $\theta$-функций в степенные ряды (упоминания этого нет
также у Альфена). Эти результаты были опубликованы Борхардтом на
основе бумаг оставшихся после Якоби \cite[{\bf
II}:\,стр.\,383--98]{jacobi}.

Якоби получил красивую динамическую систему на четыре функции
$(A,B,a,b)$ \cite[{\bf II}:\,стр.\,386]{jacobi}, которые, в наших
обозначениях, являются рациональными функциями от квадратов
$\vartheta$-констант:
$$
A=\vartheta_3^2\,,\quad B=\frac{4}{\pi^2\,\vartheta_3^2}\left\{\eta
+\frac{\pi^2}{12}\,(\vartheta_2^4-\vartheta_4^4)\right\}\,, \quad
a=4-8\,\frac{\vartheta_2^2}{\vartheta_3^2}\,,\quad
b=2\,\frac{\vartheta_2^2\,\vartheta_4^2}{\vartheta_3^4}
$$
и показал, что с помощью нее можно получать степенные разложения
$\theta$-функций.  (Системы \reff{X} и \reff{Dtau} имеют своими
коэффициентами тоже квадраты $\vartheta$-констант). Он также
заметил, что разложения получаются простыми и имеют рекурсивный вид,
если выделить экспоненциальный множитель $\re^{-\frac12 ABx^2}$
\cite[{\bf II}:\,стр.\,390]{jacobi}. Фактически
--- это  рекурренции типа Альфена
\reff{D} для $\sigma$-функции Вейерштрасса. Якоби также предъявил
два канонических преобразования переменных $(A,B,a,b)$, сохраняющие
вид уравнений.

Переход между переменными в вышеупомянутых системах не всегда
взаимно однозначен, но всегда алгебраичен. Все они являются
следствиями {\em частных\/} интегралов уравнений \reff{var}, потому
что переменные, фигурирующие в этих системах, однозначно рационально
выражаются через якобиевские константы $\vartheta, \eta$.
Родственная динамическая система возникала у Якоби еще раньше, когда
он получал свое известное уравнение на $\vartheta$-ряды \cite[{\bf
II:} стр.\,176]{jacobi} (обозначения Якоби)
\begin{equation}\label{Jacobi}
C^4\big(\!\ln
C^3C_{\tau\tau}\big)_\tau^2=16\,C^3C_{\tau\tau}-\pi^2\,,\qquad
C=\vartheta^{-2}\,.
\end{equation}
В свою очередь, в качестве следствия из \reff{Jacobi}, не трудно
показать, что логарифмические производные от $\vartheta$-рядов тоже
удовлетворяют компактному дифференциальному уравнению, которое
появится у нас в \S\,9.4. Это уравнение и его общее решение выглядят
следующим образом:
\begin{equation}\label{halphen}
\begin{array}{c}
\ds \big(X_\tau-2\,X^2 \big)X_{\tau\tau\tau}
-X_{\tau\tau}{}^{\!\!\!\!\!\!2}\;\; +16\,X^3 X_{\tau\tau}+4\,
\big(X_\tau-6\,X^2 \big)X_\tau{}^{\!\!\!2}=0\,,\\
\ds X=\frac{d}{d\tau}\ln\frac{\vartheta_k
\Big(\mfrac{a\,\tau+b}{c\,\tau+1}\Big)}{\sqrt{c\,\tau+1}}
^{\ds{}^{\ds\mathstrut}}\,.
\end{array}
\end{equation}
Общеизвестно,  что дифференциальные выражения типа
$$
\eta=-\pi\,\ri\,\frac{d}{d\tau}\ln\vartheta_2^{}
-\frac{\pi^2}{12}\big(\vartheta_3^4+\vartheta_4^4 \big)
$$
удовлетворяют уравнению Шази
\begin{equation}\label{chazy}
\pi\,\eta_{\tau\tau\tau}^{}=
12\,\ri\,\big(2\,\eta\,\eta_{\tau\tau}^{}-3\,\eta_\tau^2\big)\,.
\end{equation}
В заключение параграфа заметим, что  формулы для многократного
дифференцирования $\vartheta,\eta$-конс\-тант явно  даются как
коэффициенты рядов \reff{t1} и (\ref{series}--\ref{t234}). Эти же
коэффициенты доставляют общие выражения для величин
$\theta^{(n)}(0|\tau)$. Заметим, что соотношения между
$\theta^{(n)}(0|\tau)$ при малых $n$ часто используются в литературе
как вспомогательные тождества
\cite{bateman,tannery,weber,WW,krause}.

\section{$\theta$-функции с характеристиками}

\noindent Предыдущие результаты, используя нотацию с
характеристиками, позволяют не только унифицировать формулы, но и
могут служить предметом для дальнейших обобщений на  высшие род\'а.

\subsection{$(\alpha,\beta)$-представления}

В  $(\alpha,\beta)$-представлении может быть записан любой объект,
симметричный по $\vartheta$-константам. Например, представление
точек ветвления через $\vartheta$-константы (\ref{ee}) или
$(\alpha,\beta)$-представления для инвариантов $g_{2,3}^{}$:
$$
\begin{array}{l}
\ds g_2^{}(\tau)=\phantom{1} \frac{\pi^4}{12}\, \bbig\{
\vartheta_{\sss\alpha0}{}^{\!\!\!\!\!\!8}\;+ \spin{\alpha+\beta}\,
\vartheta_{\sss\alpha0}{}^{\!\!\!\!\!\!4}\;\;
\vartheta_{\sss0\beta}{}^{\!\!\!\!\!\!4}\;+
\vartheta_{\sss0\beta}{}^{\!\!\!\!\!\!8}\;
\bbig\}\;, \qquad\quad (\alpha,\beta)\ne(0,0)\\ \\
\ds g_3^{}(\tau)=\frac{\pi^6}{432}\,\Big\{
2\,\spin{\beta}\,\vartheta_{\sss\alpha0}{}^{\!\!\!\!\!\!12} -
3\,\vartheta_{\sss\alpha0}{}^{\!\!\!\!\!\!4} \;\;
\vartheta_{\sss0\beta}{}^{\!\!\!\!\!\!4}\;\, \big(\spin{\beta}\,
\vartheta_{\sss0\beta}{}^{\!\!\!\!\!\!4}\; -\spin{\alpha}\,
\vartheta_{\sss\alpha0}{}^{\!\!\!\!\!\!4} \;
 \big)
-2\,\spin{\alpha}\,\vartheta_{\sss0\beta}{}^{\!\!\!\!\!\!12}
 \Big\}\;.
\end{array}
$$
Другой пример --- тождество
$\vartheta_3^4=\vartheta_2^4+\vartheta_4^4$ и формула Якоби
$\vartheta_1'=2\pi\ded^3$, которые в $(\alpha,\beta)$-представлении
имеют следующий вид:
\begin{equation}\label{jacobi}
\begin{array}{rl}
\varthetaAB{\alpha}{\beta}^4\!\!\!\!&= \bbig(
\spin{\beta}\,\varthetaAB{\alpha\mbox{\tiny--}1}{0}^4+
\spin{\alpha}\,\varthetaAB{0}{\beta\mbox{\tiny--}1}^4\bbig)
\,\mfrac{\spin{\alpha\beta}+1}{2}\\ \\
\vartheta_{{\sss\!}\alpha\beta}{}{\!\!\!\!\!'}\;\:(\tau)\!\!\!\!&=
\ri^{\beta+1}\,(1-\spin{\alpha\beta})\!\cdot\!
\pi\,\ded^3(\tau)\;,\qquad \mbox{где\ \ }
\vartheta_{{\sss\!}\alpha\beta}{}{\!\!\!\!\!'}\;\:(\tau)\equiv
\theta_{{\sss\!}\alpha\beta}{}{\!\!\!\!\!'}\;\:(0|\tau)\,.
\end{array}
\end{equation}

\subsection{Преобразования производных при сдвигах}
В самом общем случае, когда $\vartheta'$-константа есть
значение производной от $\theta$-функции в некотором
$\frac12$-периоде, формула Якоби (\ref{jacobi})  обобщается в
следующее выражение:
\begin{equation}\label{jac}
\theta_{{\sss\!}\alpha\beta}{}\!\!\!\!\!'\;\, \Big(\mfrac n2+\mfrac
m2\,\tau \Big|\tau\Big)= \ri^{\sss 1-(\beta+n)m}\!\cdot\! \pi\,
\bbig\{
\ri^{\sss\beta+n}\big(1-\big\langle(\alpha+m)(\beta+n)\big\rangle
\big) \!\cdot\!\ded^3-m\,
\varthetaAB{\alpha\mbox{\tiny+}m}{\beta\mbox{\tiny+}n} \bbig\}
\,\re^{\!-\frac {\pi\ri}{4}m^2\tau}_{\mathstrut}\,.
\end{equation}
В этой формуле, в силу четности/нечетности характеристики
$\big\langle(\alpha+m)(\beta+n)\big\rangle$, всегда остается
слагаемое содержащее либо $\ded^3$, либо
$\varthetaAB{\alpha\mbox{\tiny+}m}{\beta\mbox{\tiny+}n}$.

Алгебраическая и дифференциальная замкнутость функций Якоби с
целочисленными характеристиками ведет к закону преобразования
производных от $\theta$-функций при сдвигах аргумента. А именно,
любая функция $\theta_{{\sss\!}\alpha\beta}\!\!\!\!\!'\;\;(x|\tau)$,
с аргументами, сдвинутыми на произвольные $\frac12$-периоды,
выражается через функцию $\theta_1'(x|\tau)$  и все другие функции
$\theta_{1,2,3,4}(x|\tau)$. Окончательный ответ не является простым
дифференциальным следствием формулы \reff{shift} и его следует
рассматривать как дополнение к ней:
$$
\begin{array}{rl}
\theta_{{\sss\!}\alpha\beta}{}\!\!\!\!\!'\;\, \Big(x+\mfrac
n2+\mfrac m2\,\tau \Big|\tau\Big)\!\!\!\!&= \ri^{\sss
3m(\beta+n)}_{\mathstrut}\!\cdot\! \re^{\!-\frac{\pi\ri}{4}
m(4x+m\tau)}_{\mathstrut} \,\Big\{ \big(
\theta_1'-\pi\ri\,m\,\theta_1^{}\big)\,
\thetaAB{\alpha\mbox{\tiny+}m}{\beta\mbox{\tiny+}n}-{}\\
\\& \ds\phantom{=}
{}-\big\langle(\alpha+m)\big[\textstyle\frac{\beta+n}{2}\big]
\big\rangle \,\pi\,
\varthetaAB{\alpha\mbox{\tiny+}m}{\beta\mbox{\tiny+}n}^2
\!\cdot\!\thetaAB{1\mbox{\tiny--}\alpha\mbox{\tiny--}m} {0}\,
\thetaAB{0}{1\mbox{\tiny--}\beta\mbox{\tiny--}n}
\Big\}\,\theta_1^{-1}\;.
\end{array}
$$
При $x=0$ надо воспользоваться предыдущей формулой (\ref{jac}).

\subsection{Модулярное преобразование функций $\theta_1$ и $\ded$}
Помимо теоретических рассмотрений, модулярные преобразования
являются необходимыми при вычислении $\theta,\ded$-функций, хотя
применительно к свойствам $\theta$-функций всегда указывается ``с
точностью до некоторого корня восьмой степени из единицы'', а полный
алгоритм его вычисления оставляется в стороне. Между тем из формы
$\theta$-ряда видно, что без такого алгоритма ряд, который может
быть приведен к <<гиперсходящейся>> форме, может оказаться просто не
вычислимым.

Общее $\mathrm{PSL}_2(\mathbb{Z})$-пре\-обра\-зо\-ва\-ние для
функции $\theta_1$ замкнуто в себе, что отличает ее от функций
$\theta_{2,3,4}$ (\S\,8.4). Мы приведем самодостаточную формулу,
которая, по всей видимости, так и не была представлена в литературе.
Без потери общности всегда можно считать, что в матрице
$\left(\begin{smallmatrix}a&b\\c&d
\end{smallmatrix} \right)\in\mathrm{PSL}_2(\mathbb{Z})$ величина $c>0$.
Мы тогда имеем:
$$
\left\{
\begin{array}{l}
\ds \theta_1\!\Big(\mfrac{x}{c\,\tau+d}\Big|
\mfrac{a\,\tau+b}{c\,\tau+d}\Big)= \Kappa
\!\cdot\!\sqrt{c\,\tau+d\,}\,\, \mbox{\large
$\re$}^{\frac{\scriptstyle\pi\ri \,c\, x^2} {\scriptstyle
c\,\tau+d}}_{\mathstrut}\, \theta_1(x|\tau)\\\\
\ds \Kappa= \mbox{\large $\re$}_{\mathstrut}^{3\pi\ri\,\left\{
\phantom{\sum\limits_k^{m^2}}\right.
\!\!\!\!\!\!\!\!\!\frac{a-d}{12\,c}
-\frac{d}{6}(2c-3)+\frac{c-1}{4}\mathfrak{sign}(d)-\frac14
+\frac1c\, {\textstyle\sum\limits_ {\sss k=1}^{\sss c-1}}\,k
\left[\frac{d}{c}k\right]\!\!\!\!\!\!\!\!\!
\left.\phantom{\sum\limits_k^{m^2}}\right\}}\,,\quad\Kappa^8=1\\\\
\theta_1(x|\tau+N)= \re_{\mathstrut}^{\frac{\pi\ri}{4}\mbox{\tiny
$N$}} \!\!\cdot\!\theta_1(x|\tau)
\end{array}
\right..
$$

Известно, что множитель $\Kappa$ может быть выражен через символ
Якоби $\big(\frac ab\big)$ (Эрмит \cite[{\bf
I}:\,стр.\,482--86]{hermite}; см. также
\cite[стр.\,183--93]{krazer}, \cite[{\bf II}:\,стр.\,57--8]{koenig},
\cite[стр.\,124--32]{weber}), для которого существуют
самостоятельные правила вычисления \cite{tannery}. Предложенные
формулы содержат упрощения сумм Дедекинда (о них см.
\cite{apostol}), а для самой функции $\ded$ мы имеем:
$$
\left\{
\begin{array}{rl}
\ds \ded\!\left(\mfrac{a\,\tau+b}{c\,\tau+d}\right)\!\!&=\;
\mbox{\large $\re$}_{\mathstrut}^{\pi\ri\, \left\{
\phantom{\sum\limits_k^{m^2}}\right.
\!\!\!\!\!\!\!\!\!\frac{a-d}{12\,c}
-\frac{d}{6}(2c-3)+\frac{c-1}{4}\mathfrak{sign}
(d)-\frac14+\frac1c\, {\textstyle\sum\limits_ {\sss k=1}^{\sss
c-1}}\,k \left[\frac{d}{c}k\right]\!\!\!\!\!\!\!\!\!
\left.\phantom{\sum\limits_k^{m^2}}\right\}}\,
\sqrt{c\,\tau+d\,}\,\,\ded(\tau)\\ \\
\ded(\tau+N)\!\!&=\;\mbox{\large $\re$}^{\frac{\pi\ri}{12}
\mbox{\tiny$N$}}_{\ds\mathstrut}\,\ded(\tau)
\end{array}\right.\;.
$$
Модулярное преобразование для $\theta_1'$, как одной из базисных
функций, получается взятием производной от  формулы преобразования
для $\theta_1$.

\subsection{Общие модулярные преобразования}
Важным моментом  здесь  является тот факт, что  такие преобразования
можно рассматривать и выводить не как сторонние свойства рядов, а
как следствие дифференциальных уравнений (\ref{X}--\ref{var}).
Точнее, эти уравнения допускают автоморфизмы, а их явный вид
определяется через некоторый дроб\-но-линейный анзац. То что он
дроб\-но-линейный,  можно определить из непрерывных симметрий
уравнений. Это даже не использует тот факт, что решения выражаются
через $\theta$-ряды (см. \S\,9). Более того, само наличие
дискретного автоморфизма следует из того, что каждая функция
$\theta_k$ удовлетворяет одному и тому же {\em обыкновенному\/}
дифференциальному уравнению (\S\,9.3). Аналогично
--- для функций $\vartheta_k$. Таким образом
сразу находятся два преобразования,  порождающие модулярную группу:
$\tau\mapsto\tau+1$ и $\tau\mapsto -1/\tau$. Мы опускаем подробности
к этим утверждениям, которые могут быть частично извлечены из
процедуры интегрирования в \S\,9.5.

Общий ответ, как и первый пример нетривиального использования
модулярных формул, был дан Эрмитом при получении его знаменитого
решения уравнения $x^5-x-A=0$\footnote{Решение этого уравнения
фактически выражается в $\vartheta$-константах и было получено
Эрмитом в 1858 г. См. последнюю формулу на стр.\,10 в
\cite{hermite}. Отметим неточность в этой формуле, которая,
насколько нам известно, повторяется всюду где решение
воспроизводится. Второй знак ``$-$'' там следует заменить на
``$+$''.}. Множитель $\Kappa$= $\Kappa(a,c,d)$ (см. \S\,8.3) был
записан им в виде квадратичных гауссовых сумм экспонент \cite[{\bf
I}:\,стр.\,482--86]{hermite}. Мы предлагаем здесь формулы через
экспоненту от суммы, что проще.

Общее модулярное преобразование $\tau \mapsto
\frac{a\,\tau+b}{c\,\tau+d}\in \mathrm{PSL}_2(\mathbb{Z}) $ для
$\theta$-функций Якоби определяется следующим выражением:
$$
\left\{
\begin{array}{l}
\ds \theta \raisebox{0.04em}{\mbox{$ \scalebox{0.5}[0.9]{\big[}
\!\!\raisebox{0.1em} {\scalebox{0.8}{\mbox{\tiny$
\begin{array}{c}\tilde\alpha\mbox{\tiny--}1\\\tilde\beta\mbox{\tiny--}1
\end{array}$}}}
\!\! \scalebox{0.5}[0.9]{\big]}$}}
\Big(\mfrac{x}{c\,\tau+d}\Big|\mfrac{a\,\tau+b}{c\,\tau+d}\Big)=
\Kappa\, \mbox{\large$\re$}_{\mathstrut}^ {\frac{\pi\ri}{4}\left\{
2\alpha(bc\,\beta-d+1)-c\beta(a\beta-2)-db\,\alpha^2 \right\}
}\!\cdot\! \sqrt{c\,\tau+d\,}\,\,
\mbox{\large$\re$}^{\frac{\pi\ri\,c\,x^2}{c\,\tau+d}}_{\mathstrut}\,
\theta \raisebox{0.04em}{\mbox{$ \scalebox{0.5}[0.9]{\big[}
\!\!\raisebox{0.1em} {\scalebox{0.8}{\mbox{\tiny$
\begin{array}{c}\alpha\mbox{\tiny--}1\\\beta\mbox{\tiny--}1
\end{array}$}}}
\!\! \scalebox{0.5}[0.9]{\big]}$}}
(x|\tau)\\\\
\phantom{Mmx}\theta \raisebox{0.04em}{\mbox{$
\scalebox{0.5}[0.9]{\big[} \!\!\raisebox{0.1em}
{\scalebox{0.8}{\mbox{\tiny$
\begin{array}{c}\alpha\mbox{\tiny--}1\\\beta\end{array}$}}}
\!\! \scalebox{0.5}[0.9]{\big]}$}} (x|\tau+N)=
\re_{\mathstrut}^{\frac{\pi\ri}{4}(1\mbox{\tiny--}\alpha^2)
\mbox{\tiny $N$}} \!\cdot\! \theta \raisebox{0.04em}{\mbox{$
\scalebox{0.5}[0.9]{\big[} \!\!\raisebox{0.1em}
{\scalebox{0.8}{\mbox{\tiny$
\begin{array}{c}\alpha-1\\\beta\!+\!N\alpha\end{array}$}}}
\!\! \scalebox{0.5}[0.9]{\big]}$}}(x|\tau)
\end{array}
\right.,
$$

\noindent
где характеристики вычисляются друг через друга по формулам
$$
\left\{
\begin{array}{l}
\ds \tilde\alpha=\phantom{-}d\,\alpha-c\,\beta_{{}_{\ds\mathstrut}}\\
\ds \tilde\beta=-b\,\alpha +a\,\beta
\end{array}\right.,
\qquad\qquad \left\{
\begin{array}{l}
\ds
\alpha=a\,\tilde\alpha+c\,\tilde\beta_{{}_{\ds\mathstrut}}\\
\beta=b\,\tilde\alpha +d\,\tilde\beta
\end{array}\right.\;.
$$
Замкнутость относительно $\mathrm{PSL}_2(\mathbb{Z})$-преобразований
функции $\theta_1$ есть следствие того, что выписанные формулы
перехода между характеристиками $(\alpha,\beta)$ и
$(\tilde\alpha,\tilde\beta)$ являются однородными преобразованиями.
Для справок приведем раздельные формулы преобразований для функций
$\theta_{2,3,4}$:
$$
\begin{array}{l}
\ds \theta_2\Big(\mfrac{z}{c\,\tau+d}\Big|
\mfrac{a\,\tau+b}{c\,\tau+d}\Big)= \hspace{20.8mm}\Kappa\,
\mbox{\large$\re$}_{\mathstrut}^ {\frac{\pi\ri}{4}(2-d)c}\!\cdot\!
\sqrt{c\,\tau+d\,}\,\,
\mbox{\large$\re$}^{\frac{\pi\ri\,c\,z^2}{c\,\tau+d}}_{\mathstrut}\,
\theta \raisebox{0.04em}{\mbox{$ \scalebox{0.5}[0.9]{\big[}
\!\!\raisebox{0.1em} {\scalebox{0.8}{\mbox{\tiny$
\begin{array}{c}c\!-\!1\\d\!-\!1
\end{array}$}}}
\!\! \scalebox{0.5}[0.9]{\big]}$}} (z|\tau)\,,
\\\\
\ds \theta_3\Big(\mfrac{z}{c\,\tau+d}\Big|
\mfrac{a\,\tau+b}{c\,\tau+d}\Big)= \Kappa\,
\mbox{\large$\re$}_{\mathstrut}^ {\frac{\pi\ri}{4}\, \big\{
2(a+c-ad)- ab-cd\big\} }\!\cdot\! \sqrt{c\,\tau+d\,}\,\,
\mbox{\large$\re$}^{\frac{\pi\ri\,c\,z^2}{c\,\tau+d}}_{\mathstrut}\,
\theta \raisebox{0.04em}{\mbox{$ \scalebox{0.5}[0.9]{\big[}
\!\!\raisebox{0.1em} {\scalebox{0.8}{\mbox{\tiny$
\begin{array}{c}a\!+\!c\!-\!1\\b\!+\!d\!-\!1\end{array}$}}}
\!\! \scalebox{0.5}[0.9]{\big]}$}} (z|\tau)\,, \\\\
\ds \theta_4\Big(\mfrac{z}{c\,\tau+d}\Big|
\mfrac{a\,\tau+b}{c\,\tau+d}\Big)= \hspace{10.4mm}\Kappa\,
\mbox{\large$\re$}_{\mathstrut}^ {\frac{\pi\ri}{4}\,
\big\{2(a-1)-ab\big\} }\!\cdot\! \sqrt{c\,\tau+d\,}\,\,
\mbox{\large$\re$}^{\frac{\pi\ri\,c\,z^2}{c\,\tau+d}}_{\mathstrut}\,
\theta \raisebox{0.04em}{\mbox{$ \scalebox{0.5}[0.9]{\big[}
\!\!\raisebox{0.1em} {\scalebox{0.8}{\mbox{\tiny$
\begin{array}{c}a\!-\!1\\b\!-\!1\end{array}$}}}
\!\! \scalebox{0.5}[0.9]{\big]}$}} (z|\tau)\,.
\end{array}
$$
При $x=0$, предыдущие формулы превращаются в преобразование для
$\vartheta$-констант:
$$
\left\{
\begin{array}{l}
\ds \vartheta \raisebox{0.04em}{\mbox{$ \scalebox{0.5}[0.9]{\big[}
\!\!\raisebox{0.1em} {\scalebox{0.8}{\mbox{\tiny$
\begin{array}{c}\tilde\alpha\mbox{\tiny--}1\\\tilde\beta\mbox{\tiny--}1
\end{array}$}}}
\!\! \scalebox{0.5}[0.9]{\big]}$}}
\Big(\mfrac{a\,\tau+b}{c\,\tau+d}\Big)= \Kappa\,
\mbox{\large$\re$}_{\mathstrut}^ {\frac{\pi\ri}{4}\left\{
2\alpha(bc\,\beta-d+1)-c\beta(a\beta-2)-db\,\alpha^2 \right\}
}\!\cdot\! \sqrt{c\,\tau+d\,}\, \vartheta \raisebox{0.04em}{\mbox{$
\scalebox{0.5}[0.9]{\big[}  \!\!\raisebox{0.1em}
{\scalebox{0.8}{\mbox{\tiny$
\begin{array}{c}\alpha\mbox{\tiny--}1\\\beta\mbox{\tiny--}1
\end{array}$}}}
\!\! \scalebox{0.5}[0.9]{\big]}$}}
(\tau)\\\\
\:\,\vartheta \raisebox{0.04em}{\mbox{$ \scalebox{0.5}[0.9]{\big[}
\!\!\raisebox{0.1em} {\scalebox{0.8}{\mbox{\tiny$
\begin{array}{c}\alpha\mbox{\tiny--}1\\\beta\end{array}$}}}
\!\! \scalebox{0.5}[0.9]{\big]}$}} (\tau+N)=
\re_{\mathstrut}^{\frac{\pi\ri}{4}(1\mbox{\tiny--}\alpha^2)\mbox{\tiny
$N$}} \!\cdot\! \varthetaAB{\alpha-1}{\beta\!+\!N\alpha}(\tau)
\end{array}
\right..
$$

\noindent Дополнительное упрощение   происходит когда
$\left(\begin{smallmatrix}a&b\\c&d
\end{smallmatrix}\right)\in\boldsymbol{\Gamma}(2)$.
Полагая в этом случае $\left(\begin{smallmatrix}a&b\\c&d
\end{smallmatrix}\right)=\left(\begin{smallmatrix}2n+1&2m\\2p&2q+1
\end{smallmatrix}\right)$ с целыми $(n,m,p,q)$, выписанные преобразования
приобретают разделяющийся вид: $\theta_2\mapsto\theta_2,\ldots$
Нетрудно видеть, что общую формулу всегда можно свести к
преобразованиям раздельного типа $\theta_k\mapsto\theta_k$, если
добавить неоднородные преобразования
$x\mapsto\frac{x+s\,\tau+r}{c\,\tau+d}$.

Предложенные формулы справедливы для произвольных комплексных
характеристик, а в рассматриваемом нами случае целочисленных
$(\alpha,\beta)$, т.\,е. $(\alpha,\beta)=\{0,1\}\;\mbox{mod}\;2$,
формулы замыкаются, т.\,е. тоже выполняется равенство
$(\tilde\alpha,\tilde\beta)=\{0,1\}\;\mbox{mod}\;2$. Отношение любых
двух $\vartheta,\theta$-функций не содержит множитель
$\Kappa(a,c,d)$. Эрмит использовал этот факт, чтобы построить свои
известные функции $\varphi(\tau)$, $\psi(\tau)$, $\chi(\tau)$ и
таблицы преобразования между ними  \cite{hermite,tannery}.

\subsection{Теоремы умножения}

\noindent Формулы удвоений выписывал еще Якоби \cite[{\bf
I}:\,стр.\,510]{jacobi}. Общие умножения, разбивая формулы на
сложные серии четных/нечетных  $n$, с привлечением  функций sn, cn,
dn, рассматривали Кёнигсбергер \cite[стр.\,192, 201,
203--4]{koenig}, Вебер \cite[стр.\,190--5]{weber}, Краузе
\cite[стр.\,157, 173]{krause}, а также Вейерштрасс и другие. Отметим
здесь еще детерминантную формулу умножения
Киперта--Вей\-ер\-штрас\-са \cite[\S\,15]{we2}, \cite[стр.\,332]{WW}
для функции $\sigma(nx)$ через функции $\wp,\,\wpp$.

Пусть $n$ есть произвольная комплексная величина. Тогда функции
$\theta$  удовлетворяют комплексным формулам умножения,
определяемыми следующими рекурренциями:
$$
\left\{
\begin{array}{l}\ds
\theta_1^{}(2x)=2\,\theta_1(x)\,
\frac{\theta_2(x)\,\theta_3(x)\,\theta_4(x)}
{\vartheta_2\,\vartheta_3\,\vartheta_4}\\\\
\ds \theta_1^{}(nx)=\frac{\theta_3^2(n_1^{} x)\,\theta_2^2(x)-
\theta_2^2(n_1^{}x)\,\theta_3^2(x)} {\vartheta_4^2 \cdot
\theta_1^{}\big((n-2)x\big)}\,,\qquad\mbox{где}\;\; n_1^{}\equiv n-1
\\\\
\ds \thetaAB{\alpha+1}{\beta+1}(nx)=-
\frac{\langle\beta\rangle\,\thetaAB{\alpha}{0}^2(n_1^{}x)\,
\thetaAB{\alpha}{0}^2(x)+\langle\alpha\rangle\,
\thetaAB{0}{\beta}^2(n_1^{}x)\,\thetaAB{0}{\beta}^2(x)}
{\varthetaAB{\alpha+1}{\beta+1}^2 \cdot
\thetaAB{\alpha+1}{\beta+1}\big((n-2)x\big)}
\end{array}\right..
$$
Это не единственное представление.  Другие получаются, если
использовать квадратичные тождества между $\theta$-функциями. Если
$n$
--- целое число, тогда формулы замыкаются. В частности, замкнуты
умножения для функций $\theta_{2,3,4}$:
$$
\left\{
\begin{array}{l}
\ds \theta_2^{}(nx)=\frac{\theta_3^2(n_1^{}x)\,\theta_3^2(x)-
\theta_4^2(n_1^{}x)\,\theta_4^2(x)} {\vartheta_2^2 \cdot
\theta_2^{}\big((n-2)x\big)}\\
\ds
\theta_3^{}(nx)=\frac{\theta_2^2(n_1^{\ds\mathstrut}x)\,\theta_2^2(x)+
\theta_4^2(n_1^{}x)\,\theta_4^2(x)} {\vartheta_3^2 \cdot
\theta_3^{}\big((n-2)x\big)}\\
\ds
\theta_4^{}(nx)=\frac{\theta_3^2(n_1^{\ds\mathstrut}x)\,\theta_3^2(x)-
\theta_2^2(n_1^{}x)\,\theta_2^2(x)} {\vartheta_4^2 \cdot
\theta_4^{}\big((n-2)x\big)}
\end{array}\right.\,.
$$
Существуют нерекурсивные теоремы умножения, но это будут умножения с
коэффициентами, не принадлежащими полю констант
$\vartheta_{{\sss\!}\alpha\beta}(\tau)$. Этот факт непосредственно
следует из вейерштрассовских формул для функций
$\sigma,\sigmalambda(nx)$ \cite[{\bf V}:\,стр.\,225--6, 298]{we}.

Умножение для функции $\theta_1'$ получается взятием производной.
Формулы умножения будут справедливыми также и для неканонических
функций $\theta$ из \S\,9.

\subsection{Дифференциальные уравнения}

\noindent Следующая теорема полностью описывает дифференциальное
исчисление якобиевских $\theta$-функций.

Тэта-функции Якоби с произвольными целочисленными характеристиками
$(\alpha,\beta)$, как функции двух переменных $x$ и $\tau$,
удовлетворяют расщепленным и замкнутым системам обыкновенных
дифференциальных уравнений над дифференциальным полем
$\mathbb{C}_\partial(\eta,\vartheta^2)$ вейерштрассовской константы
$\eta$ и квадратов якобиевских $\vartheta$-констант:
\small 
\begin{equation}\label{x}
\left\{
\begin{array}{rl}
\ds \frac{\partial\theta_{\alpha\beta}}{\partial x}\!\!&=\ds\;
\frac{\theta_1'}{\theta_1^{}}_{{}_{\ds\mathstrut}}
\,\theta_{\alpha\beta}-
\big\langle\alpha\,[{\sss\frac{\beta}{2}}]\big\rangle\,
\pi\,\vartheta_{\alpha\beta}^2\!\cdot\!
\frac{\thetaAB{1\mbox{\tiny--}\alpha}{0}\,
\thetaAB{0}{1\mbox{\tiny--}\beta}}{\theta_1^{}}
\\  \ds
\frac{\partial\theta_1'}{\partial x}\!\!&=\ds\;
\frac{\theta_1'^2}{\theta_1}-\pi^2\, \vartheta_3^2\,\vartheta_4^2
\!\cdot\! \frac{\theta_2^2}{\theta_1}-
4\,\Big\{\eta+\frac{\pi^2}{12}\bbig(\vartheta_3^4+
\vartheta_4^4\bbig) \Big\} \!\cdot\!\theta_1
\end{array}\right.,
\end{equation}
\medskip
\begin{equation}\label{tau}
\left\{
\begin{array}{rl}
\ds \frac{\partial\theta_{\alpha\beta}}{\partial \tau}\!\!&=\ds\;
\frac{-\ri}{4\,\pi}\, \frac{\theta_1'^2}{\theta_1^2}
\,\theta_{\alpha\beta}+\frac{\ri}{2}\,
\big\langle\alpha\,[{\sss\frac{\beta}{2}}]\big\rangle\,
 \vartheta_{\alpha\beta}^2\!\cdot \theta_1'\,
\frac{\thetaAB{1\mbox{\tiny--}\alpha}{0}\,
\thetaAB{0}{1\mbox{\tiny--}\beta}}{\theta_1^2} _{{}_{\ds\mathstrut}}
+\frac{\ri}{\pi}\Big\{\eta+\frac{\pi^2}{12} \bbig(\vartheta_3^4+
\vartheta_4^4\bbig)\Big\}\!\cdot\! \theta_{\alpha\beta}+{}
\\  &\ds\phantom{=\;} {}+\frac{\pi\,\ri}{4}
\left\{\vartheta_3^2\,\vartheta_4^2 \!\cdot\! \theta_2^2-
\vartheta_{\alpha\beta}^2
\bbig(\varthetaAB{0}{1\mbox{\tiny--}\beta}^2\!\cdot\!
\thetaAB{1\mbox{\tiny--}\alpha}{0}^2 +
\varthetaAB{1\mbox{\tiny--}\alpha}{0}^2\!\cdot\!
\thetaAB{0}{1\mbox{\tiny--}\beta}^2
\bbig)\right\}\frac{\theta_{\alpha\beta}}{\theta_1^2}
_{{}_{\mathstrut}}\\ \\
\ds \frac{\partial \theta_1'}{\partial\tau}\!\!&\ds=\;
\frac{-\ri}{4\,\pi}\,\frac{\theta_1'^3}{\theta_1^2}
+\frac{3\,\ri}{\pi}\bbig{\Big\{}
\frac{\pi^2}{4}\,\vartheta_3^2\,\vartheta_4^2\!\cdot\!
\frac{\theta_2^2}{\theta_1^2} +\eta+
\frac{\pi^2}{12}\bbig(\vartheta_3^4+\vartheta_4^4\bbig)
\bbig{\Big\}}_{{}_{\ds\mathstrut}^{}}\theta_1' -
\frac{\pi^2}{2}\,\ri\, \vartheta_2^2\,\vartheta_3^2\,\vartheta_4^2
\!\cdot\! \frac{\theta_2\,\theta_3\,\theta_4} {\theta_1^2}
\end{array}\right.
\end{equation}
\normalsize
(если рассматривать только случаи
$(\alpha,\beta)=\{0,1\}$, тогда
$\big\langle\alpha\,[{\sss\frac{\beta}{2}}]\big\rangle=1$ всегда).

$(\alpha,\beta)$-пред\-став\-ле\-ние уравнений \reff{var} на
$\vartheta, \eta$-константы имеет вид:
\begin{equation}\label{last}
\left\{
\begin{array}{rl}
\ds \frac{d\varthetaAB{\alpha}{\beta}}{d\tau}\!\!\!&\ds=
\frac{\ri}{\pi} \Big\{ \eta+\frac{\pi^2}{12}\,
\bbig(\spin{\beta}\,\varthetaAB{1\mbox{\tiny--}\alpha}{0}^4-
\spin{\alpha}\,\varthetaAB{0}{1\mbox{\tiny--}\beta}^4\bbig) \Big\}\,
\varthetaAB{\alpha}{\beta}\\
\\ \ds \frac{d\eta}{d\tau}\!\!\!&\ds= \frac{\ri}{\pi}
\Big\{2\,\eta^2- \frac{\pi^4}{72} \bbig( \varthetaAB{\alpha}{0}^8+
\spin{\alpha+\beta}\, \varthetaAB{\alpha}{0}^4
\varthetaAB{0}{\beta}^4 +\varthetaAB{0}{\beta}^8\bbig)
\Big\}\;\;\leftarrow\;(\alpha,\beta)\ne (0,0)
\end{array}\right..
\end{equation}
Решение системы уравнений (\ref{x}--\ref{tau}), сохраняющее два
классических соотношения
$$
\vartheta_2^2\,\theta_4^2-\vartheta_4^2\,\theta_2^2=
\vartheta_3^2\,\theta_1^2\,,\qquad
\vartheta_2^2\,\theta_3^2-\vartheta_3^2\,\theta_2^2=
\vartheta_4^2\,\theta_1^2\,,
$$
содержит три произвольные константы $A,B,C$ и имеет вид
\begin{equation}\label{sol}
\begin{array}{rl}
\ds \thetaAB{\alpha}{\beta}&\!\!\!=C\,\re^{\pi\ri A
(2x+A\tau)}\!\cdot\!
\thetaAB{\alpha}{\beta}(x+A\,\tau+B|\tau)\,,\\\\
\ds \theta_1'&\!\!\!= C\,\re^{\pi\ri A
(2x+A\tau)}\!\cdot\!\bbig\{\theta_1'(x+A\,\tau+B|\tau) -
2\pi\ri\,A\,\thetaAB{1}{1}(x+A\,\tau+B|\tau)\bbig\}\,.
\end{array}
\end{equation}
(Эти функции удовлетворяют уравнению теплопроводности
$4\pi\ri\,F_\tau=F_{\mathit{xx}}$). Общее решение уравнений
\reff{last}, держащее тождество Якоби \reff{jacobi}, дается
формулами
$$
\varthetaAB{\alpha}{\beta}=\sqrt[-2]{c\,\tau+d\,}\cdot
\varthetaAB{\alpha}{\beta}\!\left(\mfrac{a\,\tau+b}{c\,\tau+d}\right),
\qquad
\eta=(c\,\tau+d)^{-2}\cdot\eta\!\left(\mfrac{a\,\tau+b}{c\,\tau+d}\right)+
\frac12\frac{\pi\ri \,c}{c\,\tau+d}
$$
с тремя  константами $a\,d-b\,c=1$. Правые части в этих решениях
есть $\theta,\vartheta,\eta$-ряды. По отдельности, решения уравнений
\reff{Jacobi} и \reff{chazy} хорошо известны (Якоби, Шази).

\subsection{Замечания и обобщения}
Важным следствием предыдущих результатов является то, что
$\theta$-решения уравнения теплопроводности, как уравнения в частных
производных, специфичны, поскольку, на самом деле, они являются
решениями разделяющихся, и следовательно, обыкновенных
дифференциальных уравнений. К ним, в свою очередь, сводятся многие
решения интегрируемых нелинейных уравнений в частных производных.
Все такие решения, нестационарные и многофазные, являются следствием
интегрируемости единственной (продублированной по числу фаз) системы
обыкновенных дифференциальных уравнений на пять якобиевских функций
$\theta,\,\theta_1'$.

Динамические системы Якоби--Альфена и их обобщения тоже являются
следствиями расщепляемости уравнений на $\theta$-функции. Правила
дифференцирования (\ref{x}--\ref{last}) генерируют упомянутые
динамические системы и их решения, а
$\eta,\vartheta,\theta$-переменные являются униформизирующими для
них. Высшие производные от $\theta$-функций и $\vartheta$-констант
есть снова $\eta,\vartheta,\theta,\theta_1'$-функции, причем
рациональные полиномы от них. Это, например, тривиализирует
вычисления и проясняет результаты недавней работы \cite{pavlov}, где
было найдено приложение уравнения Шази к системам гидродинамического
типа.

Еще раз отметим, что несомненную важность представляет обобщение
предыдущих результатов на $\Theta$-функции высших родов и их
константы \cite{krazer}. См. упомянутую  работу \cite{zudilin} и
ссылки в ней, где вопрос уже рассматривался применительно к
логарифмическим производным от
$\scalebox{1.7}[1.0]{${{\vartheta}}$}$-констант. Имея в виду
результаты \S\S\,6--7, становится ясным например, что
дифференциальная замкнутость с необходимостью требует присутствия в
формулах как мероморфных абелевых интегралов ($\Rightarrow$ функций
$\Theta'$), так и их периодов, т.\,е. аналогов констант $\eta$. В
частных случаях анализ вырождается.

Рассмотрим, например, ситуацию, когда якобиан допускает редукцию
двумерной \mbox{$\Theta$-функции} к функциям Якоби:
$$
\begin{array}{rcl}
\ds\Theta \raisebox{0.01em}{\mbox{$ \scalebox{0.5}[1.0]{\big[}
\!\!\raisebox{0.1em} {\scalebox{0.8}{\mbox{\tiny
$\begin{array}{c}\alpha_1^{}\,\alpha_2^{}
\\\beta_1^{}\,\beta_2^{}\end{array}$}}}
\!\! \scalebox{0.5}[1.0]{\big]}$}} \Big(
\begin{smallmatrix}\ds u-{\scriptstyle\frac14}\alpha_2^{}\\
\ds{}^{\mathstrut}
v-{\scriptstyle\frac14}\alpha_1^{}\end{smallmatrix}\Big|
\ds\begin{smallmatrix}\ds\tau&\frac12_{}\\\frac12\,&
\ds\mu\end{smallmatrix} \Big)\!\!\!&=& \ds\!\!\!
\frac12\,\re^{-\frac{\pi\ri}{4}\alpha_1^{}\alpha_2^{}}_{\mathstrut}
\left\{ \bbig[ \thetaAB{\alpha_1^{}}{\beta_1^{}} (u|\tau)+
\ri^{\alpha_1^{}}_{\mathstrut}\,
\thetaAB{\alpha_1^{}}{\beta_1^{}\!-\!1} (u|\tau)\bbig]\,
\thetaAB{\alpha_2^{}}{\beta_2^{}}
(v|\mu)+{}\right.\\\\
\!\!\!&& \ds\!\!\!\qquad\;\;\,{}+
\ri^{\alpha_2^{}}_{\mathstrut}\left. \bbig[
\thetaAB{\alpha_1^{}}{\beta_1^{}} (u|\tau)-
\ri^{\alpha_1^{}}_{\mathstrut}\,
\thetaAB{\alpha_1^{}}{\beta_1^{}\!-\!1} (u|\tau)\bbig]\,
\thetaAB{\alpha_2^{}}{\beta_2^{}\!-\!1} (v|\mu)\right\}.
\end{array}
$$

\noindent
Соответствующие десять
$\scalebox{1.7}[1.0]{${{\vartheta}}$} \raisebox{0.04em}{\mbox{$
\scalebox{0.5}[0.9]{\big[} \!\!\raisebox{0.1em}
{\scalebox{0.8}{\mbox{\tiny
$\begin{array}{c}\boldsymbol{\alpha}\\\boldsymbol{\beta}\end{array}$}}}
\!\! \scalebox{0.5}[0.9]{\big]}$}}$-констант, будут выражаться через
величины
$$\vartheta_{2,3,4}(\tau)\,,\quad\vartheta_{2,3,4}(\mu)\,,\quad
\theta_{1,2,3,4}\!\big({\textstyle\frac14}\big|\tau\big) \,,\quad
\theta_{1,2,3,4}\!\big({\textstyle\frac14}\big|\mu\big)\,.
$$
Таким образом, для замыкания, необходимо иметь связь между
<<$\frac14$-периодами>> и обычными $\vartheta$-константами. Это
осуществляют  соотношения
$$
\begin{array}{ll}
\ds \phantom{2\,}\theta_4^{}\!\big({\textstyle\frac14}\big)=
\theta_3^{}\!\big({\textstyle\frac14}\big)\,,\qquad&
\ds\phantom{2\,} \theta_2^{}\!\big({\textstyle\frac14}\big)=
\theta_1^{}\!\big({\textstyle\frac14}\big)\,,\\\\
 \ds
2\,\theta_3^4\!\big({\textstyle\frac14}\big)=
\vartheta_4^{}\vartheta_3^3+ \vartheta_3^{}\vartheta_4^3\,,\qquad &
2\,\theta_1^4\!\big({\textstyle\frac14}\big)=
\vartheta_4^{}\vartheta_3^3- \vartheta_3^{}\vartheta_4^3
\,,
\end{array}
$$
(более подробно новые семейства $\vartheta,\theta$-тождеств будут
рассмотрены отдельно), после чего можно проинспектировать систему
производных. В данном примере, в силу \reff{X} и уравнения
теплопроводности на $\Theta \raisebox{0.04em}{\mbox{$
\scalebox{0.5}[0.9]{\big[} \!\!\raisebox{0.1em}
{\scalebox{0.8}{\mbox{\tiny
$\begin{array}{c}\boldsymbol{\alpha}\\\boldsymbol{\beta}\end{array}$}}}
\!\! \scalebox{0.5}[0.9]{\big]}$}}$, мы получаем, что имеется
дифференциальная замкнутость в виде обыкновенных дифференциальных
уравнений:
$$\left\{
\frac{\partial \Theta \raisebox{0.04em}{\mbox{$
\scalebox{0.5}[0.9]{\big[} \!\!\raisebox{0.1em}
{\scalebox{0.8}{\mbox{\tiny
$\begin{array}{c}\boldsymbol{\alpha}\\\boldsymbol{\beta}\end{array}$}}}
\!\! \scalebox{0.5}[0.9]{\big]}$}}}{\partial (u,v)}
=F_{u,v}^{\alpha\beta} \big(\Theta,\Theta'\big),\qquad
\frac{\partial \Theta' \raisebox{0.04em}{\mbox{$
\scalebox{0.5}[0.9]{\big[} \!\!\raisebox{0.1em}
{\scalebox{0.8}{\mbox{\tiny
$\begin{array}{c}\boldsymbol{\alpha}\\\boldsymbol{\beta}\end{array}$}}}
\!\! \scalebox{0.5}[0.9]{\big]}$}}}{\partial (u,v)}
=G_{u,v}^{\alpha\beta} \big(\Theta,\Theta'\big)\right\},\qquad
$$
с вычисляемыми функциями $F(\Theta,\Theta')$, $G(\Theta,\Theta')$
(вычисления громоздки)  и вытекающими отсюда уравнениями по $\tau$ и
$\mu$. Отсюда же вытекают  условия интегрируемости на константы
$\scalebox{1.7}[1.0]{${{\vartheta}}$} \raisebox{0.04em}{\mbox{$
\scalebox{0.5}[0.9]{\big[} \!\!\raisebox{0.1em}
{\scalebox{0.8}{\mbox{\tiny
$\begin{array}{c}\boldsymbol{\alpha}\\\boldsymbol{\beta}\end{array}$}}}
\!\! \scalebox{0.5}[0.9]{\big]}$}}$$(\tau,\mu)$. Этот результат сам
по себе не безинтересен, но больший интерес представляет вопрос о
том, какой вид имеют базовые дифференциальные соотношения на общие
$g$-мерные $\Theta$-функции на якобианах, их сечениях на
алгебраические кривые посредством голоморфных интегралов
(отображение Абеля) и, как следствие, условия интегрируемости этих
уравнений.

\section{Неканонические $\theta$-функции Якоби}
\noindent Дифференциальные уравнения, определяющие якобиевские
функции $\theta$, имеют  более широкий класс решений,  чем просто
$\theta$-ряды или формулы \reff{sol}. Они также могут служить одним
из независимых источников происхождения $\vartheta,\theta$-функций
вообще, так как уравнения не менее фундаментальны чем их решения. В
прикладных вопросах разновидности системы \reff{X} могут возникать
самостоятельно, поэтому параметры $\vartheta$ не обязаны быть
значениями $\theta$-рядов в нуле, а параметр $\eta$ быть периодом
мероморфного   интеграла $\zeta$ или даваться какими-либо известными
выражениями, например
$$
\eta(\tau)=-\frac{1}{12}\,\frac{\theta_1'''(0|\tau)}{\theta_1'(0|\tau)}=
-\frac{1}{4}\,\frac{\theta_3''(0|\tau)}
{\vartheta_3}-\frac{\pi^2}{12}\big(\vartheta_2^4-\vartheta_4^4\big).
$$

\subsection{Определяющие уравнения и их первые интегралы}
Стартуя с  уравнений \reff{X} с произвольными параметрами
$\vartheta, \eta$ и уравнения теплопроводности
$4\pi\ri\,\theta_\tau=\theta_{\mathit{xx}}$,  мы получаем уравнения
\reff{Dtau}. Рассматривая условия интегрируемости
$\theta_{x\tau}=\theta_{\tau x}$ для систем уравнений \reff{X} и
\reff{Dtau}, мы получаем их первые интегралы
\begin{equation}\label{A12}
\vartheta_2^2\,\theta_4^2-\vartheta_4^2\,\theta_2^2=
A_1^4\!\cdot\!\vartheta_3^2\,\theta_1^2\,,\qquad
\vartheta_2^2\,\theta_3^2-\vartheta_3^2\,\theta_2^2=
A_2^4\!\cdot\!\vartheta_4^2\,\theta_1^2
\end{equation}
и дифференциальную замкнутость поля коэффициентов, т.\,е.
динамическую систему на функции
$\vartheta=\vartheta(\tau),\,\eta=\eta(\tau)$:
\begin{equation}\label{intA}
\left\{
\begin{array}{rl}
\ds \frac{d\vartheta_2}{d\tau} _{\ds\mathstrut}&\!\!\!=
\mfrac{\ri}{\pi}\,\Big\{\eta+\mfrac{\pi^2}{12}\,
\big(\vartheta_3^4+\vartheta_4^4 \big)\Big\}\,\vartheta_2^{}\\
\ds \frac{d\vartheta_3}{d\tau}
_{\ds\mathstrut}^{\ds\mathstrut}&\!\!\!=
\mfrac{\ri}{\pi}\,\Big\{\eta+\mfrac{\pi^2}{12}\,
\big(\vartheta_3^4+\vartheta_4^4 -3A_2^4\,\vartheta_4^4\big)\Big\}\,
\vartheta_3^{}\\
\ds \frac{d\vartheta_4}{d\tau} ^{\ds\mathstrut}&\!\!\!=
\mfrac{\ri}{\pi}\,\Big\{\eta+\mfrac{\pi^2}{12}\,
\big(\vartheta_3^4+\vartheta_4^4
-3A_1^4\,\vartheta_3^4\big)\Big\}\,\vartheta_4^{}
\\
\ds \frac{d\eta}{d\tau}&\!\!\!=\mfrac{\ri}{\pi}\,2\,\eta^2-
\mfrac{\pi^3}{72}\,\ri\,\Big\{\vartheta_3^8+
\big(9A_1^4A_2^4-6A_1^4-6A_2^4+2\big)\,\vartheta_3^4\,\vartheta_4^4+
\vartheta_4^8 \Big\}^{{}^{\ds\mathstrut}}
\end{array}\right.,
\end{equation}
где $A_1,\,A_2$ являются произвольными константами, не зависящими от
$x$ или $\tau$. Прямая проверка уравнений \reff{X}, \reff{Dtau}
показывает, что соотношения \reff{A12} действительно выполняются и
являются интегралами обеих систем. Системы полностью совместны, а
сами интегралы $A_{1,2}(\theta;\vartheta)$ являются {\em
алгебраическими\/}. В свою очередь система \reff{intA} имеет один,
тоже алгебраический, интеграл $A_3(\vartheta)$, обобщающий тождество
Якоби:
\begin{equation}\label{jacobiA}
A_3^4\,\vartheta_2^4=A_1^4\,\vartheta_3^4-A_2^4\,\vartheta_4^4
\qquad \Big(\Rightarrow\;\frac{d}{d\tau}A_3=0\Big).
\end{equation}

Мы  не останавливаемся отдельно на вырождениях частных решений
уравнений \reff{intA} в элементарные функции. Эти случаи выявляются
легко и поэтому нас будет интересовать только ситуация общего
положения. Для таких неканонических <<обобщений>> $\theta$-функций
величины $A_1,\,A_2$ определяются начальными данными к уравнениям
(\ref{X}--\ref{Dtau}), а для системы \reff{intA} они являются
параметрами.

\subsection{Канонические
$\theta$-ряды и эллиптические функции} Редукция к случаю
Яко\-би--Вей\-ер\-штрас\-са осуществляется следующим образом.
Положим $A_1=A_2=1$.  Тогда функция $\eta$  будет удовлетворять
уравнению Шази \reff{chazy}, а функции $\vartheta_{3,4}$ ---
якобиевскому уравнению \reff{Jacobi}. Интеграл \reff{jacobiA}
остается произвольным. Полагая дополнительно $A_3=1$, приходим к
уравнениям \reff{var}. Заметим, что для такой симметричной формы
канонического случая уравнений \reff{intA},  алгебраический интеграл
уравнений \reff{var} выглядит уже не очевидно
\begin{equation}\label{A3}
(A_3^4-1)\!\cdot\!\vartheta_2^4\,\vartheta_3^4\,\vartheta_4^4=
(\vartheta_3^4-\vartheta_2^4-\vartheta_4^4)^3
\end{equation}
и должен трактоваться как корректная форма <<полного тождества
Якоби>> для симметризованной записи определяющих уравнений. Помимо
этого усложнения, ни одна из функций $\vartheta_{2,3,4}$ или
логарифмических производных $\ln_{\tau}\!\vartheta$ не удовлетворяет
никакому уравнению 3-го порядка, как это было в уравнениях
(\ref{Jacobi}--\ref{halphen}), хотя функция $\eta$ по-прежнему
удовлетворяет уравнению 3-го порядка \reff{chazy}. За недостатком
места мы опускаем доказательство этих утверждений.

Обозначая далее $\boldsymbol{\mathrm{P}}=\theta_2^2/\theta_1^2$, мы
находим, что имеет место уравнение
\begin{equation}\label{wp}
 \boldsymbol{\mathrm{P}}_x{}^{\!\!\!2}=4\,\pi^2\,
\big(\vartheta_4^2\!\cdot\!\boldsymbol{\mathrm{P}}+A_1^4\vartheta_3^2\big)
\big(\vartheta_3^2\!\cdot\!\boldsymbol{\mathrm{P}}+A_2^4\vartheta_4^2\big)\,
\boldsymbol{\mathrm{P}}\,,
\end{equation}
а значит $\boldsymbol{\mathrm{P}}$ выражается через $\wp$-функцию
Вейерштрасса, которая, как известно, пропорциональна отношению
якобиевских $\theta$-рядов
$$
\wp(2x|\mu)=\frac{\pi^2}{12}\,
\big\{\vartheta_3^4(\mu)+\vartheta_4^4(\mu)\big\}+ \frac{\pi^2}{4}\,
\vartheta_3^2(\mu)\,\vartheta_4^2(\mu)\!\cdot\!
\frac{\theta_2^2(x|\mu)}{\theta_1^2(x|\mu)}\,.
$$
Здесь и далее под $\vartheta_{3,4}(\mu)$ понимаются
$\vartheta$-ряды. Приводя уравнение \reff{wp} к канонической форме
Вейерштрасса и прилагая стандартный вычислительный аппарат теории
\cite{WW}, мы получаем, что эллиптическим модулем для уравнения
\reff{wp} будет величина $\mu$, определяемая из следующего
трансцендентного уравнения:
\begin{equation}\label{J}
J(\mu)=\frac{1}{54}\frac{\big(A_1^8\,\vartheta_3^8+A_2^8\,\vartheta_4^8+
A_3^8\,\vartheta_2^8\big)^3}
{A_1^8\,A_2^8\,A_3^8\,\vartheta_2^8\,\vartheta_3^8\,\vartheta_4^8}\,,
\end{equation}
где, для симметрии, мы воспользовались интегралом \reff{jacobiA}. (В
явном виде решение  уравнения \reff{J} дается, что хорошо известно
\cite{WW,a}, как отношение определенных гипергеометрических рядов
\cite[стр.\,87--8]{bateman}). Таким образом, должна иметь место
следующая формула для отношения $\boldsymbol{\mathrm{P}}$:
$$
\frac{\theta_2^2}{\theta_1^2}=
\frac{\vartheta_3^4(\mu)+\vartheta_4^4(\mu)}
{3\,\vartheta_3^2\,\vartheta_4^2}
-\frac{A_1^4\vartheta_3^4+A_2^4\vartheta_4^4}
{3\,\vartheta_3^2\,\vartheta_4^2} +
\frac{\vartheta_3^2(\mu)\,\vartheta_4^2(\mu)}
{\vartheta_3^2\,\vartheta_4^2}\!\cdot\!
\frac{\theta_2^2(x+x_0^{}|\mu)}{\theta_1^2(x+x_0^{}|\mu)}\,.
$$
Аналогичные формулы  получаются для дробей $\theta_k/\theta_j$ без
квадратов. Нетрудно видеть, что такая ва\-ри\-а\-ция привела бы к
эллиптическим функциям Якоби $\mathrm{sn}\sim \theta_1/\theta_4$ и
т.\,д.:
$$
\left(\frac{\theta_1}{\theta_4}\right)_{\!x}^{\!2} =\pi^2
\left\{\!A_1^4\vartheta_3^2\!\cdot\!\Big(\frac{\theta_1}{\theta_4}\Big)^{\!2}
-\vartheta_2^2\right\}\!
\left\{\!A_3^4\vartheta_2^2\!\cdot\!
\Big(\frac{\theta_1}{\theta_4}\Big)^{\!2} -\vartheta_3^2 \right\}.
$$
Итак, отношение любых двух $\theta$-решений уравнений \reff{X},
\reff{Dtau} пропорционально отношению канонических $\theta$-рядов с
нетривиальным и новым модулем $\mu$. Покажем, что уравнения
интегрируются через канонические $\theta$-ряды полностью.

\subsection{Расширение канонических $\theta$-функций}
Система \reff{X} имеет пятый порядок, в то время как
вейерштрассовский базис $(\sigma,\zeta,\wp)$ --- третий.  В
вейерштрассовском случае мы имеем следующее легко выводимое
дифференциальное уравнение 3-го порядка на любую из функций Якоби
\begin{equation}\label{w}
\ds F_x{}^{\!\!2}= 4\,\bbig\{F+4\,\eta+
{\textstyle\frac{\pi^2}{3}}(\vartheta_3^4+\vartheta_4^4)\bbig\}
\bbig\{F+4\,\eta+
{\textstyle\frac{\pi^2}{3}}(\vartheta_2^4-\vartheta_4^4)\bbig\}
\bbig\{F+4\,\eta-
{\textstyle\frac{\pi^2}{3}}(\vartheta_2^4+\vartheta_3^4)\bbig\}\,,
\end{equation}
где $F=\ln_{\mathit{xx}}\! \theta_k(x|\tau)$, а
$\vartheta=\vartheta(\tau),\,\eta=\eta(\tau)$. В неканоническом
случае $A_1\ne 1\ne A_2$ это уравнение должно иметь свой аналог в
виде дифференциального уравнения 5-го порядка. Анализируя систему
\reff{X}, мы получаем, что таковым является уравнение\footnote{Это
важное уравнение, как уравнение 3-го порядка на функцию $F$,
является обобщением канонического уравнения Вейерштрасса
$F_x^2=4F^3-g_2^{}F-g_3^{}$, а точнее, его следствия
$F_{\mathit{xxx}}=12FF_x$. Это следствие {\em не является
редукцией\/} уравнения \reff{F}, хотя решением   \reff{F} тоже
является $\wp$-функция: $F=\wp(\om|\om,\omp)-\wp(x+c|\om,\omp)$.
Различие состоит в том, что до тех пор пока не затрагивается
дифференциальная замкнутость $\theta_k$, достаточно одного уравнения
Вейерштрасса \reff{w}. В обоих случаях периоды $2\,(\om,\omp)$ и
модуль $\tau$ возникают как константы интегрирования.}
\begin{equation}\label{F}
\begin{array}{c}
F^2F_{\mathit{xxx}}-2\,F F_x F_{\mathit{xx}}+F_x{}^{\!\!3}+
(F^4)_x=0\,,\\\\
F=(\ln\theta)_{\mathit{xx}}^{}+4\,\bbig\{\eta+\mfrac{\pi^2}{12}
(\vartheta_3^4+\vartheta_4^4)\bbig\}\,,
\end{array}
\end{equation}
которому удовлетворяет произвольное решение $\theta=\theta_k(x)$
уравнений \reff{X}, какие бы ни были величины $\eta,\vartheta$.
Опуская вычисления и вводя обозначение
$$
M(\varkappa,\mu)=\varkappa^2\bbig\{\eta(\mu)+ \mfrac{\pi^2}{12}
\big(\vartheta_3^4(\mu)+\vartheta_4^4(\mu)\big)\!\bbig\}
-\bbig\{\eta+\mfrac{\pi^2}{12} \big(\vartheta_3^4+\vartheta_4^4)
\bbig\}\,,
$$
выпишем  общий интеграл уравнений \reff{X}:
\begin{equation}\label{Xgeneral}
\begin{array}{lr}
\ds
\pm\,\theta_1=&\ds\!\!\!\frac{\vartheta_2\vartheta_3\vartheta_4}{2\,\ded^3(\mu)}
\cdot C\,\theta_1(\varkappa\,x+B|\,\mu)\,\re_{\ds\mathstrut}^{2M(x+A)^2}\,,\\
\ds \pm\,\theta_2=&\ds\!\!\!
\frac{\varkappa\,\vartheta_2}{\vartheta_2(\mu)}\cdot
C\,\theta_2(\varkappa\,x+B|\,\mu)
\,\re_{\ds\mathstrut}^{2M(x+A)^{2^{\ds\mathstrut}}} \,,\\
\ds
\pm\,\theta_3=&\ds\!\!\!\frac{\varkappa\,\vartheta_3}{\vartheta_3(\mu)}
\cdot C\,\theta_3(\varkappa\,x+B|\,\mu)
\,\re_{\ds\mathstrut}^{2M(x+A)^{2^{\ds\mathstrut}}} \,,\\
\ds
\pm\,\theta_4=&\ds\!\!\!\frac{\varkappa\,\vartheta_4}{\vartheta_4(\mu)}
\cdot C\,\theta_4(\varkappa\,x+B|\,\mu)
\,\re_{\mathstrut}^{2M(x+A)^{2^{\ds\mathstrut}}} \,,
\end{array}
\end{equation}
где $\{A,B,C,\varkappa,\mu\}$ --- пять  констант интегрирования.
Формула для $\theta_1'$ --- это производная от первой формулы в
\reff{Xgeneral}. Знаки $\pm$ можно произвольно менять у любой пары
$(\theta_j,\theta_k)$, что легко устанавливается непосредственно из
системы \reff{X}.

Решение \reff{Xgeneral} показывает, что зависимость от $\varkappa$ и
$\mu$ нетривиальна, в отличие от <<угадываемой>> зависимости по
константам $B,C$ и линейной экспоненты $\re^{Ax}$. Отсюда также
следует, что если оставить в стороне  мультипликативные константы в
\reff{Xgeneral}, то зависимость этих решений и самого уравнения
\reff{F} от параметров уравнений $\vartheta,\eta$ осуществляется
через один существенный параметр
\begin{equation}\label{param}
\frac14\,\Lambda=\eta+
\frac{\pi^2}{12}\big(\vartheta_3^4+\vartheta_4^4\big)\,.
\end{equation}
Алгебраические интегралы \reff{A12} в
$(\varkappa,\mu)$-представлении \reff{Xgeneral} принимают вид
следующих обобщений тождеств Якоби:
\begin{equation}\label{omega}
\vartheta_2^2\,\theta_4^2-\vartheta_4^2\,\theta_2^2=\varkappa^2
\frac{\vartheta_3^4(\mu)}{\vartheta_3^4}\!\cdot\!
\vartheta_3^2\,\theta_1^2\,, \qquad
\vartheta_2^2\,\theta_3^2-\vartheta_3^2\,\theta_2^2=\varkappa^2
\frac{\vartheta_4^4(\mu)}{\vartheta_4^4}\!\cdot\!
\vartheta_4^2\,\theta_1^2\,.
\end{equation}
При $\varkappa=1$ и таких $\vartheta,\eta$, что найдется такое
$\mu$, что $\vartheta=\vartheta(\mu)$ и $\eta=\eta(\mu)$, мы снова
приходим к случаю Якоби--Вейерштрасса.

Особо подчеркнем, что переход к неканоническому случаю, хотя и
осуществляется элементарной функцией $\re^{2Mx^2}$, не элементарно
зависит от параметров $(\varkappa,\mu)$ и является {\em
трансцендентным расширением\/}, так как канонические $\sigma$- и
$\theta$-функции определены с точностью до линейной экспоненты
$\re^{Ax}$ (см. формулы \reff{sol}). Экспоненциальный множитель
$\re^{2Mx^2}$ этого расширения тоже зависит от параметров
интегрирования $\varkappa$ и $\mu$. Свойства квазипериодичности для
неканонических $\theta$-расширений, т.\,е. аналоги формул
\reff{shift}, легко устанавливаются из самих решений \reff{Xgeneral}
и мы их не приводим.

\subsection{Перенормировка $\theta$-функций}
\noindent Решение \reff{Xgeneral} подсказывает  сделать
перенормировку $\theta\mapsto\boldsymbol\theta$ (она не
единственна):
$$
\boldsymbol\theta_1^{}=\theta_1^{}\,, \qquad
\boldsymbol\theta_2^{}=\pi\,\vartheta_3^{}\vartheta_4^{}
\!\cdot\!\theta_2^{}\,, \qquad
\boldsymbol\theta_3^{}=\pi\,\vartheta_2^{}\vartheta_4^{}
\!\cdot\!\theta_3^{}\,, \qquad
\boldsymbol\theta_4^{}=\pi\,\vartheta_2^{}\vartheta_3^{}
\!\cdot\!\theta_4^{}\,,
$$
после которой уравнения \reff{X} и \reff{F} станут содержать лишь
параметр \reff{param}:
\begin{equation}\label{Xbold}
\left\{\begin{array}{ll} \ds
\frac{\partial\boldsymbol\theta_2^{}}{\partial x}=
\frac{\boldsymbol\theta_1'}{\boldsymbol\theta_1^{}}\,\boldsymbol\theta_2^{}-
\frac{\boldsymbol\theta_3^{}\boldsymbol\theta_4^{}}{\boldsymbol\theta_1^{}}\,,
&\qquad \ds \frac{\partial\boldsymbol\theta_4^{}}{\partial x}=
\frac{\boldsymbol\theta_1'}{\boldsymbol\theta_1^{}}\,\boldsymbol\theta_4^{}-
\frac{\boldsymbol\theta_2^{}\boldsymbol\theta_3^{}}{\boldsymbol\theta_1^{}}
_{{}_{\ds\mathstrut}}\,, \qquad\quad
\frac{\partial\boldsymbol\theta_1^{}}{\partial
x}=\boldsymbol\theta_1'
\\
\ds\frac{\partial\boldsymbol\theta_3^{}}{\partial x}=
\frac{\boldsymbol\theta_1'}{\boldsymbol\theta_1^{}}\,\boldsymbol\theta_3^{}-
\frac{\boldsymbol\theta_2^{}\boldsymbol\theta_4^{}}{\boldsymbol\theta_1^{}}\,,&
\qquad\ds \frac{\partial\boldsymbol\theta_1'}{\partial x}=
\frac{\boldsymbol\theta_1'^2}{\boldsymbol\theta_1^{}}-
\frac{\boldsymbol\theta_2^2}{\boldsymbol\theta_1^{}}- \Lambda
\!\cdot\!\boldsymbol\theta_1^{}
\end{array}\right..
\end{equation}
Отсюда сразу следует, что если теперь привлечь  зависимость от
$\tau$, тогда $\tau$-диф\-фе\-рен\-ци\-ро\-ва\-ния $\theta$-функций
\reff{Dtau} (и вообще все уравнения) тоже значительно упрощаются.
Вводя для удобства  замену $\tau=4\pi\ri\,\boldsymbol\tau$, мы
получаем следующую систему уравнений:
\begin{equation}\label{Dtaubold}
\left\{
\begin{array}{rl}
\ds \frac{\partial\boldsymbol\theta_1}{\partial
\boldsymbol\tau}&\ds\!\!\!= \frac{\boldsymbol\theta_1'^2}
{\boldsymbol\theta_1^{}}-\frac{\boldsymbol\theta_2^2}
{\boldsymbol\theta_1^{}}-
\Lambda\!\cdot\!\boldsymbol\theta_1\,,\qquad

\frac{\partial\boldsymbol\theta_1'^{}}{\partial \boldsymbol\tau}=
\frac{\boldsymbol\theta_1'^3} {\boldsymbol\theta_1^2}-
3\big(\boldsymbol\theta_2^2+
\Lambda\!\cdot\!\boldsymbol\theta_1^2\big)
\,\frac{\boldsymbol\theta_1'}{\boldsymbol\theta_1^2}+
2\,\frac{\boldsymbol\theta_2^{}\boldsymbol\theta_3^{}
\boldsymbol\theta_4^{}} {\boldsymbol\theta_1^2}_{{}_{\ds\mathstrut}}
\\
\ds \frac{\partial\boldsymbol\theta_2}{\partial
\boldsymbol\tau}&\ds\!\!\!=
\frac{\boldsymbol\theta_1'^2}{\boldsymbol\theta_1^2}\,
\boldsymbol\theta_2-2\,\boldsymbol\theta_1'\,
\frac{\boldsymbol\theta_3^{}\boldsymbol\theta_4^{}}
{\boldsymbol\theta_1^2}-
\big(\boldsymbol\theta_2^2-\boldsymbol\theta_3^2-
\boldsymbol\theta_4^2\big)\,
\frac{\boldsymbol\theta_2^{}}{\boldsymbol\theta_1^2}
_{{}_{\ds\mathstrut}}-
\bbig\{\Lambda-\ln_{\boldsymbol\tau}(\vartheta_3\vartheta_4)
\bbig\}\!\cdot\!\boldsymbol\theta_2\\

\ds \frac{\partial\boldsymbol\theta_3}{\partial
\boldsymbol\tau}&\ds\!\!\!= \frac{\boldsymbol\theta_1'^2}
{\boldsymbol\theta_1^2}\,\boldsymbol\theta_3-
2\,\boldsymbol\theta_1'\,
\frac{\boldsymbol\theta_2^{}\boldsymbol\theta_4^{}}
{\boldsymbol\theta_1^2}+\boldsymbol\theta_4^2
\,\frac{\boldsymbol\theta_3^{}}{\boldsymbol\theta_1^2}
_{{}_{\ds\mathstrut}}-
\bbig\{\Lambda-\ln_{\boldsymbol\tau}(\vartheta_2\vartheta_4)\bbig\}\!\cdot\!
\boldsymbol\theta_3\\

\ds \frac{\partial\boldsymbol\theta_4}{\partial
\boldsymbol\tau}&\ds\!\!\!= \frac{\boldsymbol\theta_1'^2}
{\boldsymbol\theta_1^2}\,\boldsymbol\theta_4-
2\,\boldsymbol\theta_1'\,
\frac{\boldsymbol\theta_2^{}\boldsymbol\theta_3^{}}
{\boldsymbol\theta_1^2}+\boldsymbol\theta_3^2\,
\frac{\boldsymbol\theta_4^{}}{\boldsymbol\theta_1^2}{\sss\,}-
\bbig\{\Lambda-\ln_{\boldsymbol\tau}(\vartheta_2\vartheta_3)\bbig\}\!\cdot\!
\boldsymbol\theta_4
\end{array}\right..
\end{equation}
Через каноническую форму (\ref{Xbold}--\ref{Dtaubold}) уравнений на
$\theta$-функции механизм интегрирования и весь сопутствующий анализ
становится  прозрачным.

В уравнения входят не величины $\vartheta$, а их логарифмические
производные  и поэтому условиями совместности будут {\em
алгебраические\/} соотношения между функциями $\boldsymbol\theta$,
коэффициентами $\Lambda,\,\ln_{\boldsymbol\tau}\!\vartheta$ и, с
другой стороны, единственное дифференциальное соотношение,
содержащее $\Lambda_{\boldsymbol\tau}\equiv\dot\Lambda$. Мы
получаем:
\begin{equation}\label{Lambda}
\begin{array}{c}\ds
\frac{\dot\vartheta_2}{\vartheta_2}+\Lambda=0\,,
\qquad
\frac{\dot\vartheta_3}{\vartheta_3}+\Lambda=
\frac{\boldsymbol\theta_3^2-\boldsymbol\theta_2^2}
{\boldsymbol\theta_1^2}\,,
\qquad
\frac{\dot\vartheta_4}{\vartheta_4}+\Lambda=
\frac{\boldsymbol\theta_4^2-\boldsymbol\theta_2^2}
{\boldsymbol\theta_1^2}_{{}_{\ds\mathstrut}}\,,
\\ \ds
\dot\Lambda-2\!\left(\frac{\dot\vartheta_3}{\vartheta_3}+
\frac{\dot\vartheta_4}{\vartheta_4}\right)\!\Lambda-2\,
\frac{\dot\vartheta_3}{\vartheta_3}\,
\frac{\dot\vartheta_4}{\vartheta_4}=0\,.
\end{array}
\end{equation}
Отметим что в качестве первичного объекта выступают не симметричные
уравнения \reff{var}, а условия интегрируемости \reff{Lambda}
($\Rightarrow$\reff{intA}). При этом один из параметров
$\ln_{\boldsymbol\tau}\!\vartheta_2$ и $\Lambda$ (или $\eta$)
фактически входил в теорию фиктивно. Поскольку все системы
нелинейны, для исключения функций $\boldsymbol\theta$ из
\reff{Lambda}, необходимо дополнительное дифференцирование.  Для
симметрии мы оставляем три величины $\vartheta$ и, обозначив
$(X,Y,Z)=\bbig{(}\frac{\dot\vartheta_2}{\vartheta_2},
\frac{\dot\vartheta_3}{\vartheta_3},
\frac{\dot\vartheta_4}{\vartheta_4}{\sss\!}\bbig{)}$, получаем
следующие уравнения:
$$
\frac12\,\dot X=(Y+Z)\,X-YZ\,,\quad \frac12\,\dot
Y=(X+Z)\,Y-XZ\,,\quad \frac12\,\dot Z=(X+Y)\,Z-XY\,.
$$
Это широко известная система Дарбу--Альфена \cite[{\bf
I}:\,стр.\,331]{halphen}, \cite{tah}, \cite{zudilin},
\cite[стр.\,577]{conte}, а ее следствием является уравнение
\reff{halphen} из \S\,7:
$$
\big(\dot X-2\,X^2\big)\,\raisebox{0.09em}{\mbox{$\dddot{X}$}} -
\ddot X^2 +16\,X^3\ddot X+4\big(\dot X-6\,X^2\big)\dot X^2=0\,.
$$
Таким образом все уравнения интегрируются в канонических рядах
Якоби.

{\em Замечание\/}. Перенормировку можно было бы продолжить
$\boldsymbol{\theta}\mapsto\boldsymbol{\theta}\,\re_{\mathstrut}^{\frac12\Lambda
x^2}$, тем самым удалив параметр $\Lambda$ из формул. В \reff{Xbold}
тогда надо положить $\Lambda=0$, а условия интегрируемости
\reff{Lambda} превратятся  в простые выражения
$$
Y=X+\pi^2\varkappa^2\vartheta_4^4(\mu)\,,\qquad
Z=X+\pi^2\varkappa^2\vartheta_3^4(\mu)
$$
и уравнение Риккати на функцию $X(\boldsymbol\tau)$ с переменными
коэффициентами $\varkappa(\tau),\,\mu(\tau)$:
$$
\frac12\,\dot X=X^2+\pi^2\varkappa^2
\big\{\vartheta_3^4(\mu)+\vartheta_4^4(\mu)\big\}\!\cdot\!
X-\pi^4\varkappa^4\vartheta_3^4(\mu)\vartheta_4^4(\mu)\,.
$$
Это уравнение, как и предыдущие, тоже интегрируется.

\subsection{Общие интегралы}
Обозначив константы интегрирования  уравнений \reff{intA} как
$(a,b,c,\boldsymbol{d})$ и обозначив
$\boldsymbol{\mathrm{T}}\equiv\frac{a\,\tau+b}{c\,\tau+d}$, мы
получаем общий интеграл  для \reff{intA}:
\begin{equation}\label{T}
\begin{array}{c}
\ds\vartheta_2=
\boldsymbol{d}\,\frac{\vartheta_2(\boldsymbol{\mathrm{T}})}
{\sqrt{c\,\tau+d}}\,,\quad \vartheta_3=
\frac{1}{A_1}\,\frac{\vartheta_3
(\boldsymbol{\mathrm{T}})}{\sqrt{c\,\tau+d}}\,,\quad \vartheta_4=
\frac{1}{A_2}\,
\frac{\vartheta_4(\boldsymbol{\mathrm{T}})}{\sqrt{c\,\tau+d}}\,,\\\\
\ds\eta=\frac{1}{(c\,\tau+d)^2}\left\{
\eta(\boldsymbol{\mathrm{T}})+ \mfrac{\pi^2}{12}
\Big(\mfrac{A_1^4-1}{A_1^4}\,
\vartheta_3^4(\boldsymbol{\mathrm{T}})+ \mfrac{A_2^4-1}{A_2^4}\,
\vartheta_4^4(\boldsymbol{\mathrm{T}})\Big)\!\right\} +
\frac{1}{2}\,\frac{\pi\ri\,c}{c\,\tau+d}\,,
\end{array}
\end{equation}
где под
$\vartheta(\boldsymbol{\mathrm{T}}),\,\eta(\boldsymbol{\mathrm{T}})$
понимаются ряды Якоби из \S\,2 и, как обычно, $a\,d-b\,c=1$.
Заметим, что главный параметр теории
$$
\frac14\,\Lambda=\frac{1}{(c\,\tau+d)^2}\,\Big\{
\eta(\boldsymbol{\mathrm{T}})+ \mfrac{\pi^2}{12}
\bbig(\vartheta_3^4(\boldsymbol{\mathrm{T}})+
\vartheta_4^4(\boldsymbol{\mathrm{T}})\bbig)\Big\}+\frac{1}{2}\,
\frac{\pi\ri\,c}{c\,\tau+d}
$$
не зависит от того, берется канонический случай или нет, и
удовлетворяет уравнению 3-го порядка \reff{halphen}, в которое надо
подставить $X=\frac{\ri}{4\pi}\Lambda$. К сказанному в \S\,9.2
добавим, что симметризация канонического случая уравнений
\reff{intA} в \reff{var}, которую <<напрашивается наложить
определением>>, была бы не вполне корректна. В качестве условий
интегрируемости, уравнения \reff{var} не возникают и  это не зависит
от правил перенормировки и выбора интегралов $A_{1,2}$.
Алгебраический интеграл \reff{A3} является {\em нелинейным\/}   по
переменным $\vartheta^4$, а общий интеграл уравнений \reff{var}
найти пока не удается.

Положим теперь, что величины $\eta,\vartheta$ в уравнениях \reff{X}
есть функции от $\tau$ в силу уравнений \reff{intA}, а константы
$(A,B,C,\varkappa,\mu)$ есть неизвестные функции от $\tau$.
Подставляя тогда формулы \reff{Xgeneral} в \reff{Dtau}, мы получаем
систему дифференциальных уравнений на $(A,B,C,\varkappa,\mu)$.
Вычисления можно заранее упростить, так как у нас уже есть два
интеграла \reff{A12} ($\Leftrightarrow$ \reff{omega}) этих
уравнений:
\begin{equation}\label{intO}
\varkappa\,\frac{\vartheta_3^2(\mu)}{\vartheta_3^2}=A_1^2\,, \qquad
\varkappa\,\frac{\vartheta_4^2(\mu)} {\vartheta_4^2}=A_2^2\,,
\end{equation}
из которых мы можем получить и сами уравнения. Обозначая точкой
производную по $\tau$, указанные уравнения выглядят следующим
образом:
$$
\dot\mu=\varkappa^2\,, \quad
\pi\ri\,\frac{\dot\varkappa}{\varkappa}=2\,M(\varkappa,\mu)\,,\quad
\dot A=0\,,\quad \dot B=A\,\dot\varkappa\,,\quad \frac{\dot
C}{C}=-\frac{\dot \varkappa}{\varkappa}\,.
$$
Первые два уравнения имеют решение
$$
\frac{\vartheta_3(\mu)}{\vartheta_4(\mu)}= \boldsymbol{D}\!\left[
\frac{\vartheta_3(\boldsymbol{\mathrm{T}})}
{\vartheta_4(\boldsymbol{\mathrm{T}})}\right]^{\boldsymbol{E}},\qquad
\varkappa=\frac{\sqrt{\boldsymbol{E}}}{c\,\tau+d}\,
\frac{\vartheta_2^2(\boldsymbol{\mathrm{T}})}{\vartheta_2^2(\mu)}\,,
$$
а оставшиеся, вслед за ними, интегрируются элементарно:
$$
A=\boldsymbol{A}\,,\quad
B=\frac{\boldsymbol{A}\sqrt{\boldsymbol{E}}}{c\,\tau+d}\,
\frac{\vartheta_2^2(\boldsymbol{\mathrm{T}})}{\vartheta_2^2(\mu)}
+\boldsymbol{B}\,,\quad
C=\frac{\boldsymbol{C}}{\sqrt{\boldsymbol{E}}}\,(c\,\tau+d)\,
\frac{\vartheta_2^2(\mu)}{\vartheta_2^2(\boldsymbol{\mathrm{T}})}\,.
$$
В силу уравнений  \reff{intO} и \reff{T}, мы должны положить
$\boldsymbol{D}=\boldsymbol{E}=1$. Такое сокращение числа констант
обусловлено тем, что интегралы \reff{intO} являются  интегралами как
$x$- так и $\tau$-уравнений, а сами уравнения нелинейны. В итоге мы
получаем решение \cite{br}:
$$
\frac{\vartheta_3(\mu)}{\vartheta_4(\mu)}=
\frac{\vartheta_3(\boldsymbol{\mathrm{T}})}
{\vartheta_4(\boldsymbol{\mathrm{T}})}\quad\Rightarrow\quad
\mu=\widehat{\boldsymbol{\Gamma}(4)}(\boldsymbol{\mathrm{T}})\,.
$$
Выбирая для удобства самый простой вариант
$\mu=\boldsymbol{\mathrm{T}}$ (после чего сразу вычисляются
$\varkappa$ и $A,B,C$), мы получаем искомый совместный общий
интеграл уравнений (\ref{X}--\ref{Dtau}), в которых коэффициенты
$\eta,\vartheta$ определяются по формулам \reff{T}:
$$
\begin{array}{lr}
\ds \pm\,\theta_1=&\ds\!\!\! \frac{1}{A_1A_2}\,
\frac{\boldsymbol{d}\,\boldsymbol{C}}{\sqrt{c\,\tau+d}}\,
\theta_1\!\Big(\mfrac{x+\boldsymbol{A}}{c\,\tau+d}+\boldsymbol{B}
\Big|\mfrac{a\,\tau+b}{c\,\tau+d}
\Big)\,\re_{\ds\mathstrut}^{\frac{-\pi\ri\, c}
{c\tau+d}(x+\boldsymbol{A})^{2^{\ds\mathstrut}}}\,,\\
\ds \pm\,\theta_2=&\ds\!\!\!
\frac{\boldsymbol{d}\,\boldsymbol{C}}{\sqrt{c\,\tau+d}}\,
\theta_2\!\Big(\mfrac{x+\boldsymbol{A}}{c\,\tau+d}+\boldsymbol{B}
\Big|\mfrac{a\,\tau+b}{c\,\tau+d}
\Big)\,\re_{\ds\mathstrut}^{\frac{-\pi\ri\, c}
{c\tau+d}(x+\boldsymbol{A})^{2^{\ds\mathstrut}}}\,,
\\\pm\,\theta_3=&
\ds\!\!\! \frac{1}{A_1}\,
\frac{\phantom{\boldsymbol{d}\,}\boldsymbol{C}}{\sqrt{c\,\tau+d}}\,
\theta_3\!\Big(\mfrac{x+\boldsymbol{A}}{c\,\tau+d}+\boldsymbol{B}
\Big|\mfrac{a\,\tau+b}{c\,\tau+d}
\Big)\,\re_{\ds\mathstrut}^{\frac{-\pi\ri\, c}
{c\tau+d}(x+\boldsymbol{A})^{2^{\ds\mathstrut}}}\,,
\\
\ds \pm\,\theta_4=&\ds\!\!\! \frac{1}{A_2}\,
\frac{\phantom{\boldsymbol{d}\,}\boldsymbol{C}}{\sqrt{c\,\tau+d}}\,
\theta_4\!\Big(\mfrac{x+\boldsymbol{A}}{c\,\tau+d}+\boldsymbol{B}
\Big|\mfrac{a\,\tau+b}{c\,\tau+d}
\Big)\,\re_{\ds\mathstrut}^{\frac{-\pi\ri\, c}
{c\tau+d}(x+\boldsymbol{A})^{2^{\ds\mathstrut}}}\,.
\end{array}
$$
Формула для $\theta_1'$ получается взятием производной от формулы
для $\theta_1$. Для редукции к каноническому случаю \reff{sol}
достаточно положить $A_1=A_2=\boldsymbol{d}=1$ и выбрать
преобразование $\left(\begin{smallmatrix}a&b\\c&d
\end{smallmatrix} \right)=\left(\begin{smallmatrix}1&0\\2&1
\end{smallmatrix} \right)$, так как группа ${\boldsymbol{\Gamma}(2)}$ не
переставляет функции $\theta$ или $\vartheta$.

\section{Приложение. Уравнение Пенлеве-VI}
\noindent Связь эллиптических функций с шестым уравнением Пенлеве
\begin{equation}\label{P6}
\begin{array}{l}\ds
y_{\mathit{xx}}^{}=\frac12\! \left(\frac1y+\frac{1}{y-1}+
\frac{1}{y-x}\right) y_x^2- \left(\frac1x+\frac{1}{x-1}+
\frac{1}{y-x}\right)y_x^{}+{}\\\\
\ds \phantom{y_{\mathit{xx}}^{}=}{} +\frac{y(y-1)(y-x)}{x^2(x-1)^2}
\left\{\alpha-\beta\,\frac{x}{y^2}+\gamma\,\frac{x-1}{(y-1)^2}-
\bbig(\delta-\mfrac12\bbig)\,\frac{x(x-1)}{(y-x)^2} \right\}
\end{array}
\end{equation}
(как впрочем появление и самого уравнения \reff{P6}) была
установлена Р.\,Фуксом в работе \cite{fuchs}. Вскоре после него
Пенлеве \cite{painleve} придал этому уравнению замечательный вид
(конвертируя в используемые здесь обозначения)
\begin{equation}\label{P6wp}
-\frac{\pi^2}{4}\,\frac{d^2 z}{d\tau^2}=\alpha\,\wpp(z|\tau)+
\beta\,\wpp(z-1|\tau)+\gamma\,\wpp(z-\tau|\tau)+
\delta\,\wpp(z-1-\tau|\tau)\,,
\end{equation}
сделав трансцендентную замену переменных $(y,x)\mapsto(z,\tau)$
\begin{equation}\label{subs}
x=\frac{\vartheta_4^4(\tau)}{\vartheta_3^4(\tau)}\,,\qquad
y=\frac13+\frac13\frac{\vartheta_4^4(\tau)}{\vartheta_3^4(\tau)}
-\frac{4}{\pi^2} \frac{\wp(z|\tau)}{\vartheta_3^4(\tau)}\,.
\end{equation}
К этим  известным фактам необходимо добавить следующий комментарий.

Уравнение \reff{P6} и подстановка \reff{subs} отражают главные
свойства как самого уравнения так и его решений. Если уравнение
имеет второй порядок, линейно по $y_{\mathit{xx}}^{}$, рационально
по $y_x^{}$ и имеет только неподвижные точки ветвления в решениях
(свойство Пенлеве), тогда, как мы видим,  число таких точек не
превосходит три. Их всегда можно поместить в $x_j=\{0,1,\infty\}$, а
уравнение будет иметь вид уравнения Фукса--Пенлеве \reff{P6} либо
некоторый его предельный случай \cite{painleve}. Далее, плоскость
$(x)$  можно конформно и взаимно-однозначно отобразить на
фундаментальный 4-угольник группы $\boldsymbol{\Gamma}(2)$ в
плоскости новой переменной $(\tau)$ с помощью модулярной функции
$x=k'^2(\tau)$ (первая формула в \reff{subs}). В свою очередь,
поведение функции $k'^2(\tau)$ в прообразах точек $x_j$ имеет
экспоненциальный характер по локальному параметру $\tau$ и поэтому
ветвление произвольного степенного  или логарифмического характера в
окрестности $x_j$ перейдет в локально однозначную зависимость по
переменной $\tau$. В этом мы видим частичное объяснение
происхождения подстановки \reff{subs}. Через модулярную функцию она
появилась у Пенлеве \cite{painleve}, а через уравнение Лежандра ---
у Фукса \cite{fuchs}. Остальные особые точки решений являются
полюсами, могут быть подвижными, но их расположение и глобальное
поведение решений полностью определяется логарифмическими
производными от {\em целых\/} функций \cite{painleve},
\cite[стр.\,77--180]{conte}. Построение таких целых трансцендентных
функций, как это подчеркивалось самим Пенлеве (1902), завершает
процедуру интегрирования уравнения (<<int\'egration parfaite>> в его
терминологии). Для уравнения \reff{P6} характер полюсов известен
(Пенлеве), а его решение имеет структуру
\begin{equation}\label{entire}
y\sim \frac{d}{dx}\mathrm{Ln}\,
\frac{\mbox{\Large$\boldsymbol\tau$}_{\!\!2}^{}}
{\mbox{\Large$\boldsymbol\tau$}_{\!\!1}^{}}
\end{equation}
с целыми функциями $\mbox{\Large$\boldsymbol\tau$}_{\!\!1,2}^{}(x)$,
имеющими, быть может, неподвижные критические особенности. Полный
список таких формул см. в \cite[стр.\,165]{conte}.

Известно только два случая\footnote{С учетом  автоморфизмов в
пространстве параметров, которые хорошо известны в литературе по
уравнениям Пенлеве. См. например статью Громака в \cite{conte},
работу \cite{korotkin} и ссылки в них.}, когда для уравнения
\reff{P6} выписывается его общий интеграл. Это решения Пикара (1889)
и Хитчина (1995), но ни для одного из этих случаев форма Пенлеве
\reff{entire} не известна. В этой связи представляет интерес
продемонстрировать нетривиальную ситуацию,  когда функции
$\mbox{\Large$\boldsymbol\tau$}_{\!\!1,2}^{}(x)$ выписываются явно.
Ниже, мы приводим такой пример. Он соответствует параметрам
$\big(\alpha=\beta=\gamma=\delta=\frac18\big)$, которые найдены
Хитчиным при описании SU(2)-инвариантных анти-самодуальных метрик
уравнений Эйнштейна \cite{hitchin}. Уравнение \reff{P6wp} тогда
приводится к виду
$$
-\pi^2\,\frac{d^2 z}{d\tau^2}=4\,\wpp(2z|\tau) \quad
\Leftrightarrow\quad\frac{d^2 z}{d\tau^2}=4\pi\,\ded^9(\tau)\,
\frac{\theta_1^{}(2z|\tau)}{\theta_1^4(z|\tau)}\,,
$$
а  его общий интеграл, для параметрической формы \reff{subs}, найден
в \cite[стр.\,74, 78]{hitchin}:
$$
\wp(z|\tau)=\wp(A\tau+B|\tau)+\frac12\, \frac{\wpp(A\tau+B|\tau)}
{\zeta(A\tau+B|\tau)-(A\tau+B)\,\eta(\tau)+\frac{\pi\ri}{2}\, A}\,.
$$

\subsection{Форма Пенлеве}
Описанное выше дифференциальное исчисление вейерштрассовских и
якобиевских функций фактически автоматизирует  вычисления, связанные
с любыми  <<эллиптическими>> решениями уравнений
(\ref{P6}--\ref{subs}) и их вырождениями.   Сложность решения,
однако, отмечалась  как в первоначальной работе
\cite[стр.\,75]{hitchin}, так и позже (см. замечание 5.1 в
\cite{korotkin})\footnote{В $\theta$-функциональных решениях работ
\cite{hitchin,korotkin} присутствуют функции
$\theta_1,\theta_1',\theta_1'',\theta_1''',\vartheta_1',\vartheta_1'',
\vartheta_1'''$.}. Используя результаты \S\,5 рассматриваемому
решению можно придать следующий вид:
$$
\wp(z|\tau)=\frac{\pi}{2\,\ri}\,\frac{d}{d\tau}\mathrm{Ln}
\frac{\zeta(A\tau+B|\tau)-(A\tau+B)\,\eta(\tau)+\frac{\pi\ri}{2}\,A}
{\ded^2(\tau)}\,.
$$
Переводя теперь все это в $\theta$-функции, получаем параметрическую
форму решения:
\begin{equation}\label{Tau}
y=\frac{2\,\ri}{\pi}\frac{1}{\vartheta_3^4(\tau)}\, \frac{d}{d\tau}
\mathrm{Ln}\frac{\theta_{1\!\!\!}'(A\tau+B|\tau)+2\pi\ri\,A\,
\theta_1^{}\!(A\tau+B|\tau)}
{\vartheta_2^2(\tau)\,\theta_1^{}\!(A\tau+B|\tau)}\,, \qquad
x=\frac{\vartheta_4^4(\tau)}{\vartheta_3^4(\tau)}\,.
\end{equation}
Остается переписать это решение в исходных переменных $(x,y)$
используя формулы обратного перехода $\tau\mapsto x$:
\begin{eqnarray}
\ds
\frac{d}{d\tau}&\!\!\!\!=&\!\!\!\pi\ri\,x\,(x-1)\,\vartheta_3^4(\tau)\,
\frac{d}{dx}\,,\qquad
\vartheta_2^2(\tau)=\frac{2}{\pi}\,\sqrt{1-x}\,K'(\sqrt{x})\,,
\nonumber\\\nonumber\\
\ds
\re^{\pi\ri\tau}&\!\!\!\!=&\!\!\!\exp\Big\{\!-\!\pi\mfrac{K(\sqrt{x})}
{K'(\sqrt{x})}\Big\}=\label{weier}
\\\nonumber\\\ds
\phantom{e^{\pi\ri\,\tau}} &\!\!\!\!=&\!\!\!
\Big(\mfrac{1-x}{16}\Big)+8 \Big(\mfrac{1-x}{16}\Big)^{\!2}
+84\Big(\mfrac{1-x}{16}\Big)^{\!3} +
992\Big(\mfrac{1-x}{16}\Big)^{\!4}
+12514\Big(\mfrac{1-x}{16}\Big)^{\!5}+\cdots\,. \nonumber
\end{eqnarray}
Упростив, и меняя $\ri A\mapsto  A$, получим окончательный ответ
\begin{equation}\label{hitchin}
y=x\,(1-x)\,\frac{d}{dx}\mathrm{Ln}
\mfrac{\left[\mbox{\normalsize$\theta_1'$}\!\!\left(\!
A\frac{K(\sqrt{x})}{K'(\sqrt{x})}+B\bbig|\frac{\ri\,
K(\sqrt{x})}{K'(\sqrt{x})}\right)+ \mbox{\normalsize$2\pi
A\!\cdot\!\theta_1$}\!\!\left(\!
A\frac{K(\sqrt{x})}{K'(\sqrt{x})}+B\bbig|\frac{\ri\,
K(\sqrt{x})}{K'(\sqrt{x})}\right) \right]^2}
{\mbox{\normalsize$(1-x)\,\theta_1^2$} \!\left(\!
A\frac{K(\sqrt{x})}{K'(\sqrt{x})}+B\bbig|\frac{\ri\,
K(\sqrt{x})}{K'(\sqrt{x})}\right)
\mbox{\normalsize$K'^2(\sqrt{x})$}}\,,
\end{equation}
в котором $K$ и $K'$ берутся, в зависимости от предпочтений, как
классические полные эллиптические интегралы или, например,
гипергеометрические функции \cite{a,bateman,WW}:
$$
K(\sqrt{x})=\frac{\pi}{2}\cdot
{}_2F_1\Big(\mfrac12,\mfrac12;1\bbig|x \Big)\,,\qquad
K'(\sqrt{x})=K(\sqrt{1-x})\,.
$$
Хорошее упражнение
--- проверить решение \reff{hitchin} прямой подстановкой. Для этого
достаточно дополнить изложенный выше аппарат известными правилами
дифференцирования  полных эллиптических интегралов $K,K'$ и $E,E'$
\cite{a}:
\begin{equation}\label{EK}
\begin{array}{rlrl}
\ds
2\,\frac{d}{dx}K(\sqrt{x})&\!\!\!\ds=\frac{E(\sqrt{x})}{x\,(1-x)}-
\frac{K(\sqrt{x})}{x}\,,&\ds
2\,\frac{d}{dx}K'(\sqrt{x})&\!\!\!\ds=\frac{E'(\sqrt{x})}{x\,(x-1)}+
\frac{K'(\sqrt{x})}{x-1}\,,\\\\
\ds 2\,\frac{d}{dx}E(\sqrt{x})&\!\!\!\ds=\frac{E(\sqrt{x})}{x}-
\frac{K(\sqrt{x})}{x}\,, &\ds \quad
2\,\frac{d}{dx}E'(\sqrt{x})&\!\!\!\ds=\frac{E'(\sqrt{x})}{x-1}+
\frac{K'(\sqrt{x})}{x-1}
\end{array}
\end{equation}
и соотношением Лежандра $\etap=\tau\,\eta-\frac{\pi}{2}\ri$,
записанным в <<$x$-представлении>>:
$$
EK'+E'K-KK'=\frac{\pi}{2}\,.
$$
Строго говоря, в  \reff{hitchin} следует еще сократить общий <<не
целый>> множитель $\exp(\frac{\pi\ri}{4}\tau)$, присутствующий в
рядах $\theta_1^{}$ и $\theta_1'$. Он, впрочем, является
<<неподвижным>> как и немероморфная особенность решения,
определяемая множителем $K'^2(\sqrt{x})$. Из \reff{hitchin} и
\reff{EK} не трудно видеть, что  решение есть простая сумма
<<неподвижной критической особенности>> логарифмического типа и
накапливающихся <<подвижных полюсов>>:
\begin{equation}\label{52}
y=\frac{E'(\sqrt{x})}{K'(\sqrt{x})}+2\,x\,(1-x)\,\frac{d}{dx}\mbox{Ln}
\left\{\frac{\theta_1'}{\theta_1}\!\!\left(\! \textstyle
A\frac{K(\sqrt{x})}{K'(\sqrt{x})}+B\bbig|\frac{\ri\,
K(\sqrt{x})}{K'(\sqrt{x})}\right)+2\pi A\right\},
\end{equation}
где первое слагаемое имеет конечную критическую особенность  только
в точке $x=0$.

Ряды, аналогичные рядам типа \reff{weier}, Вейерштрасс выписывал в
связи с рассмотрением эллиптической модулярной задачи обращения для
модулярной функции $k^2(\tau)$ \cite[стр.\,53--4, 56, 58]{we2},
\cite[стр.\,367]{WW}. Подобным способом выписываются разложения в
окрестности любых других точек, а сами формулы \reff{Tau} или
\reff{hitchin}, \reff{52} доставляют аналитический ответ к степенным
разложениям работы \cite{hitchin} на стр.\,89--92, 108--9.

\subsection{Распределение полюсов}
На {\sc Рис.\,1} частично представлено  трансцендентное
распределение  полюсов $x_{\mathit{mn}}$, принадлежащих одной из
двух серий. А именно, серии, которая допускает явную параметризацию
своих решений:
\begin{equation}\label{poles}
\theta_1 \!\!\left(\!
A\mfrac{K(\sqrt{x})}{K'(\sqrt{x})}+B\Big|\mfrac{\ri\,
K(\sqrt{x})}{K'(\sqrt{x})}\right)=0\quad\Rightarrow\quad
x_{\mathit{mn}}=
\frac{\vartheta_4^4}{\vartheta_3^4}\Big(\mfrac{m-B}{n+A}\Big)\,,
\quad n,m\in\mathbb{Z}\,.
\end{equation}
Сюда необходимо добавить естественное  условие на $(n,m)$ вида
$\boldsymbol\Im\bbig(\frac{m-B}{n+A} \bbig)>0$. Данные распределения
начинают деформироваться при изменении начальных данных $(A,B)$,
порождая разнообразные картины, но все они имеют точки накопления
полюсов в неподвижных  особенностях $x_j=\{0,1,\infty \}$.

\begin{figure}[htbp]
\centering
%
%
\raisebox{16em}{\mbox{\hspace{20em}$(x)$}}\hspace{-21em}%
\includegraphics[width=10 cm,angle=0]{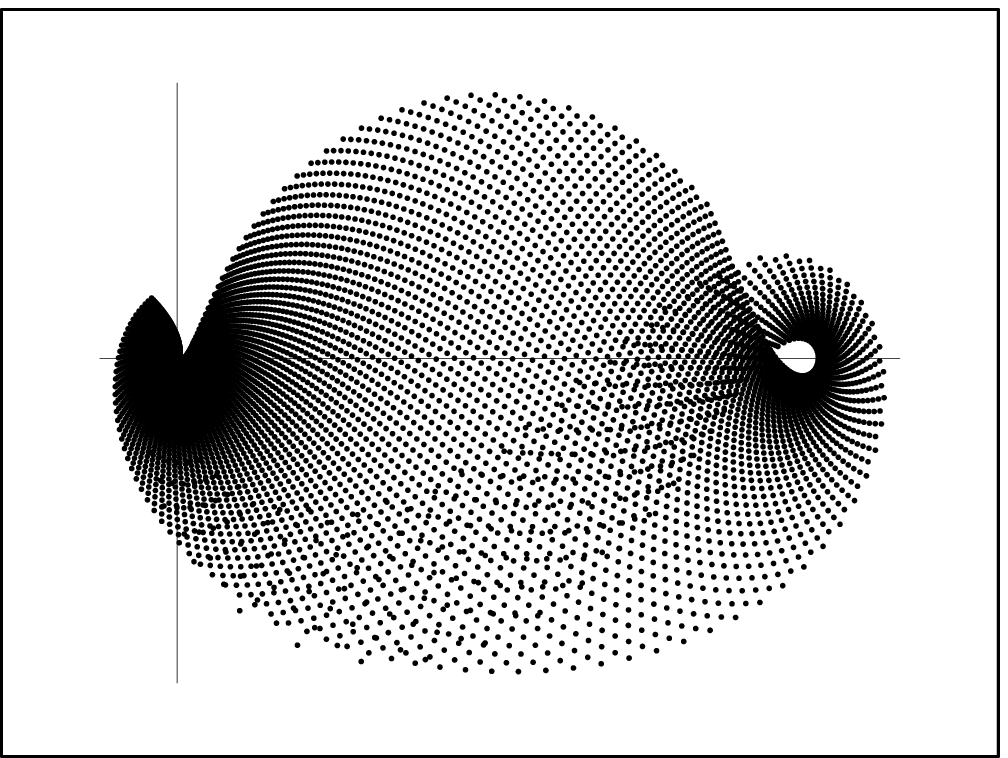}
\caption{Полюса $x_{\mathit{mn}}$ вида \reff{poles} решений
Пикара--Хитчина  \reff{hitchin}, \reff{pic} при
$A=125.45-103.29\,\ri$, $B=36.710-69.980\,\ri$ и $(n,m)=-30\ldots
70$.}
\end{figure}

Зная строение функции $k'^2(\tau)$ и фундаментальной области группы
ее автоморфизмов, т.\,е. $\boldsymbol{\Gamma}(2)$, не трудно также
описать вещественные или чисто комплексные полюса из серии
\reff{poles}. Заметим, что эта же серия полюсов \reff{poles}
исчерпывает все полюса решения уравнения \reff{P6} в случае Пикара:
$\alpha=\beta=\gamma=\delta=0$. В самом деле, в силу \reff{P6wp}
$\ddot z=0\;\Rightarrow\;z=a\tau+b$ и, конвертируя подстановку
\reff{subs} в $\theta$-функции, мы получим
\begin{equation}\label{pic}
y=-\frac{\vartheta_3^2(\tau)}{\vartheta_4^2(\tau)}\,
\frac{\theta_2^2\big(\frac z2|\tau \big)}{\theta_1^2\big(\frac
z2|\tau \big)}\qquad\Rightarrow\qquad
y_{\mbox{\tiny\sc Pic}}^{}=-\frac{1}{\sqrt{x}}\,
\mfrac{\mbox{\normalsize$\theta_2^2$}\left(\!
A\frac{K(\sqrt{x})}{K'(\sqrt{x})}+B\bbig|\frac{\ri\,
K(\sqrt{x})}{K'(\sqrt{x})}\right)}{\mbox{\normalsize$\theta_1^2$}
\left(\! A\frac{K(\sqrt{x})}{K'(\sqrt{x})}+B\bbig|\frac{\ri\,
K(\sqrt{x})}{K'(\sqrt{x})}\right)}\,.
\end{equation}

Решения Пикара--Хитчина  есть функции от $\sqrt{x}$. С другой
стороны,  сама величина $\sqrt{x}=s$, через интегралы
$K,K'(\sqrt{x})$ и подстановку \reff{subs}, связана с линейным
дифференциальным уравнением Фукса 2-го порядка с четырьмя
особенностями $s=\{0,1,-1,\infty \}$. Следовательно, это будет
частный случай уравнения Хойна \cite{br}:
\begin{equation}\label{heun}
Y_{ss}=-\frac12\,\frac{(s^2+1)^2}{s^2\,(s^2-1)^2}\,Y\quad\Rightarrow\quad
\frac{Y_2}{Y_1}=\ri\,\frac{K(s)}{K'(s)}=\tau\,.
\end{equation}
Поскольку функции $\zeta,\wp$  определяются неопределенными
эллиптическими интегралами, мы получаем отсюда, что {\em
рассматриваемые случаи уравнения Пенлеве являются интегрируемыми в
эллиптических квадратурах над дифференциальным полем, определяемым
функциями Хойна вида \reff{heun} или, эквивалентно,
гипергеометрическими функциями Лежандра от квадратичного расширения
рациональных функций по $x$\/}.

Мы наблюдаем также, что известное преобразование Окамото $y\mapsto
\tilde y$ \cite{korotkin}
$$
\tilde
y=\frac{x\,(x-1)\,y\,y_x^{}-x\,y\,(y-1)}{x\,(x-1)\,y_x^{}-y\,(y-1)}\,,
$$
переводящее решения Пикара и Хитчина друг в друга (известный факт),
оставляет инвариантной полюса \reff{poles}, но
добавляет/унич\-то\-жа\-ет вторую серию полюсов, определяемую нулями
уравнения
$$
\zeta(A\tau+B|\tau)=A\,\zeta(\tau|\tau)+B\,\zeta(1|\tau)\,.
$$
Это эквивалентно нахождению  $A$-точек трансцендентной аналитической
функции
\begin{equation}\label{mero}
f(x;A,B):\qquad
\frac{1}{2\pi}\mfrac{\mbox{\normalsize$\theta_1'$}\left(\!
A\frac{K(\sqrt{x})}{K'(\sqrt{x})}+B\bbig|\frac{\ri\,
K(\sqrt{x})}{K'(\sqrt{x})}\right)}{\mbox{\normalsize$\theta_1$}
\left(\! A\frac{K(\sqrt{x})}{K'(\sqrt{x})}+B\bbig|\frac{\ri\,
K(\sqrt{x})}{K'(\sqrt{x})}\right)}\,,
\end{equation}
т.\,е. уравнению $f(x;A,B)=A$. Функция \reff{mero} есть не что иное
как канонический нормализованный мероморфный эллиптический интеграл
$I(z+B|\tau)$, имеющий полюс в точке $z=-B$ и рассматриваемый на
прямолинейном сечении $z=A\tau$ пространства, получаемого как прямое
произведение эллиптической кривой $(z)$ и 1-мерного пространства
модулей эллиптических кривых $(\tau)$.

\mbox{\Large$\boldsymbol\tau$}-функции являются важными объектами в
теории уравнений Пенлеве,  поскольку генерируют их гамильтонианы и
другие объекты \cite{conte,korotkin}. Они, однако, необходимы не
столько для гамильтонианов, которые известны для всех уравнений
Пенлеве, сколько для представления решений через них в форме Пенлеве
\reff{entire}. См. например пояснения Окамото на стр.\,740--741 в
\cite{conte} в связи с построением
\mbox{\Large$\boldsymbol\tau$}-функций, соответствующих
положительной и отрицательной серии вычетов для решений уравнения
Пенлеве-II. \mbox{\Large$\boldsymbol\tau$}-функция, построенная в
\cite{korotkin}, соответствует только пикаровской серии полюсов
\reff{poles} и совпадает,  с точностью до <<неподвижного
немероморфного>> множителя, с нашей функцией
$$
\mbox{\Large$\boldsymbol\tau$}_{\!\!1}^{}(x;A,B)\sim\theta_1
\!\!\left(\! A\mfrac{K(\sqrt{x})}{K'(\sqrt{x})}+B\Big|\mfrac{\ri\,
K(\sqrt{x})}{K'(\sqrt{x})}\right) \,.
$$
Это означает, что существует второй гамильтониан, соответствующий
второй $\mbox{\Large$\boldsymbol\tau$}$-функ\-ции и полюсам с
противоположным знаком вычетов. Она дается формулой
$$
\mbox{\Large$\boldsymbol\tau$}_{\!\!2}^{}(x;A,B)=
\mbox{\Large$\boldsymbol\tau$}_{\!\!1}^{}(x;A,B)\,\frac{d}{dB}
\mbox{Ln}\,\big\{
\mbox{\Large$\boldsymbol\tau$}_{\!\!1}^{}(x;A,B)\,\re^{2\pi A
B}_{\mathstrut}\big\}\,.
$$
Выбор такого представления объясняется  тем, что мероморфные абелевы
интегралы (интегралы 2-го рода) типа функций \reff{mero} могут быть
представлены как производные от логарифмических интегралов
(интегралы 3-го рода) по параметру, определяющему расположение
одного из двух логарифмических полюсов (константа $B$). Этот факт
может быть предметом дальнейших обобщений, но мы их здесь не
обсуждаем.

За недостатком места, мы не приводим также другие формулы, которые
здесь уместны, поскольку они с очевидностью выводятся из
аналитических ответов как в $x$-пред\-став\-ле\-нии \reff{hitchin},
\reff{52}, так и в $\tau$-представлении \reff{Tau}. Например,
простые выражения для космологических метрик Тода--Хитчина
\cite{tod,babich,hitchin} или явный вид дифференциальных уравнений
3-го порядка на {\Large$\boldsymbol\tau$}-функции. На существование
таких уравнений неоднократно указывал Пенлеве (см. например
\cite[стр.\,1114]{painleve}). Они являются следствиями уравнения
3-го порядка на $\sigma(x|\tau)$, упомянутого в \S\,5.

\begin{figure}[htbp]
\centering
\includegraphics[scale=0.89]{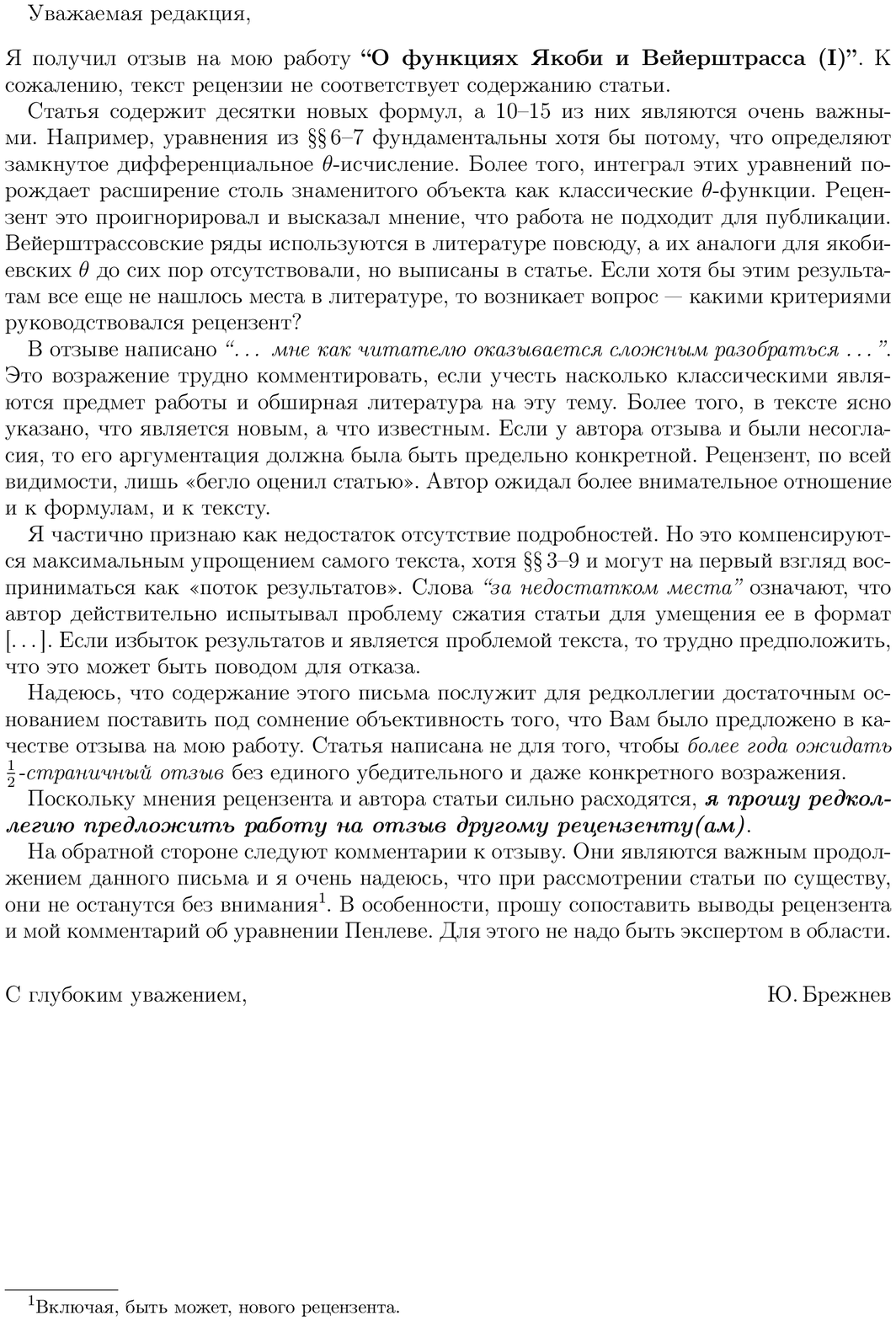}
\end{figure}
\begin{figure}[htbp]
\centering
\includegraphics[scale=0.89]{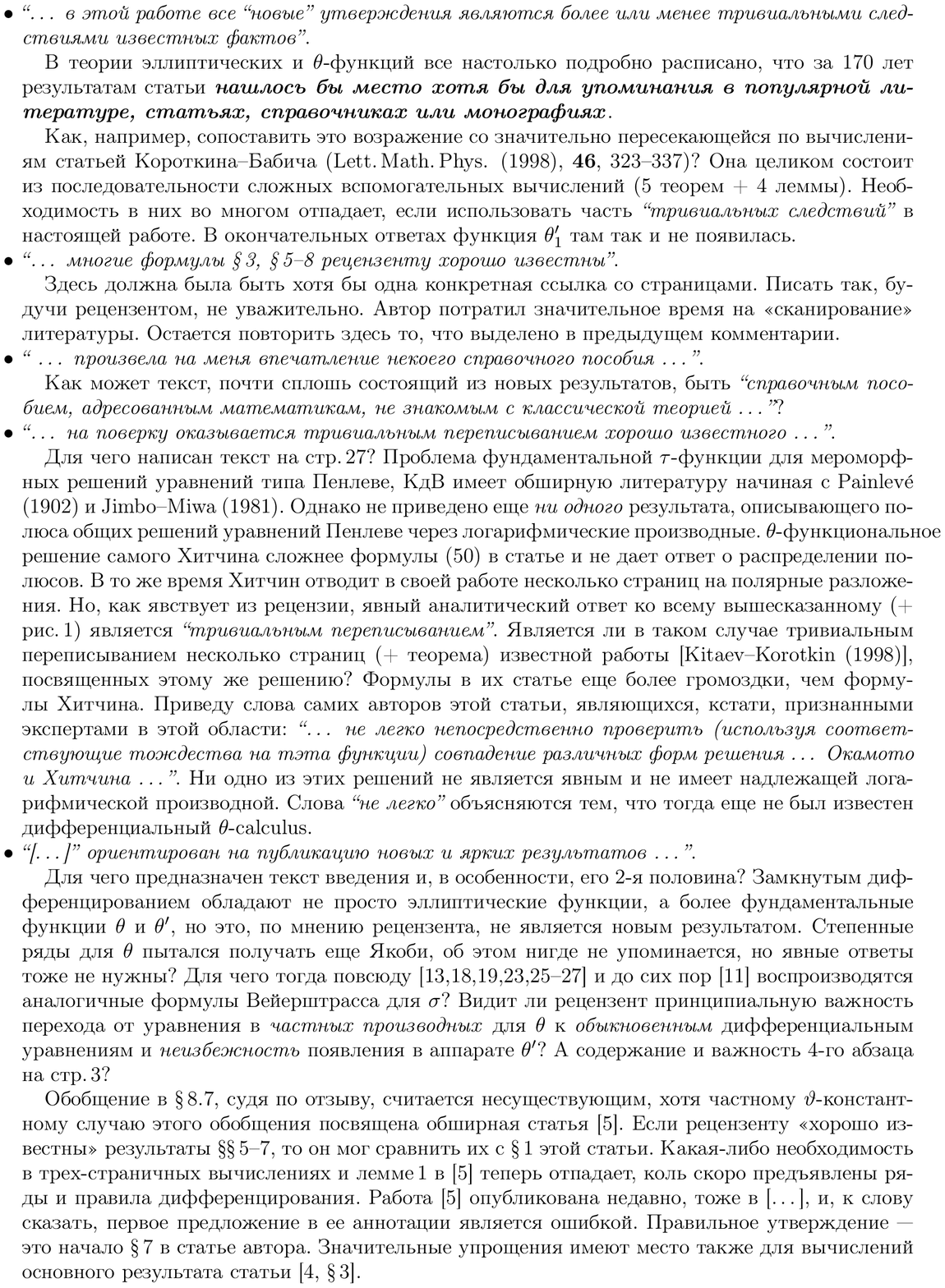}
\end{figure}


\thebibliography{99}

\bibitem{hitchin2}\mbox{\sc Атья, М., Хитчин, Н.}
{\em Геометрия и динамика магнитных монополей\/}. Москва (1991).

\bibitem{a}\mbox{\sc Ахиезер, Н.\,И.}
{\em Элементы теории эллиптических функций\/}. Москва (1970).

\bibitem{bateman}\mbox{\sc Бейтмен, Г., Эрдейи, А.}
{\em Высшие трансцендентные функции. Эллиптические и автоморфные
функции. Функции Ламе и Матье\/}. Москва (1967).

\bibitem{pavlov}\mbox{\sc Бухштабер, В.\,М., Лейкин, Д.\,В.,
Павлов, М.\,В.} {\em Егоровские гидродинамические цепочки, уравнение
Шази и операторы, аннулирующие сигма-функцию\/}. Функц.\,Анализ и
его Прил. (2003), {\bf 37}, 13--26.

\bibitem{zudilin}\mbox{\sc Зудилин, В.\,В.}
{\em Тэта-константы и дифференциальные уравнения\/}. Мат.\,Сборник
(2000), {\bf 191}(12), 77--122.

\bibitem{WW}\mbox{\sc Уиттекер, Э.\,Т., Ватсон, Дж.\,Н.}
{\em Курс современного анализа\/}.  {\bf II}. Москва (1963).

\bibitem{apostol}\mbox{\sc Apostol, N.\,M.} {\em Modular functions and
Dirichlet Series in Number Theory\/}. Springer--Verlag (1976).

\bibitem{babich}\mbox{\sc Babich, M.\,V., Korоtkin, D.\,A.}
{\em Self-dual SU(2) invariant Einstein metrics and modular
dependence of theta-functions\/}. Lett.\,Math.\,Phys. (1998), {\bf
46}, 323--337.

\bibitem{br}\mbox{\sc Brezhnev, Yu.\,V.}
{\em On uniformization of algebraic curves\/}. Moscow Math.\,Journ.
(2008), {\bf 8}(2), 1--39.

\bibitem{conte}\mbox{\sc Conte, R.} (Ed.)
{\em The Painlev\'e property. One century later\/}. CRM Series in
Mathematical Physics (1999).

\bibitem{eilbeck}\mbox{\sc Eilbeck, J.\,C., Enol'skii, V.\,Z.}
{\em Bilinear operators and power series for the Weierstrass
$\sigma$-function\/}. Journ.\,Phys. {\bf A}: Math.\,Gen. (2000),
{\bf 33}, 791--794.


\bibitem{fuchs}\mbox{\sc Fuchs, R.} {\em Sur quelques \'equations
diff\'erentielles lin\'eares du second ordre\/}.
Compt.\,Rend.\,Acad.\,Sci. (1905), {\bf CXLI}(14), 555--558.

\bibitem{halphen}\mbox{\sc Halphen, G.-H.}
{\em Trait\'e des Fonctions Elliptiques et de Leurs Applications\/}.
{\bf I--III}. Gauthier--Villars: Paris (1886--1891).

\bibitem{hermite}\mbox{\sc Hermite, C.} {\em \OE uvres}. {\bf II}.
Gauthier--Villars: Paris (1908).

\bibitem{hitchin}\mbox{\sc Hitchin, N.} {\em Twistor spaces,
Einstein metrics and isomonodromic deformations.}
Journ.\,Diff.\,Geom. (1995), {\bf 42}(1), 30--112.

\bibitem{jacobi}\mbox{\sc Jacobi, C.} {\em Gesammelte Werke\/}.
{\bf I, II}. Verlag von G.\,Reimer: Berlin (1882--1891).

\bibitem{korotkin}\mbox{\sc Kitaev, A.\,V., Korotkin, D.\,A.}{\em On solutions of the
Schlesinger equations in terms of theta-functions\/}. Intern. Math.
Research Notices (1998), {\bf 17}, 877--905.

\bibitem{koenig}\mbox{\sc Koenigsberger, L.} {\em
Vorlesungen \"uber die theorie der Elliptischen Functionen\/}. {\bf
II}. Leipzig: Druck und Verlag von B.\,G.\,Teubner (1874).


\bibitem{krause}\mbox{\sc Krause, M.}
{\em Theorie der Doppeltperiodischen Functionen\/}. {\bf I}. Verlag
von B.\,G.\,Teubner: Leipzig (1895).

\bibitem{krazer}\mbox{\sc Krazer, A.}
{\em Lehrbuch der Thetafunktionen\/}. Leipzig: Teubner (1903).
Reprint: New York, Chelsea (1970).


\bibitem{painleve}\mbox{\sc Painlev\'e, P.} {\em Sur les \'equations
diff\'erentialles du second ordre \`a points critiques fixes\/}.
Compt.\,Rend.\,Acad.\,Sci. (1906), {\bf CXLIII}(26), 1111--1117.

\bibitem{tah}\mbox{\sc Takhtajan, L.\,A.} {\em A simple example of
modular forms as tau-functions for integrable equations\/}.
Теор.\,Мат.\,Физика (1992), {\bf 93}(2), 330--341.

\bibitem{tannery}\mbox{\sc Tannery, J., Molk, J.}
{\em Elements de la theorie des fonctions elliptiques\/}. {\bf
I--IV}. Gauthier--Villars: Paris (1893--1902).


\bibitem{tod}\mbox{\sc Tod, K.\,P.} {\em Self-dual Einstein metrics
from the Painlev\'e VI equation\/}. Phys.\,Lett.\,{\bf A} (1994),
{\bf 190}, 221--224.

\bibitem{weber}\mbox{\sc Weber, H.} {\em Lehrbuch der Algebra. III.
Elliptische Funktionen und algebraische Zahlen\/}. F.\,Vieweg \&
Sohn: Braunschweig (1908).

\bibitem{we}\mbox{\sc Weierstrass, K.} {\em Mathematische Werke\/}.
{\bf II, V}. Mayer \& M\"uller: Berlin (1894).

\bibitem{we2}\mbox{\sc Weierstrass, K.}
{\em Formeln und Lehrs\"atze zum Gebrauche der elliptischen
Functionen\/}. Bearbeitet und herausgegeben von H.\,A.\,Schwarz.
G\"ottingen (1885).
\end{document}